\documentclass[12pt,twoside]{article}

\usepackage{graphicx}
\usepackage{amsfonts}

\def\o{\overline}
\def\v{\vskip.5ex}
\def\vv{\vskip1ex}
\def\vvv{\vskip2ex}
\def\vvvv{\vskip4ex}

\def\picill#1by#2(#3)
{\vbox to #2
{\hrule width #1 height 0pt depth 0pt
\vfill\epsffile{#3}}}

\makeatletter
\long\def\UR#1{\leavevmode\setbox\@tempboxa\hbox{#1}\@tempdima\fboxrule
    \advance\@tempdima \fboxsep \advance\@tempdima \dp\@tempboxa
   \hbox{\lower \@tempdima\hbox
  {\vbox{\hrule \@height \fboxrule
          \hbox{  \hskip\fboxsep
          \vbox{\vskip\fboxsep \box\@tempboxa\vskip\fboxsep}\hskip
                 \fboxsep\vrule \@width \fboxrule}%
                  }}}}
\makeatother

\makeatletter
\long\def\LR#1{\leavevmode\setbox\@tempboxa\hbox{#1}\@tempdima\fboxrule
    \advance\@tempdima \fboxsep \advance\@tempdima \dp\@tempboxa
   \hbox{\lower \@tempdima\hbox
  {\vbox{
          \hbox{  \hskip\fboxsep
          \vbox{\vskip\fboxsep \box\@tempboxa\vskip\fboxsep}\hskip
                 \fboxsep\vrule \@width \fboxrule}%
                 \hrule \@height \fboxrule}}}}
\makeatother

\makeatletter
\long\def\UL#1{\leavevmode\setbox\@tempboxa\hbox{#1}\@tempdima\fboxrule
    \advance\@tempdima \fboxsep \advance\@tempdima \dp\@tempboxa
   \hbox{\lower \@tempdima\hbox
  {\vbox{\hrule \@height \fboxrule
          \hbox{\vrule \@width \fboxrule \hskip\fboxsep
          \vbox{\vskip\fboxsep \box\@tempboxa\vskip\fboxsep}\hskip
                 \fboxsep }%
                  }}}}
\makeatother

\makeatletter
\long\def\LL#1{\leavevmode\setbox\@tempboxa\hbox{#1}\@tempdima\fboxrule
    \advance\@tempdima \fboxsep \advance\@tempdima \dp\@tempboxa
   \hbox{\lower \@tempdima\hbox
  {\vbox{
          \hbox{\vrule \@width \fboxrule \hskip\fboxsep
          \vbox{\vskip\fboxsep \box\@tempboxa\vskip\fboxsep}\hskip
                 \fboxsep }%
                 \hrule \@height \fboxrule}}}}
\makeatother

\let \ttorg \tt \def \tt{\ttorg \obeyspaces}

\begin{document}


\title{\Large\bf Knot Diagrammatics}

\author{Louis H. Kauffman \\ Department of Mathematics, Statistics and Computer
Science \\ University of Illinois at Chicago \\ 851 South Morgan Street\\
Chicago, IL, 60607-7045}

\maketitle

\thispagestyle{empty}

\subsection*{\centering Abstract}

{\em This paper is a survey of knot theory and invariants of knots and links from the point of view of categories of diagrams.
The topics range from foundations of knot theory to virtual knot theory and topological quantum field theory. }

\section{Introduction}

This paper is an exploration of the theme of knot diagrams. I have deliberately
focused on  basics in a number of interrelated domains. In most cases some
fundamentals are done in a new or in a more concise way. Some parts of the paper
are expository to fill in the context. This exposition is an outgrowth of the
series of lectures \cite{KD} that I gave in Tokyo at the {\em Knots 96} conference in the
summer of 1996. The present exposition goes considerably farther than those lectures, and includes material
on virtual knot theory and on links that are undetectable by the Jones polynomial.
\vspace{3mm}

The paper is divided into seven sections (counting from section 2 onward).  
Section 2, on the Reidemeister moves,
gives a proof of Reidemeister's basic theorem (that the three Reidemeister moves
on diagrams generate ambient isotopy of links in three-space). A discussion on
graph embeddings extends Reidemeister's theorem to graphs and proves the
appropriate moves for topological and rigid vertices. We hope that this
subsection fills in some gaps in the literature. Section 3 discusses Vassiliev
invariants and invariants of rigid vertex graphs. This section is expository,
with discussions of the four-term relations, Lie algebra weights, relationships
with the Witten functional intergal, and combinatorial constructions for some
Vassiliev invariants. The discussion raises some well-known problems about
Vassiliev invariants. The section on the functional integral introduces a useful
abstract tensor notation that helps in understanding how the Lie algebra weight
systems are related to the functional integral. Sections 4 and 5 are based on a
reformulation of the Reidemeister moves so that they work with diagrams arranged
generically transverse to a special  direction in the plane. We point out how the
technique by which we proved Reidemeister's Theorem (it is actually
Reidemeister's original technique) generalises to give these moves as well. The
moves with respect to a vertical are intimately related to quantum link
invariants and to Hopf algebras. Section 4 is a quick exposition of quantum link
invariants , their relationship with Vassiliev invariants, classical Yang-Baxter
equation and infinitesimal braiding relations. Again, this provides the context
to raise many interesting questions.  Section 5 is a very concise introduction to
the work of the author, David Radford and Steve Sawin  on invariants of
three-manifolds from finite dimensional Hopf algeras. We touch on the question of
the relationship of this work to the Kuperberg  invariants. Section 6 is a
discussion of the Temperley-Lieb algebra. Here we give a neat proof of the
relation structure in the Temperley-Lieb monoid via piecewise linear diagrams.
The last part of this section explains the relationship of the Temperley-Lieb
monoid to parenthesis structures and shows how this point of view can be used to 
relate parentheses to the pentagon and the Stasheff polyhedron. 
Section seven discusses virtual knot theory. Section eight discusses
the construction of links that while linked, have the same Jones polynomial as the unlink.
\vspace{3mm}

\noindent {\bf Acknowledgement.} Most of this effort was sponsored by the Defense
Advanced Research Projects Agency (DARPA) and Air Force Research Laboratory, Air
Force Materiel Command, USAF, under agreement F30602-01-2-05022. Some of this
effort was also sponsored by the National Institute for Standards and Technology
(NIST). The U.S. Government is authorized to reproduce and distribute reprints
for Government purposes notwithstanding any copyright annotations thereon. The
views and conclusions contained herein are those of the authors and should not be
interpreted as necessarily representing the official policies or endorsements,
either expressed or implied, of the Defense Advanced Research Projects Agency,
the Air Force Research Laboratory, or the U.S. Government. (Copyright 2004.) 
It gives the author great pleasure to acknowledge support from NSF Grant DMS-0245588.
\bigbreak

\section{Reidemeister Moves}

Reidemeister \cite{Reidemeister} discovered a simple set of moves on link
diagrams that captures the concept of ambient isotopy of knots in three-dimensional 
space. There are three basic Reidemeister moves. 
Reidemeister's theorem states that two diagrams represent ambient isotopic knots (or links) if
and only if there is a sequence of Reidemeister moves taking one diagram to the
other. The Reidemeister moves are illustrated in Figure 1. \vspace{3mm}

\centerline{\includegraphics[scale=1.0]{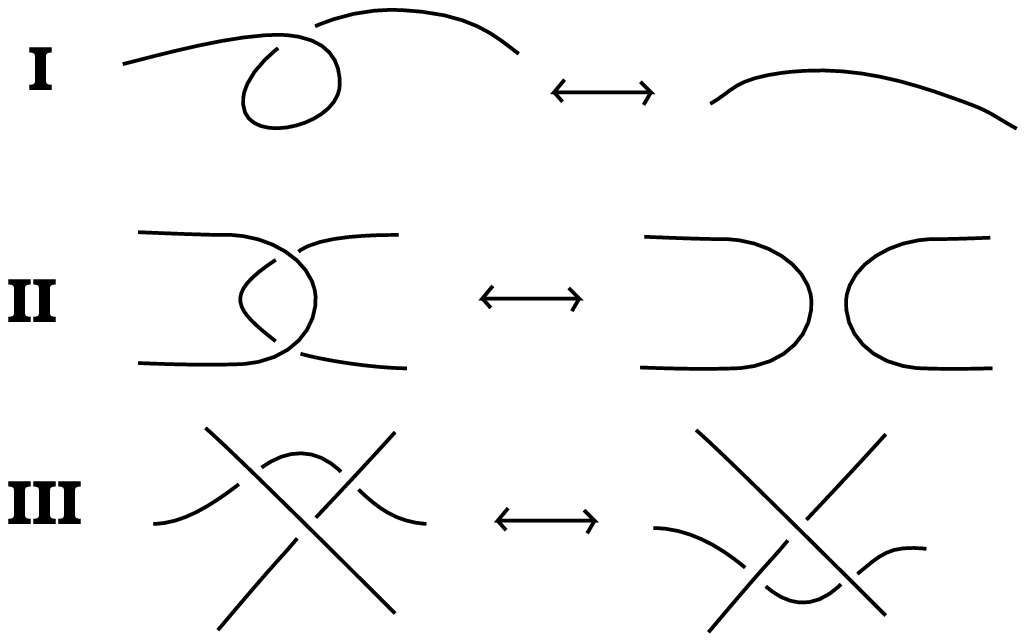}}
\v

\begin{center}
{\bf Figure 1 - Reidemeister Moves}
\end{center}
\vvvv

Reidemeister's three moves are  interpreted as  performed on a larger diagram in
which the small diagram shown is a literal part. Each move is performed without
disturbing the rest of the diagram.  Note that this means that each move occurs,
up to topological deformation, just as it is shown in the diagrams in Figure 1. 
There are no extra lines in the local diagrams. For example, the equivalence (A)
in Figure 2 is {\bf not} an instance of a single first Reidemeister move. Taken
literally, it factors into a move $II$ followed by a \mbox{move $I$.}
\vvv

\centerline{\includegraphics[scale=1.0]{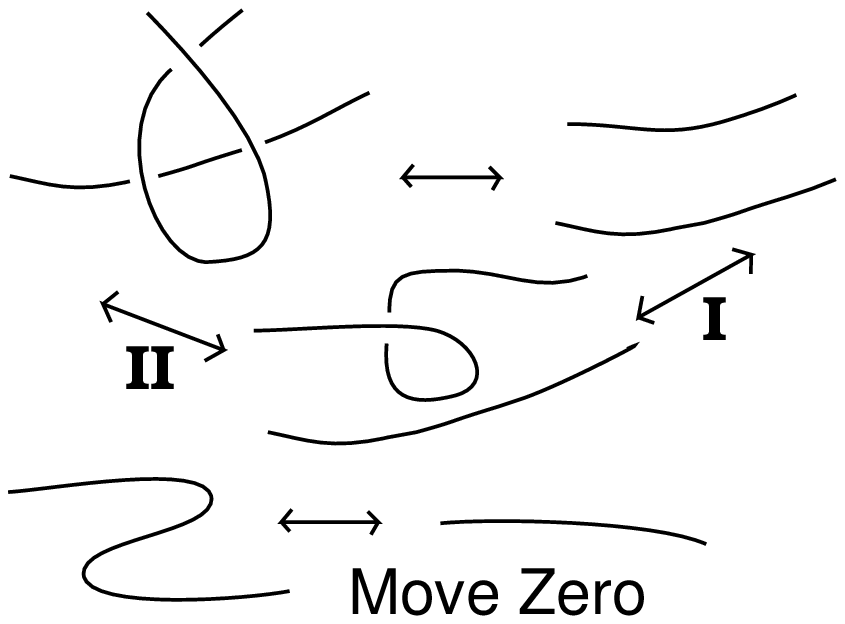}}
\vvvv\vvv

\begin{center}
{\bf Figure 2 - Factorable Move, Move Zero}
\end{center}


Diagrams are always subject to topological deformations in the plane that
preserve the structure of the crossings. These deformations could be designated
as ``Move Zero". See Figure 2. \vspace{3mm}

A few exercises with the Reidemeister moves are in order.  First of all, view the
diagram in Figure 3. It is unknotted and you can have a good time finding a
sequence of Reidemeister moves that will do the trick. Diagrams of this type are
produced by tracing a a curve and always producing an undercrossing at each
return crossing. This type of knot is called a {\em standard unknot}.  Of course
we see clearly that a standard unknot is unknotted by just {\em pulling} on it,
since it has the same structure as a coil of rope that is wound down onto a flat
surface.
\vvvv\vvvv

\centerline{\includegraphics[scale=1.1]{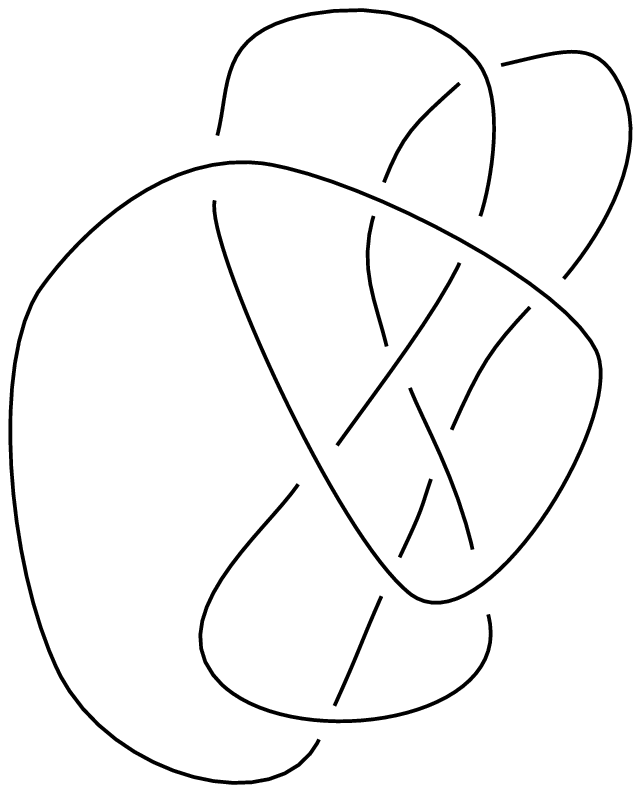}}
\vvv

\begin{center}
{\bf Figure 3 - Standard Unknot}
\end{center}
\vvvv\vvvv

Can one recognise unknots by simply looking for sequences of Reidemeister moves
that undo them?  This would be easy if it were not for the case that there are
examples of unknots that require some moves that increase the number of crossings
before they can be subsequently decreased. Such an demonic example is illustrated
in Figure 4.
\newpage
\vspace*{-.5in}

\centerline{\includegraphics[scale=1.1]{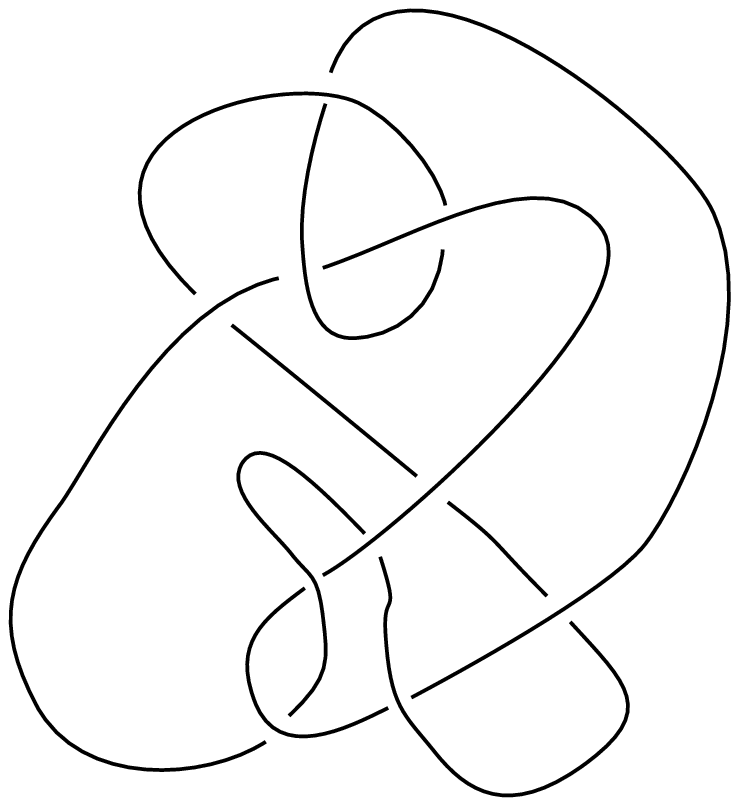}}
\vvvv

\begin{center}
{\bf Figure 4 - A Demon}
\end{center}
\vvvv\vvvv

It is generally not so easy to recognise unknots. However, here is a tip: Look
for {\em macro moves} of the type shown in Figure 5. In a macro move, we identify
an arc that passes entirely under some piece of the diagram (or entirely over)
and shift this part of the arc, keeping it under (or over) during the shift. 
In Figure 5, we illustrate a macro move on an arc that passes under a piece of the diagram 
that is indicated by arcs going into a circular region. A more general marcro move is possible where the 
moving arc moves underneath one layer of diagram, and at the same time, over another layer of diagram.
Macro moves often allow a reduction in the number of crossings even though the
number of crossings will increase during a sequence of Reidemeister moves that
generates the macro move. \vspace{3mm}


As shown in Figure 5, the macro-move includes as a special case both the second
and the third Reidemeister moves, and it is not hard to verify that a macro move
can be generated by a sequence of type II and type III Reidemeister moves.  It is
easy to see that the type I moves can be left to the end of any deformation.  The
demon of Figure 4  is easily demolished by macro moves, and from the point of
view of macro moves the diagram never gets more complicated. \vspace{3mm}
\newpage

\centerline{\includegraphics[scale=1.0]{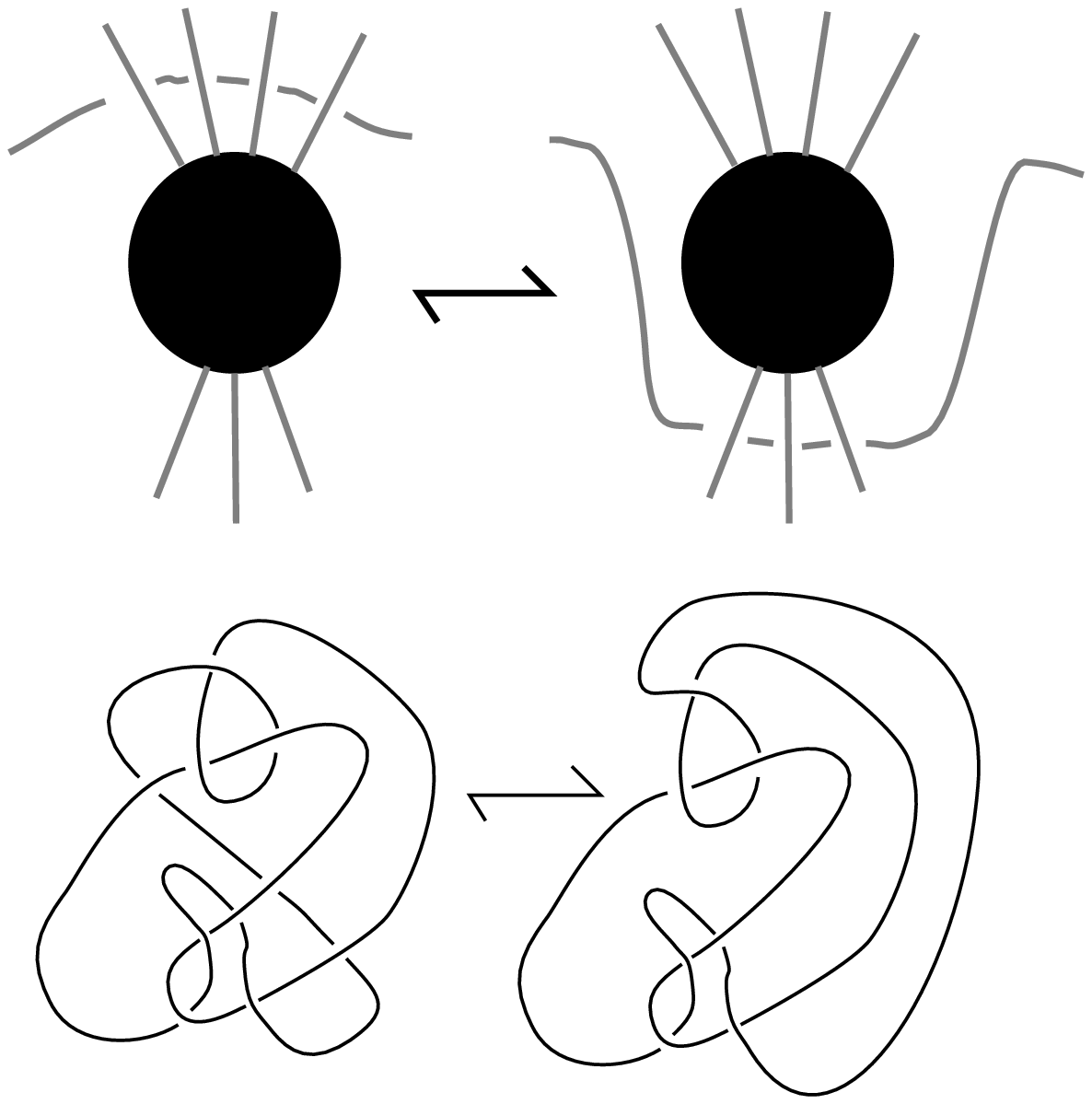}}
\vv

\begin{center}
{\bf Figure 5 - Macro Move }
\end{center}
\vvvv

Let's say that an knot can be {\em reduced} by a set of moves if it can be
transformed by these moves to the unknotted circle diagram through diagrams that
never have more crossings than the original diagram. Then we have shown that
there are diagrams representing the unknot that cannot be reduced by the
Reidemeister moves.  On the other hand, I do not know whether unknotted diagrams
can always be reduced by the macro moves in conjunction with the first
Reidemeister move. If this were true it would give a combinatorial way to
recognise the unknot. \vspace{3mm}

\noindent {\bf Remark.} In fact, there is a combinatorial way to recognise the unknot based on 
a diagrams and moves. In \cite{Dynnikov} I. A. Dynnikov finds just such a result, using piecewise 
linear knot diagrams with all ninety degree angles in the diagrams, and all arcs in the diagram either 
horizontal or vertical. The interested reader should consult his lucid paper.
\bigbreak

\subsection{Reidemeister's Theorem}

We now  indicate how Reidemeister proved his Theorem. \vspace{3mm}

An embedding of a knot or link in three-dimensional space is said to be {\em
piecewise linear} if it consists in a collection of straight line segments joined
end to end.  Reidemeister started with a {\em single} move in three-dimensional
space for piecewise linear knots and links.   Consider a point in the complement
of the link, and an edge in the link such that the surface of the triangle formed
by the end points of that edge and the new point is not pierced by any other edge
in the link. Then one can replace the given edge on the link by the other two
edges of the triangle, obtaining a new link that is ambient isotopic to the
original link. Conversely, one can remove two consecutive edges in the link and
replace them by a new edge that goes directly from initial to final points,
whenever the triangle spanned by the two consecutive edges is not pierced by any
other edge of the link. This triangle replacement constitutes Reidemeister's
three-dimensional move. See Figure 6.  It can be shown that two piecewise linear
knots or links are ambient isotopic in three-dimensional space if and only if
there is a  sequence of Reidemeister triangle moves from one to the other. This
will not be proved here. At the time when Reidemeister wrote his book,
equivalence via three-dimensional triangle moves was taken as the definition of
topological equivalence of links. \vspace{3mm}
\vvvv\vvv

\centerline{\includegraphics[scale=1.1]{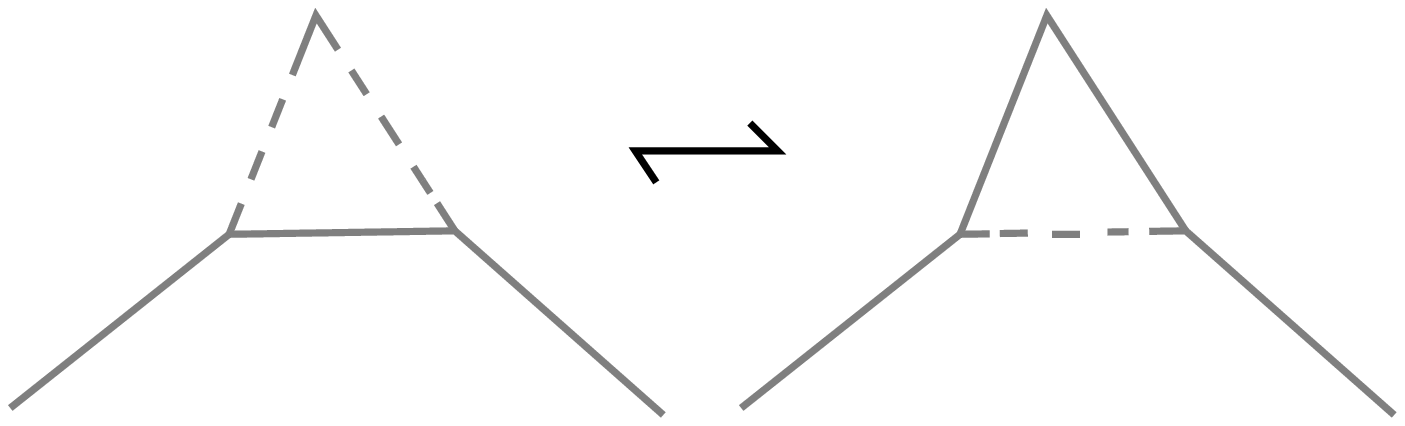}}
\vvvv\vvv

\begin{center}
{\bf Figure 6 - Triangle Move }
\end{center}
\vvvv\vvv

It can also be shown that tame knots and links have piecewise linear
representatives in their ambient isotopy class. It is sufficient for our purposes
to work with piecewise linear knots and links. Reidemeister's planar moves then
follow from an analysis of the shadows projected  into the plane by Reidemeister
triangle moves in space.  Figure 7  gives a hint of this analysis. The result is
a reformulation of the three-dimensional problems of knot theory to a
combinatorial game in the plane.
\newpage

\centerline{\includegraphics[scale=1.0]{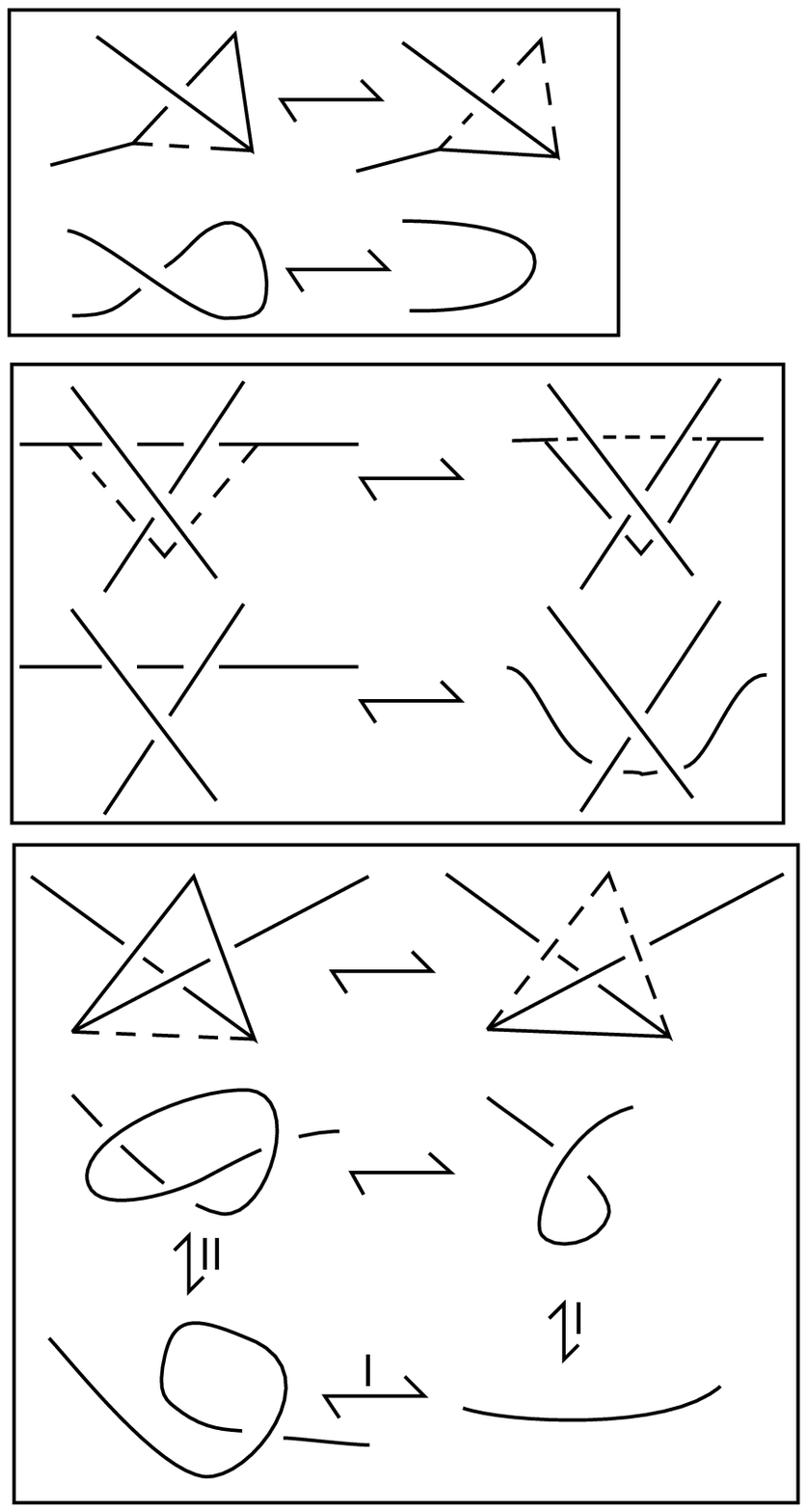}}
\vvvv

\begin{center}
{\bf Figure 7 - Shadows}
\end{center}

To go beyond the hint in Figure 7 to a complete proof that Reidemeister's planar
moves suffice involves preliminary remarks about subdivision.  The simplest
subdivision that one wants to be able to perform on a piecewise linear link is
the placement of a new vertex at an interior point of an edge - so that that edge
becomes two edges in the subdivided link. Figure 8 shows how to accomplish this
subdivision via triangle moves.
\vvvv\vvvv\vvvv

\centerline{\includegraphics[scale=1.3]{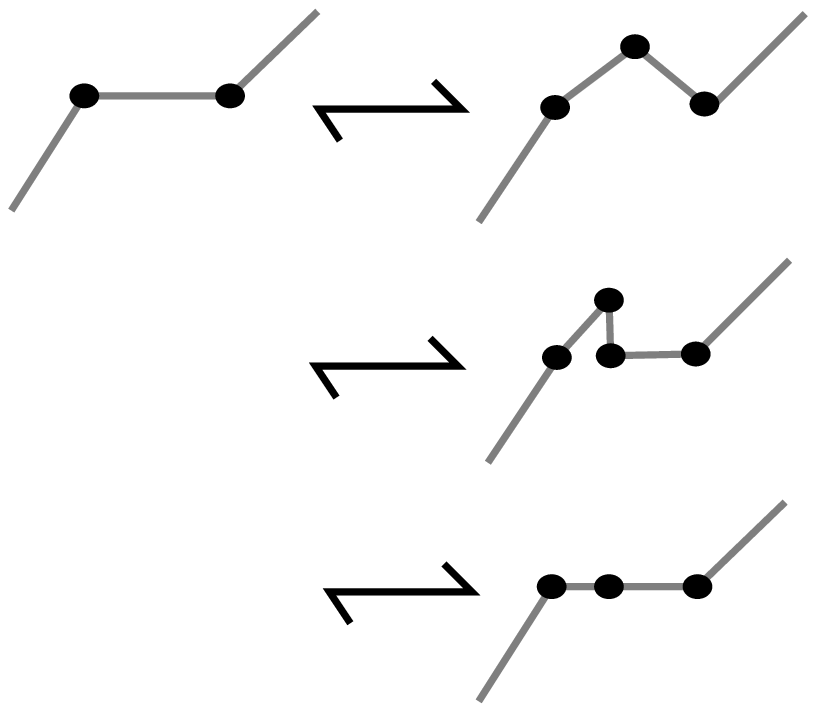}}
\vvvv\vvvv

\begin{center}
{\bf Figure 8 - Subdivision of an Edge}
\end{center}
\vvvv\vvvv\vvvv

Any triangle move can be factored into a sequence of smaller triangle moves
corresponding to a simplicial subdivision of that triangle.  This is obvious,
since the triangles in the subdivision of the large triangle that is unpierced by
the link are themselves unpierced by the link. \vspace{3mm}

To understand how the Reidemeister triangle move behaves on diagrams it is
sufficient to consider a projection of the link in which the triangle is
projected to a non-singular triangle in the plane. Of course, there may be many
arcs of the link also projected upon the interior of the projected triangle.
However, by using subdivision, we can assume that the cases of the extra arcs are
as shown in Figure 9. In Figure 9 we have also shown how each of these cases can
be accomplished by (combinations of) the three Reidemeister moves. This proves
that a projection of a single triangle move can be accomplished by a sequence of
Reidemeister diagram moves.
\newpage
\vspace*{-.5in}

\centerline{\includegraphics[scale=0.95]{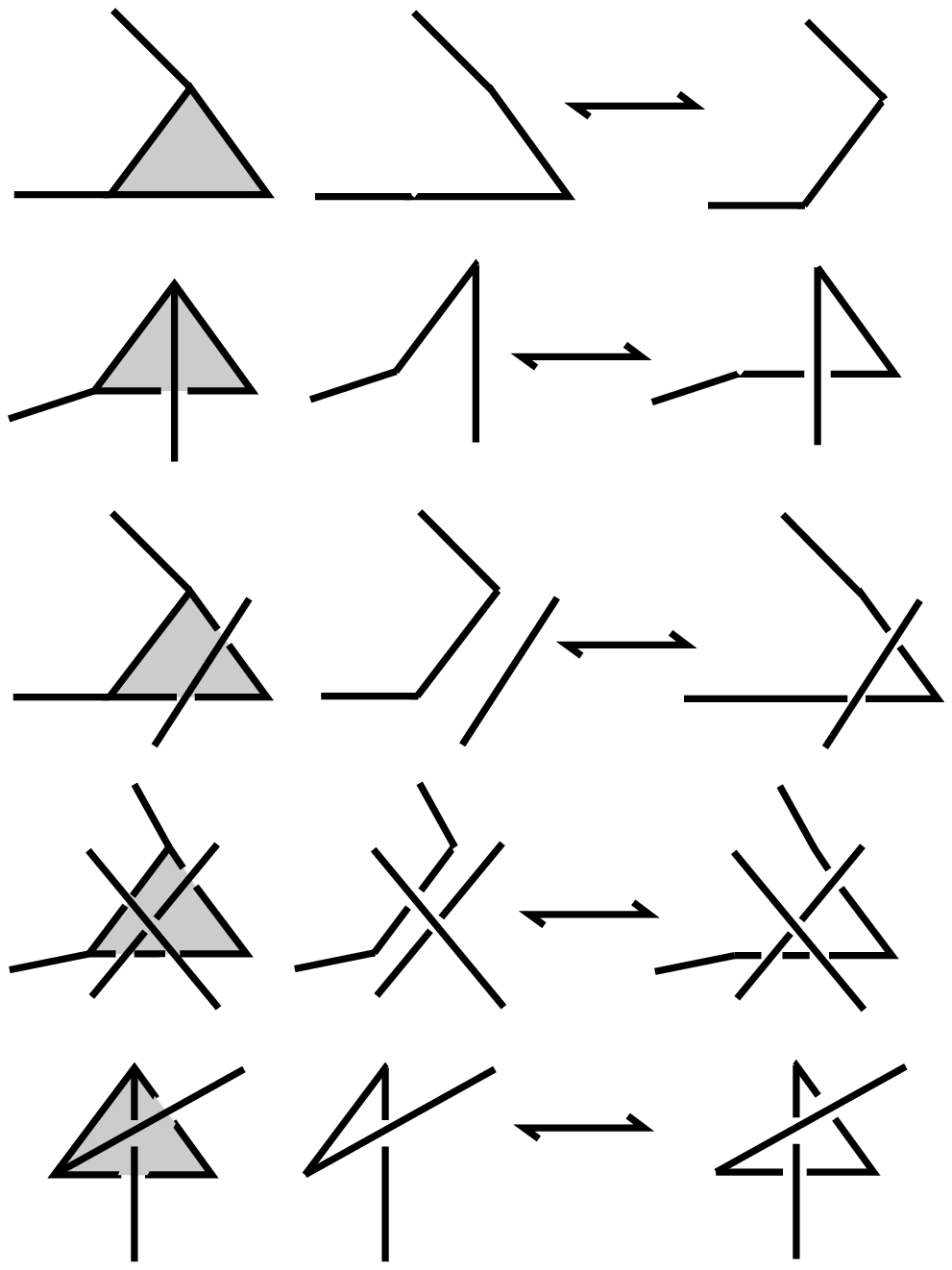}}
\vvv

\begin{center}
{\bf Figure 9 - Projections of Triangle Moves}
\end{center}
\vvv

A piecewise linear isotopy consists in a finite sequence of triangle moves. There
exists a direction in three-dimensional space that makes a non-zero angle with
each of theses triangles and is in general position with the link diagram.
Projecting to the plane along this direction makes it possible to perform the
entire ambient isotopy in the language of projected triangle moves. Now apply the
results of the previous paragraph and we conclude \vspace{3mm}

\noindent {\bf Reidemeister's Theorem.}  If two links are piecewise linearly
equivalent (ambient isotopic), then there is a sequence of Reidemeister diagram
moves taking a projection of one link to a projection of the other. \vspace{3mm}

Note that the proof tells us that the two diagrams can be obtained from one
spatial projection direction for the entire spatial isotopy. It is obvious that
diagrams related by Reidemeister moves represent ambient isotopic links.
Reidemeister's Theorem gives a complete combinatorial description of the topology
of knots and links in three-dimensional space.
\newpage

\subsection{Graph Embeddings} Let $G$ be a (multi-)graph.  That is, $G$ is a
finite abstract graph with, possibly, a multiplicity of edges between any two of
its vertices. Now consider the embeddings of $G$ in Euclidean three-space
$R^{3}.$  In the category of topological embeddings, any edge of $G$ can acquire
local knotting as shown in Figure 10. On top of this there is the possibility of
global knotting that results from the structure of the graph as a whole.
\vvvv\vvvv\vvvv

\centerline{\includegraphics[scale=1.2]{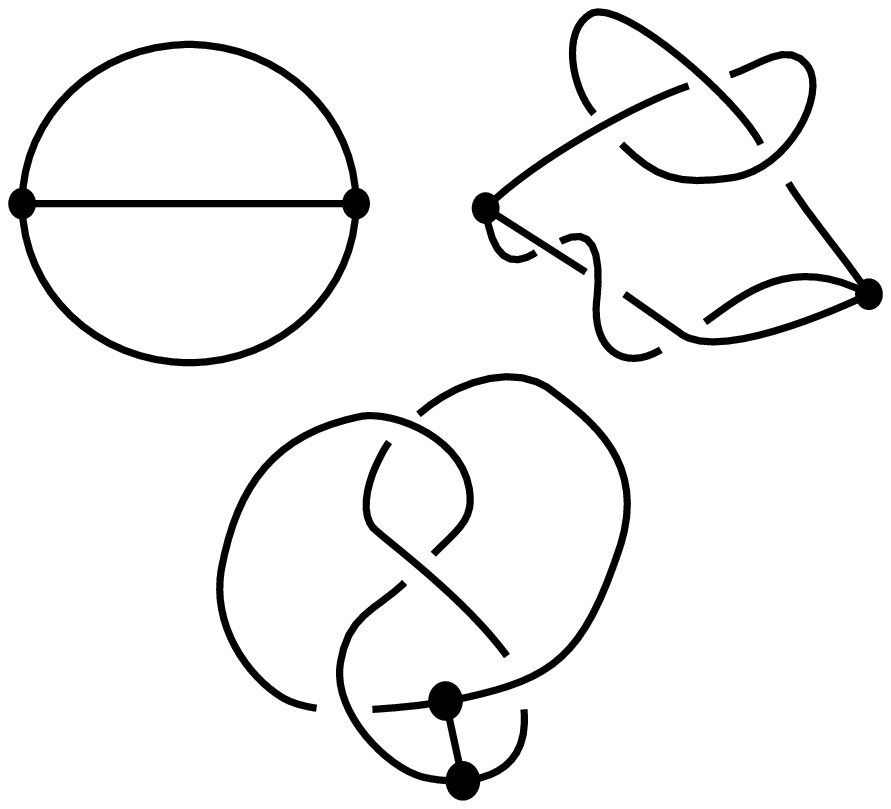}}
\vvvv\vvvv

\begin{center}
{\bf Figure 10 - Graph Embedding}
\end{center}
\vvvv\vvvv\vvvv

Topological or piecewise linear ambient isotopy of graph embeddings is
complicated by the fact that arbitrary braiding can be created or destroyed at a
vertex, as illustrated in Figure 11.
\newpage

\centerline{\includegraphics[scale=1.2]{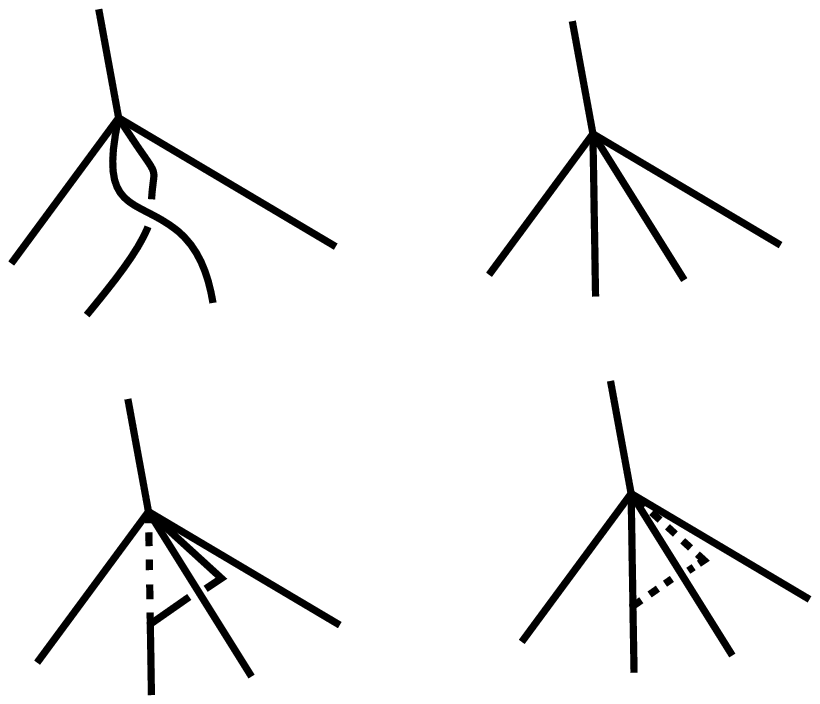}}
\vvvv\vvv

\begin{center}
{\bf Figure 11 - Braiding at a Vertex}
\end{center}
\vvvv\vvvv

For this reason, it is useful to consider ways to restrict the allowed movement
in the neighborhood of a vertex. One way to accomplish this is to decree that
each vertex will come equipped with a specific cyclic order of the edges meeting
the vertex. This cyclic order can be instantiated on the boundary of a disk, and
the graph replaced by a configuration of disks with cyclic orders of marked
points along their boundaries. The edges of the original graph are replaced by
edges that go from one disk to another terminating in the marked points. Call
such an arrangement a {\em rigid vertex graph} $G$.  If $G$ is a rigid vertex
graph, then we consider embeddings of $G$ where the disks are embedded metrically
while the (graphical) edges are embedded topologically.  A {\em rigid vertex
isotopy} of one $RV$ ($RV$ will stand for rigid vertex) embedding $G$ to another
$G'$  is  a combination of ambient isotopies of the embedded edges of the graph
(the {\em strings} of the graph) relative to their endpoints on the disks,
coupled with affine motions of the disks (carrying along the strings in ambient
isotopy).  An affine motion of a disk is a combination of parallel translations
of the disk along a given direction in three-space and rotations of the disk
about an axis through its center.  We can think of a given disk as embedded
inside a standard three-ball with the strings from the disk emanating straight to
the boundary of the three-ball. Each basic affine motion is assumed to leave the
points on the boundary of the containing three-ball fixed.  Thus the types of
affine motion are as illustrated in Figure 12.
\newpage

\centerline{\includegraphics[scale=1.0]{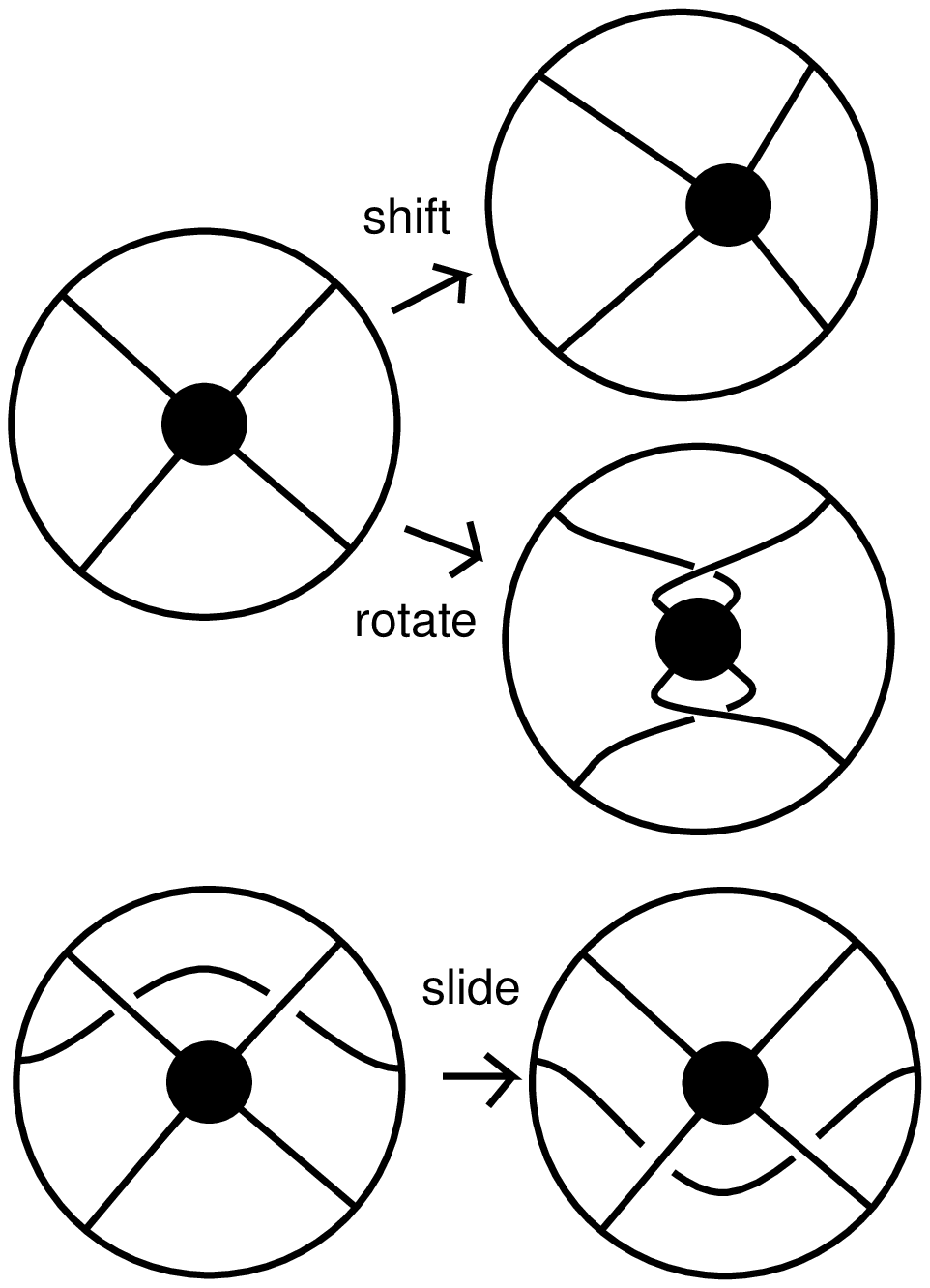}}
\vvvv

\begin{center}
{\bf Figure 12 - Rigid Vertex Graphs and Affine Motions}
\end{center}

We will give versions of the Reidemeister moves for both topological isotopy and
rigid vertex isotopy of embedded graphs. In the topological case the extra moves
are illustrated in Figure 13. Here we have indicated the elementary braiding at a
vertex and slide moves that take an edge underneath a vertex.  The proof that
these moves suffice is a generalisaton of our original proof of the Reidemeister
moves. That is, we model the graph embeddings by piecewise linear embeddings.
This may entail subdividing the edges of the original graph so that those edges
can have enough flexibility to sustain a given topological conformation. Thus,
when we speak of a piecewise linear embedding of a given graph, we mean a
piecewise linear embedding of a graph that is obtained from the given graph by
subdividing some of its edges. \vspace{3mm}

\centerline{\includegraphics[scale=0.85]{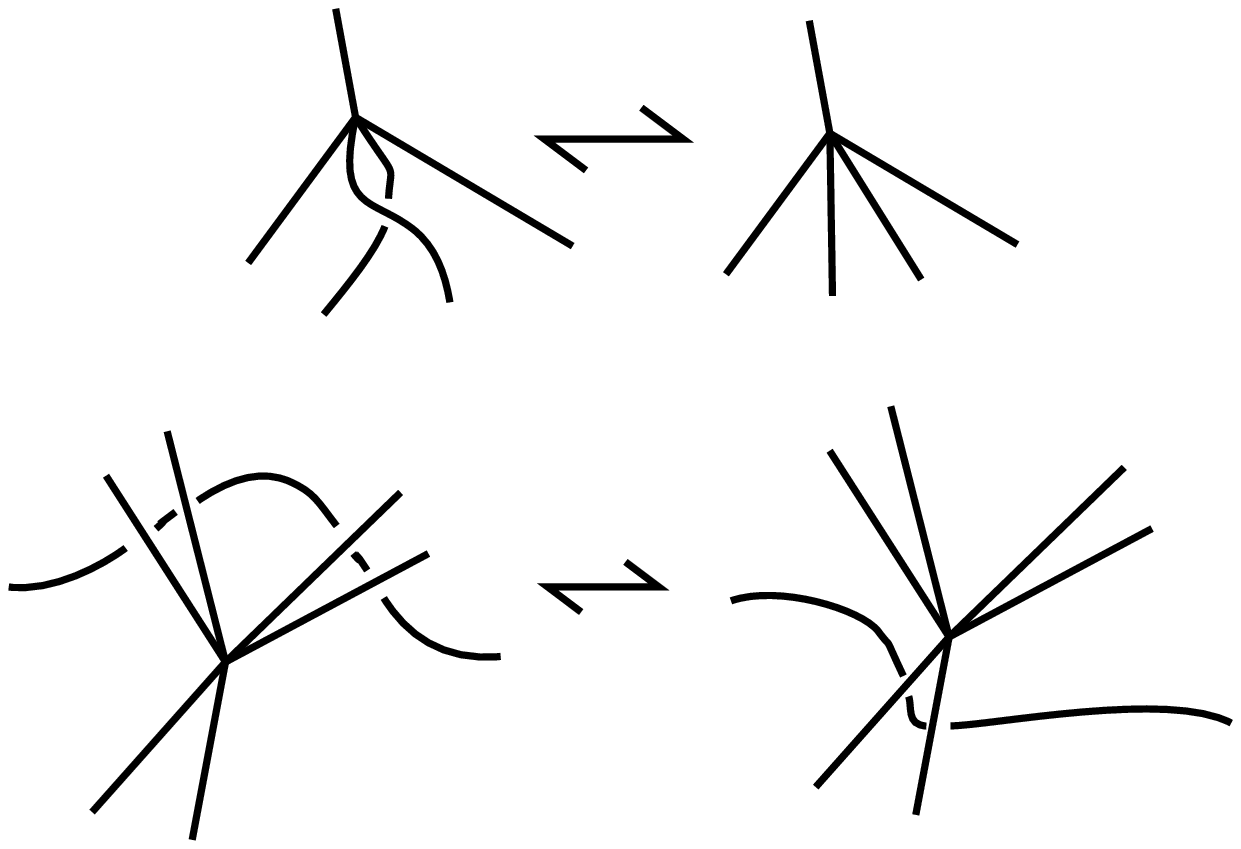}}
\vvvv

\begin{center}
{\bf Figure 13 - Extra Moves For Topological Isotopy of Graphs}
\end{center}

Piecewise ambient isotopy of   graph embeddings is defined exactly as in the case
of piecewise linear isotopy for knots and links. The same projection arguments
apply and the extra moves are obtained from the three-dimensional triangle move
as illustrated in Figure 14. This  completes the proof of our assertion about the
topological  Reidemeister moves for graphs. \vspace{3mm}

\centerline{\includegraphics[scale=0.85]{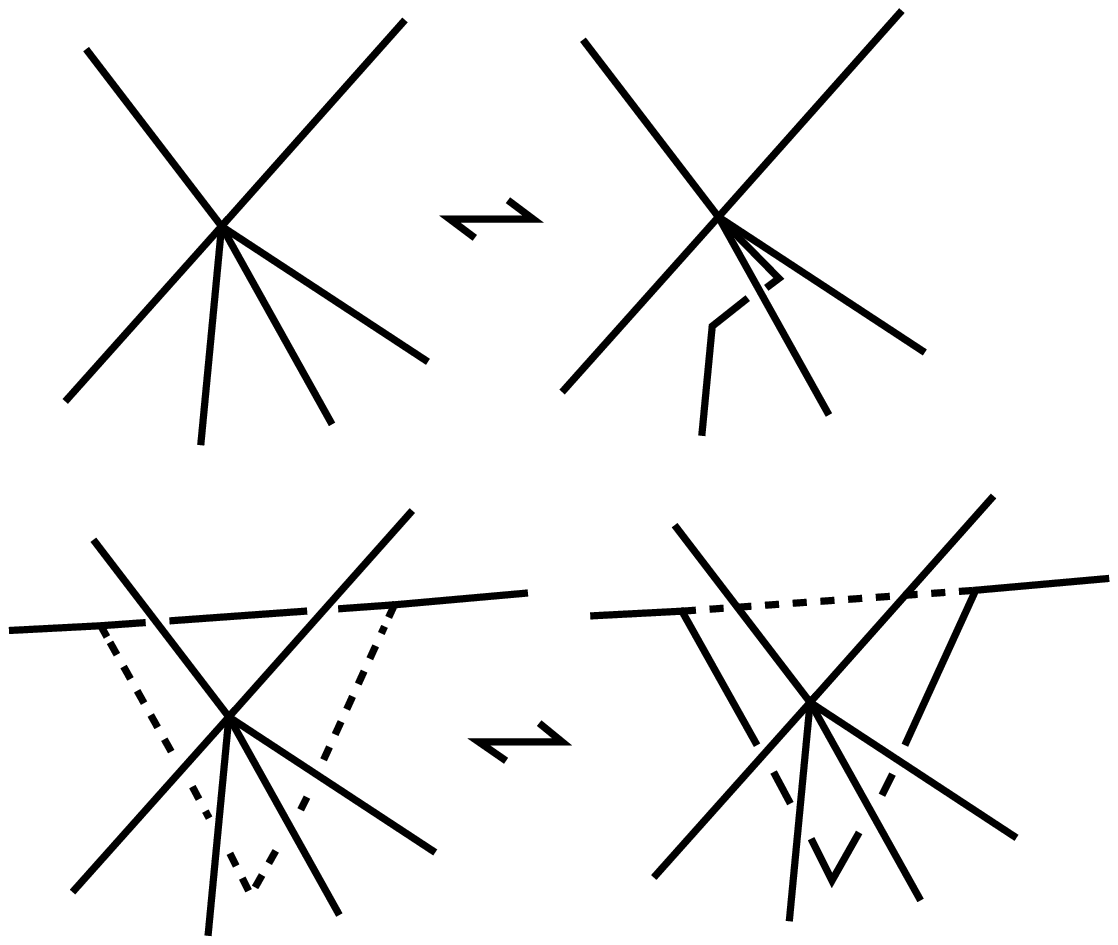}}
\vvvv

\begin{center}
{\bf Figure 14 - PL Isotopy Inducing  Topological Graphical Moves}
\end{center}
\newpage
\vspace*{-.25in}

Consider rigid vertex isotopy of rigid vertex ($RV$) graphs. We will assume that
the topological moves are performed in the piecewise linear setting. Thus
subdivisions of the edges of the graph can be produced. Basic translational
affine moves of the embedded disks can have piecewise linear starting and ending
states by drawing straight lines from the marked points on the disk boundaries to
the corresponding points in the containing balls.  Rotatory moves with the center
of a disk as axis can also have piecewise linear starting and ending states by
taking the braiding that is induced by the rotation and suitably subdividing it.
These remarks show that RV isotopy can be achieved in the PL category.

The next point to consider is the result of projection of an RV isotopy on the
corresponding diagrams. A sequence of elementary RV isotopies from a graph $G$ to
a graph $G'$ has associated with it a direction of projection so that each PL
triangle move has its triangle projected to a non-singular triangle in the plane
and each affine move has its disk projected to a non-singular disk in the plane.
In the case of the affine moves we can assume that the before to after appearance
of the disk and its corresponding containing ball will represent either a
topological identity map (albeit an affine shift) or a rotation about the disk
axis by $\pi$ radians. (Higher multiples of $\pi$ can be regarded as iterates of
a $\pi$ rotation.) Therefore the basic $\pi$ rotation can be schematized as shown
in Figure 15.   Figure 15 illustrates the moves that we need to add to the
Reidemeister moves to obtain a planar diagram version of RV isotopy. The
remaining moves in Figure 15 follow from the same projection arguments that we
have used earlier in this section.   This completes the construction of the
diagrammatic calculus for RV isotopy.
\vvv

\centerline{\includegraphics[scale=0.75]{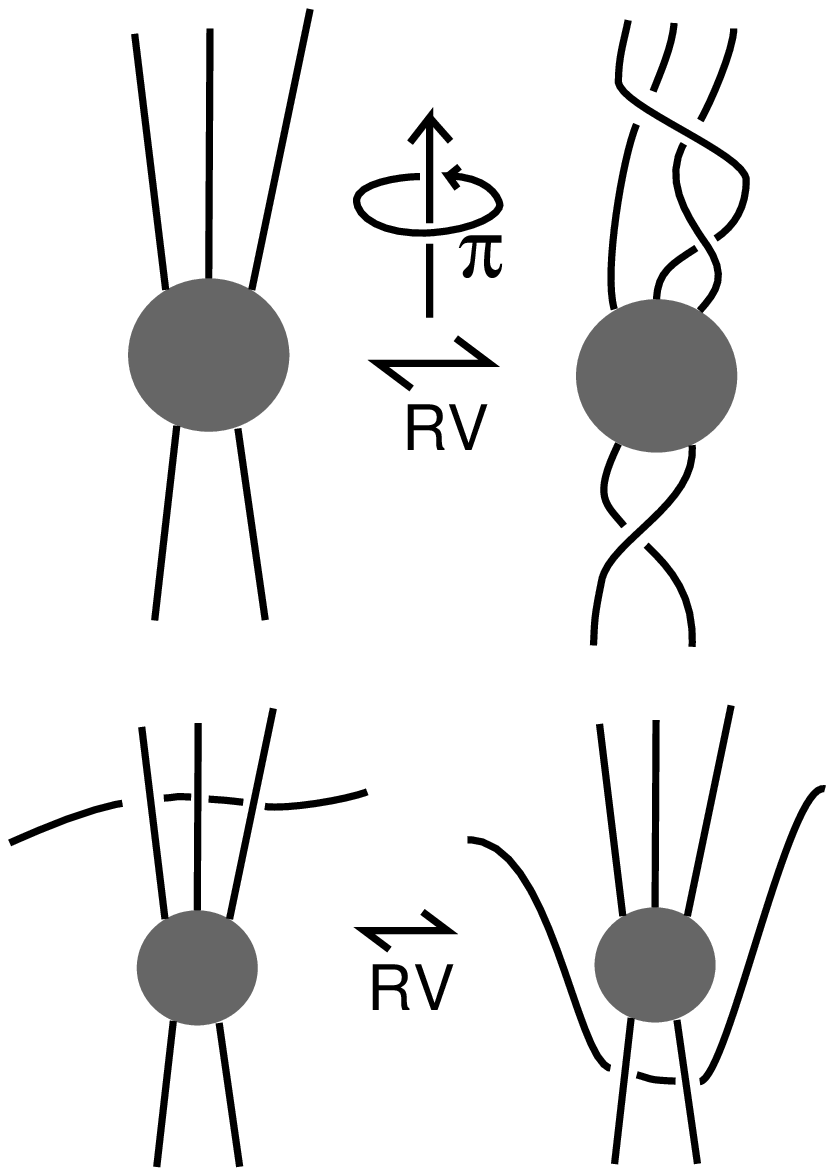}}
\vv

\begin{center}
{\bf Figure 15 - Diagrammatic Rigid Vertex Isotopy}
\end{center}
\newpage
\vspace*{-.2in}

Note that the generating moves for rigid vertex graph isotopy are almost the same
as the generating moves for topological graph isotopy, except that the braiding
at the vertex  in the rigid vertex case comes from the twisting the disk as a
whole. This circumstance makes the construction of invariants of rigid vertex
graphs much easier.  We will discuss constructions of such invariants in the next
section. In section 4 we will return to the Reidemeister moves and reformulate
them once again for the sake of quantum link invariants.
\vvv

\section{Vassiliev Invariants and Invariants of Rigid Vertex Graphs}

If  $V(K)$ is a  (Laurent polynomial valued,  or more generally - commutative
ring valued) invariant of knots,  then it can be naturally extended to an
invariant of rigid vertex graphs by defining the invariant of graphs in terms of
the knot invariant via an ``unfolding"  of the vertex. That is, we can regard the
vertex as a ``black box" and replace it by any tangle of our choice. Rigid vertex
motions of the graph preserve the contents of the black box, and hence entail
ambient isotopies of the link obtained by replacing the black box by its
contents. Invariants of knots and links that are evaluated on these replacements
are then automatically rigid vertex invariants of the corresponding graphs. If we
set up a collection of multiple replacements at the vertices with standard
conventions for the insertions of the tangles, then a summation over all possible
replacements can lead to a graph invariant with new coefficients corresponding to
the different replacements.  In this way each invariant of knots and links
implicates a large collection of graph invariants. See \cite{Kauffman-Graph},
\cite{Kauffman-Vogel}.
\vvvv

\centerline{\includegraphics[scale=1.0]{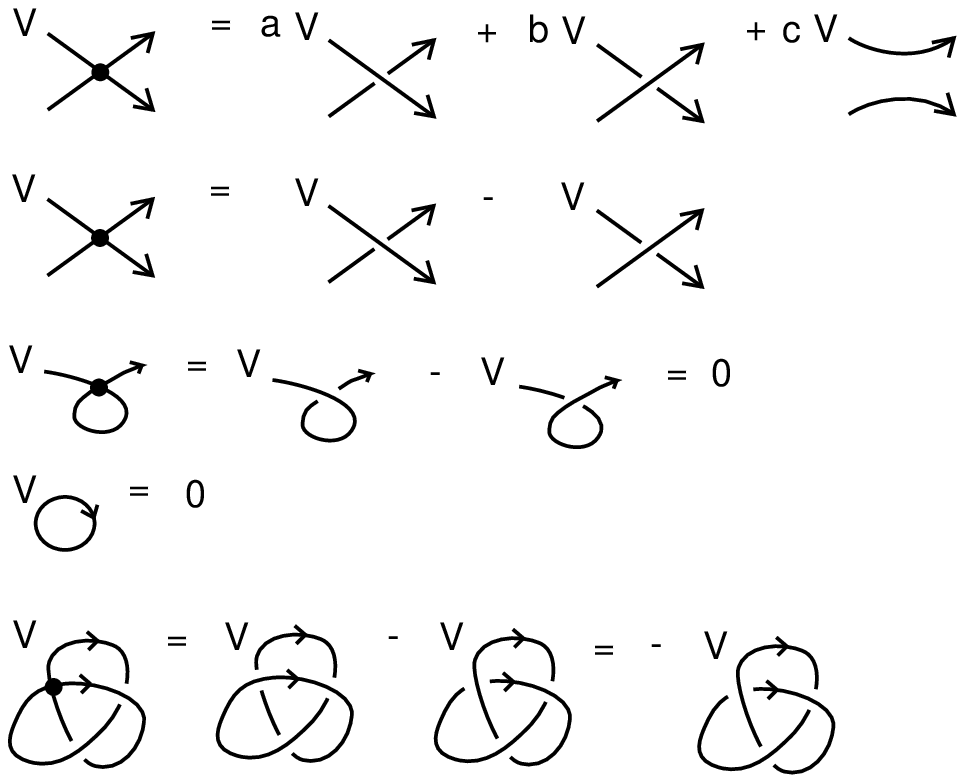}}
\vvv

\begin{center}
{\bf Figure 16 - Graphical Vertex Formulas}
\end{center}
\newpage

The simplest tangle replacements for a 4-valent vertex are the two crossings,
positive and negative, and the oriented smoothing. Let $V(K)$ be any invariant of
knots and links. Extend $V$ to the category of rigid vertex embeddings of 4-valent
graphs by the formula (See Figure 16) $$V(K_{*}) = aV(K_{+}) + bV(K_{-}) +
cV(K_{0})$$ \vspace{3mm}

Here   $K_{*}$   indicates  an embedding with a transversal 4-valent vertex. This
formula means that we define  $V(G)$  for an embedded 4-valent graph  $G$  by
taking the sum

$$V(G) = \sum_{S} a^{i_{+}(S)}b^{i_{-}(S)}c^{i_{0}(S)}V(S)$$

\noindent with the summation over  all knots and links $S$ obtained from  $G$ by replacing
a node of $G$ with either a crossing of positive or negative type, or with  a
smoothing (denoted $0$).  Here $i_{+}(S)$ denotes the number of positive crossings in the
replacement, $i_{-}(S)$ the number of negative crossings in the replacement, and 
$i_{0}(S)$ the number of smoothings in the replacement. It is not hard to see that if $V(K)$  is an  ambient
isotopy invariant of knots, then,  this extension is a rigid vertex isotopy
invariant of graphs.  In rigid vertex isotopy the cyclic order at the vertex is
preserved, so that the vertex behaves like a rigid disk with flexible strings
attached to it at specific points.  See the previous section.
\v

There is a rich class of graph invariants that can be studied in this manner. The
Vassiliev Invariants \cite{Vassiliev}, \cite{Birman and Lin}, \cite{Bar-Natan}
constitute the important special case of these graph invariants where $a=+1$,
$b=-1$ and $c=0.$    Thus  $V(G)$  is a Vassiliev invariant if

$$V(K_{*}) = V(K_{+})  -  V(K_{-}).$$ Call this formula the {\em exchange
identity} for the Vassiliev invariant $V.$ $V$  is said to be of  finite type 
$k$  if  $V(G) = 0$  whenever  $|G| >k$  where $|G|$  denotes the number of
4-valent nodes in the graph $G.$ The notion of finite type is of paramount
significance in studying these invariants. One reason for this is the following
basic Lemma. \vspace{3mm}

\noindent {\bf Lemma.} If a graph $G$ has exactly $k$ nodes, then the value of a
Vassiliev invariant $v_{k}$ of type $k$ on $G$, $v_{k}(G)$, is independent of the
embedding of $G$. \vspace{3mm}

\noindent {\bf Proof.} The different embeddings of $G$ can be represented by link
diagrams with some of the 4-valent vertices in the diagram corresponding to the
nodes of $G$. It suffices to show that the value of $v_{k}(G)$ is unchanged under
switching of a crossing.  However, the exchange identity for $v_{k}$ shows that
this difference is equal to the evaluation of $v_{k}$ on a graph with $k+1$ nodes
and hence is equal to zero. This completes the proof.\vspace{3mm}

The upshot of this Lemma is that Vassiliev invariants of type $k$ are intimately
involved with certain abstract evaluations of graphs with $k$ nodes. In fact,
there are  restrictions (the four-term relations) on these evaluations demanded
by the topology (we shall articulate these restrictions shortly) and it follows
from results of Kontsevich \cite{Bar-Natan} that such abstract evaluations
actually determine the invariants. The invariants derived from classical Lie
algebras are all built from Vassiliev invariants of finite type. All this is
directly related to Witten's functional integral \cite{Witten}. \vspace{3mm}

\noindent {\bf Definition.} Let $v_{k}$ be a Vassiliev invariant of type $k$. The
{\em top row} of $v_{k}$ is the set of values that $v_{k}$ assigns to the set of
(abstract) 4-valent graphs with $k$ nodes. If we concentrate on Vassiliev
invariants of knots, then these graphs are all obtained by marking $2k$ points on
a circle, and choosing a pairing of the $2k$ points. The pairing can be indicated
by drawing a circle and connecting the paired points with arcs. Such a diagram is
called a {\em chord diagram}. Some examples are indicated in Figure 17.
\vspace{3mm}

\centerline{\includegraphics[scale=1.0]{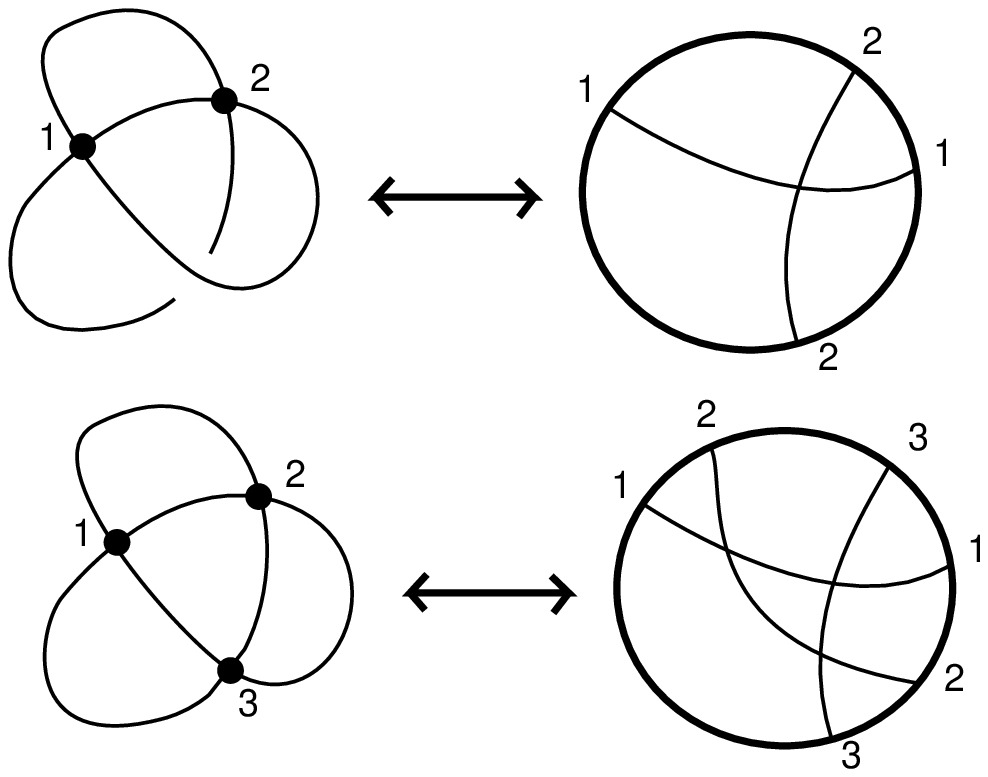}}
\vvv

\begin{center}
{\bf Figure 17 - Chord Diagrams}
\end{center}
\vvv

Note that a top row diagram cannot contain any isolated pairings since this would
correspond to a difference of local curls on the corresponding knot diagram (and
these curls, being isotopic, yield the same Vassiliev invariants. \vspace{3mm}

\noindent{\bf The Four-Term Relation}. (Compare \cite{Stanford}.) Consider a
single embedded graphical node in relation to another embedded arc, as
illustrated in Figure 18. The arc underlies the lines incident to the node at
four points and can be slid out and isotoped over the top so that it overlies the
four nodes. One can also switch the crossings one-by-one to exchange the arc
until it overlies the node.  Each of these four switchings gives rise to an
equation, and the left-hand sides of these equations will add up to zero,
producing a relation corresponding to the right-hand sides. Each term in the
right-hand side refers to the value of the Vassiliev invariant on a graph with
two nodes that are neighbors to each other. See Figure 18. \vspace{3mm}

\centerline{\includegraphics[scale=1.0]{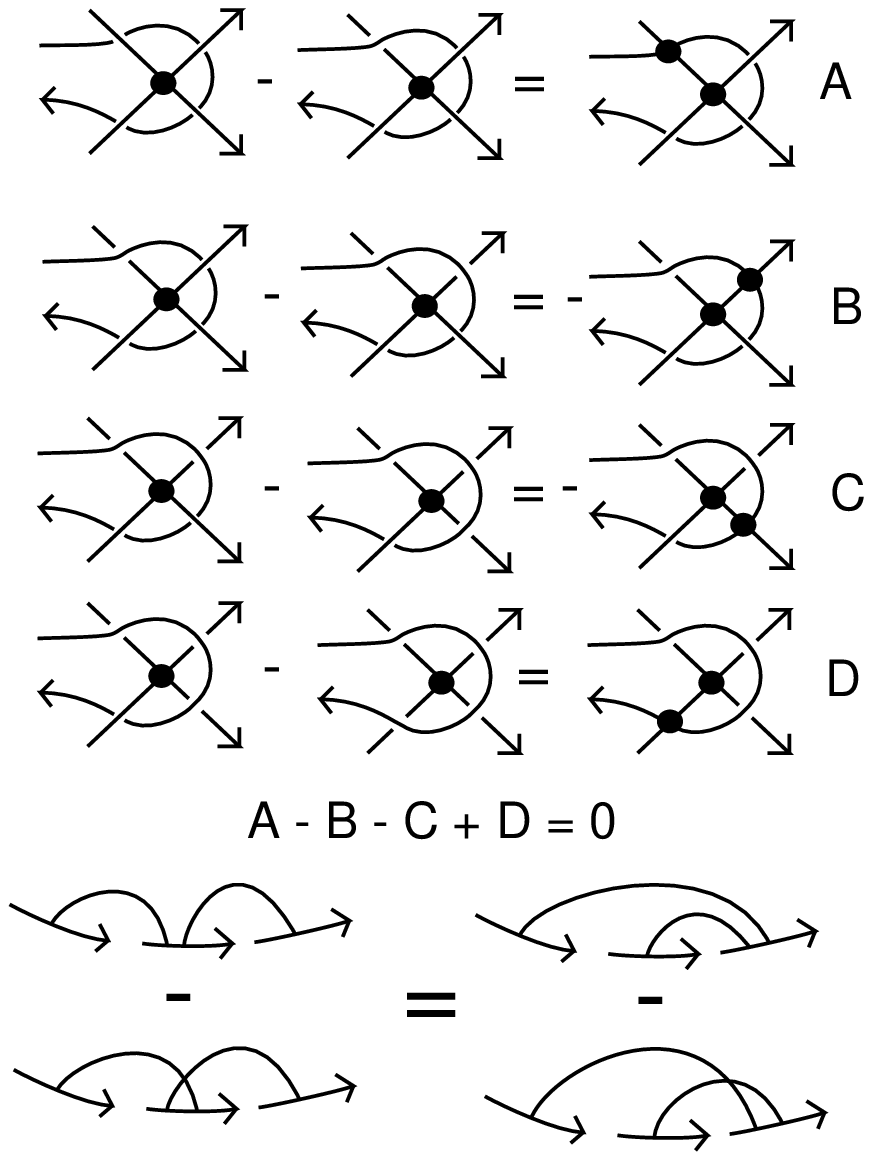}}
\vvvv\vv

\begin{center}
{\bf Figure 18 - The Four Term Relation}
\end{center}
\vvvv\vvv

There is a corresponding 4-term relation for chord diagrams. This is the 4-term
relation for the top row. In chord diagrams the relation takes the form shown at
the bottom of Figure 18. Here we have illustrated only those parts of the chord
diagram that are relevant to the two nodes in question (indicated by two pairs of
points on the circle of the chord diagram). The form of the relation shows the
points on the chord diagram that are immediate neighbors. These are actually
neighbors on any  chord diagram that realizes this form. Otherwise there can be
many other pairings present in the situation. \vspace{3mm}

As an example, consider the possible chord diagrams for a Vassiliev invariant of
type $3.$  There are two possible diagrams as shown in Figure 19. One of these
has the projected pattern of the trefoil knot and we shall call it the {\em
trefoil graph}. These diagrams satisfy the 4-term relation. This shows that one
diagram must have twice the evaluation of the other. Hence it suffices to know
the evaluation of one of these two diagrams to know the top row of a Vassiliev
invariant of type 3. We can take this generator to be the trefoil graph
\newpage
\vspace*{-.4in}

\centerline{\includegraphics[scale=0.85]{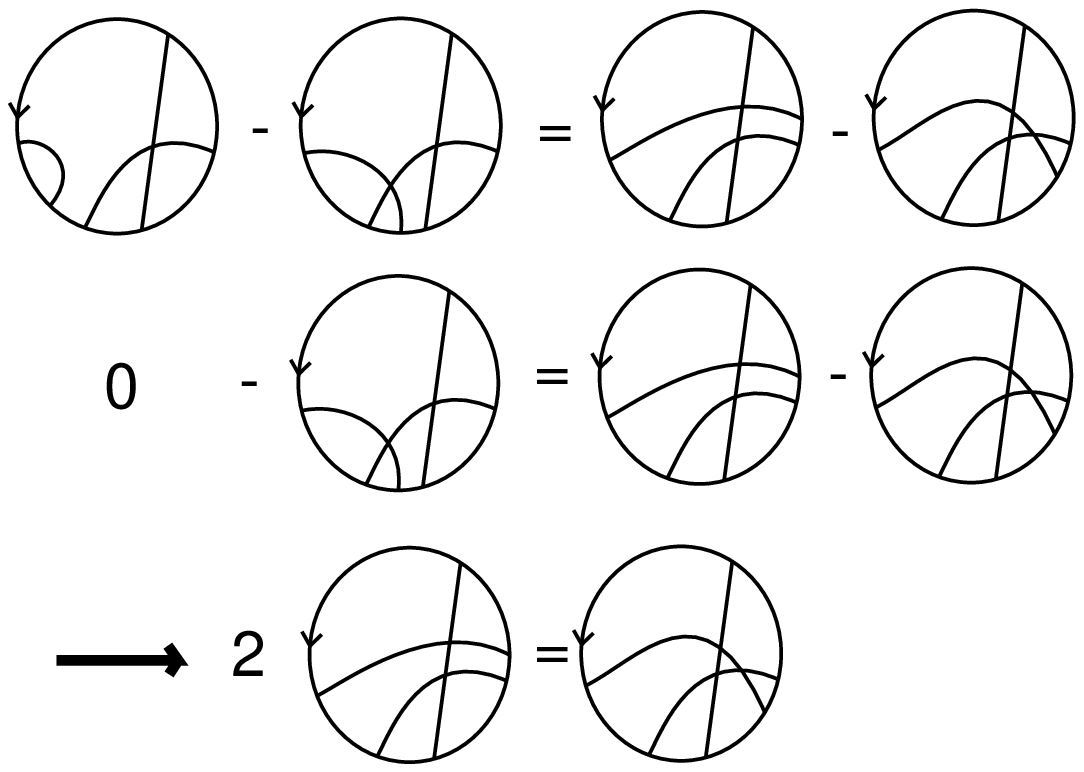}}
\vvvv

\begin{center}
{\bf Figure 19 - Four Term Relation For Type Three Invariant}
\end{center}
\vvv

Now one more exercise: Consider any Vassiliev invariant $v$ and let's determine
its value on the trefoil graph as in Figure 20.
\vv

\centerline{\includegraphics[scale=0.9]{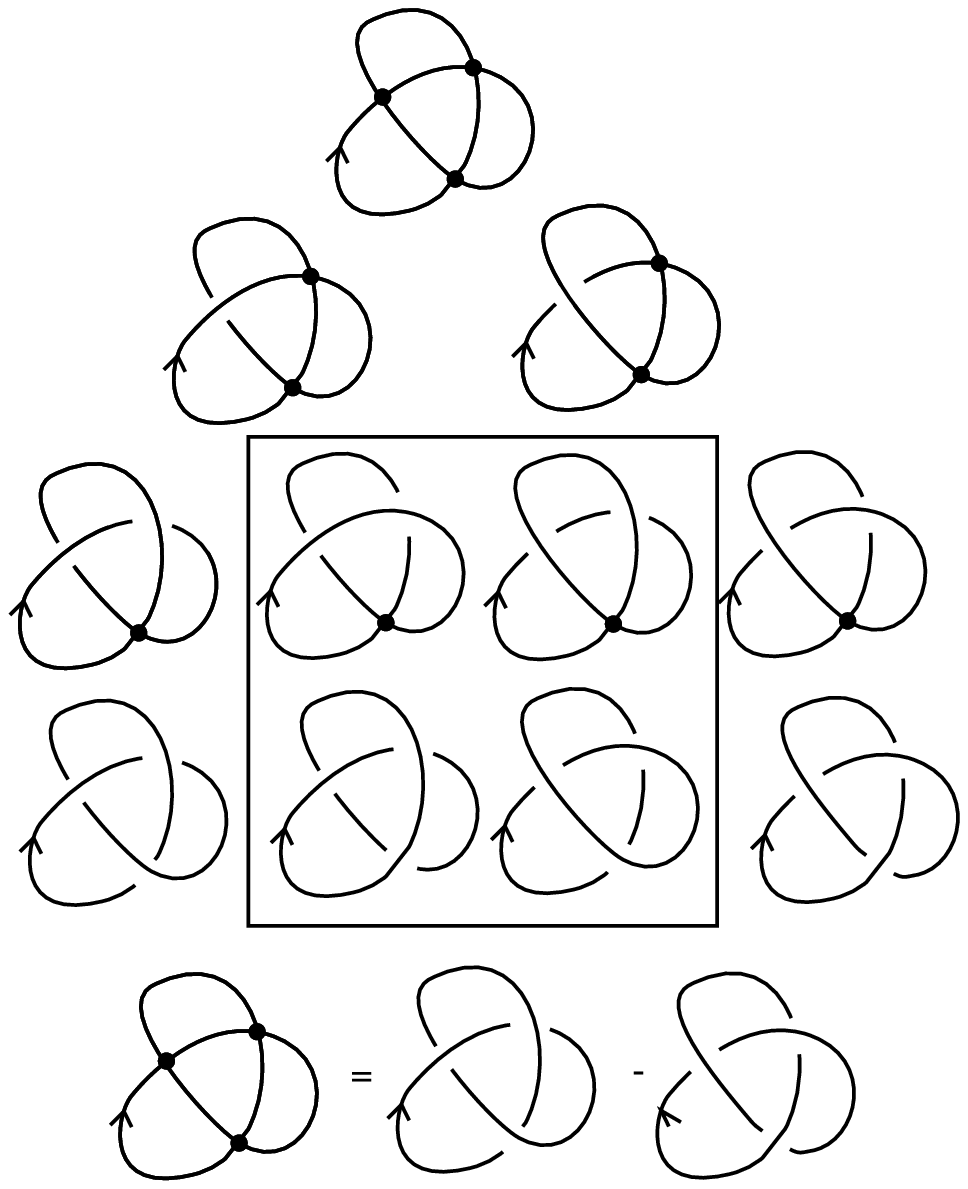}}
\vv

\begin{center}
{\bf Figure 20 - Trefoil Graph}
\end{center}
\newpage

The value of this invariant on the trefoil graph is equal to the difference
between its values on the trefoil knot and its mirror image. Therefore any
Vassiliev invariant that assigns a non-zero value to the trefoil graph can tell
the difference between the trefoil knot and its mirror image. \vspace{3mm}

\noindent {\bf Example.} This example shows how the original Jones polynomial is
composed of Vassiliev invariants of finite type. Let $V_{K}(t)$ denote the
original Jones polynomial \cite{Jones}. Recall the  oriented state expansion for
the Jones polynomial \cite{KNOTS} with the basic formulas ($\delta$ is the loop
value.)

$$V_{K_{+}} = -t^{1/2}V_{K_{0}} -tV_{K_{\infty}}$$ $$V_{K_{-}} =
-t^{-1/2}V_{K_{0}} -t^{-1}V_{K_{\infty}}.$$ $$\delta = -(t^{1/2} + t^{-1/2}).$$

Let $t=e^{x}.$ Then

$$V_{K_{+}} = -e^{x/2}V_{K_{0}} -e^{x}V_{K_{\infty}}$$ $$V_{K_{-}} =
-e^{-x/2}V_{K_{0}} -e^{-x}V_{K_{\infty}}.$$ $$\delta = -(e^{x/2} + e^{-x/2}).$$

Thus

$$V_{K_{*}} = V_{K_{+}}  -  V_{K_{-}} = -2sinh(x/2)V_{K_{0}}
-2sinh(x)V_{K_{\infty}}.$$

Thus  $x$ divides $V_{K_{*}}$, and therefore  $x^{k}$ divides $V_{G}$ whenever
$G$ is a graph with at least $k$ nodes.  Letting $$V_{G}(e^{x}) = \sum_{k=0}^{
\infty} v_{k}(G)x^{k},$$ we see that this condition implies that $v_{k}(G)$
vanishes whenever $G$ has more than $k$ nodes.  Hence {\em the coefficients of
the powers of $x$ in the expansion of $V_{K}(e^{x})$ are Vassiliev invariants of
finite type!} This result was first observed by Birman and Lin \cite{Birman and
Lin} by a  different argument. \vspace{3mm}

Let's look a little deeper and see the structure of the top row for the Vassiliev
invariants related to the Jones polynomial. By our previous remarks the top row
evaluations correspond to the leading terms in the power series expansion. Since
$$ \delta = -(e^{x/2}+e^{-x/2}) = -2 +[higher],$$ $$-e^{x/2}+e^{-x/2} = -x +
[higher],$$ $$-e^{x}+e^{-x} = -2x + [higher],$$ it follows that the top rows for
the Jones polynomial are computed by the recursion formulas $$v(K_{*}) =
-v(K_{0}) -2V(K_{\infty})$$ $$v([loop])=-2.$$ \vspace{3mm}

The reader can easily check that this recursion formula for the top rows of the
Jones polynomial implies that $v_{3}$ takes the value $24$ on the trefoil graph
and hence it is the Vassiliev invariant of type $3$ in the Jones polynomial that
first detects the difference between the trefoil knot and its mirror image.

This example gives a good picture of the general phenomenon of how the Vassiliev
invariants become building blocks for other invariants. In the case of the Jones
polynomial, we already know how to construct the invariant and so it is possible
to get a lot of information about these particular Vassiliev invariants by
looking directly at the Jones polynomial. This, in turn, gives insight into the
structure of the Jones polynomial itself. \vspace{3mm}

\subsection{Lie Algebra Weights}

Consider the  diagrammatic relation shown in Figure 21 . Call it (after Bar-Natan
\cite{Bar-Natan}) the $STU$ relation. \vspace{3mm}

\centerline{\includegraphics[scale=1.0]{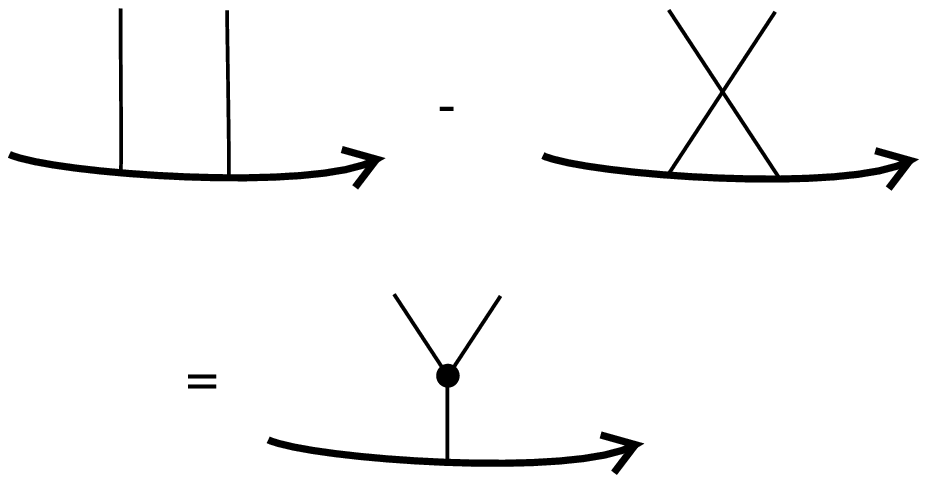}}
\vvvv

\begin{center}
{\bf Figure 21 - The $STU$ Relation}
\end{center}

\noindent {\bf Lemma.} $STU$ implies the 4-term relation. \vspace{1mm}

\noindent {\bf Proof.} View Figure 22. \vspace{3mm}

\centerline{\includegraphics[scale=1.2]{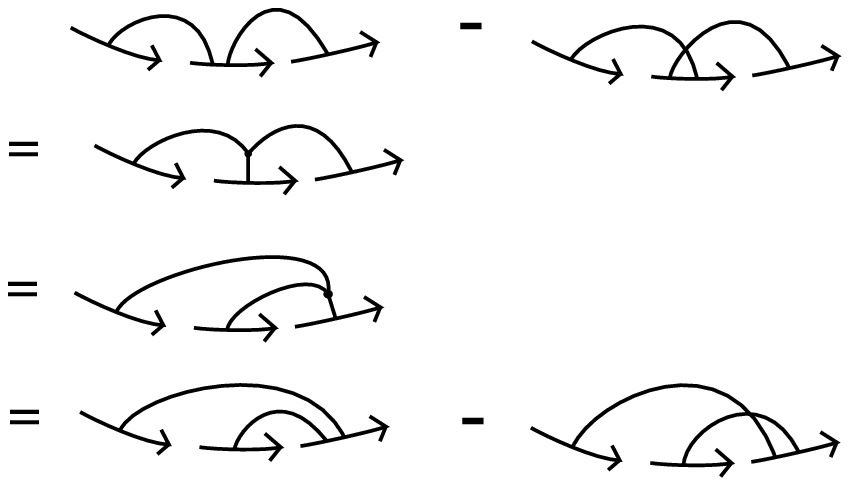}}
\vvvv

\begin{center}
{\bf Figure 22 - A Diagrammatic Proof}
\end{center}
\vvvv\vvv

$STU$ is the smile of the Cheshire cat. That smile generalizes the idea of a
Lie algebra. Take a (matrix) Lie algebra with generators $T^{a}$. Then $$T^{a}T^{b} -
T^{b}T^{a} = if_{abc}T^{c}$$ expresses the closure of the Lie algebra under
commutators. Translate this equation into diagrams as shown in Figure 23, and see
that this translation is $STU$ with Lie algebraic clothing!
\vvvv

\centerline{\includegraphics[scale=1.0]{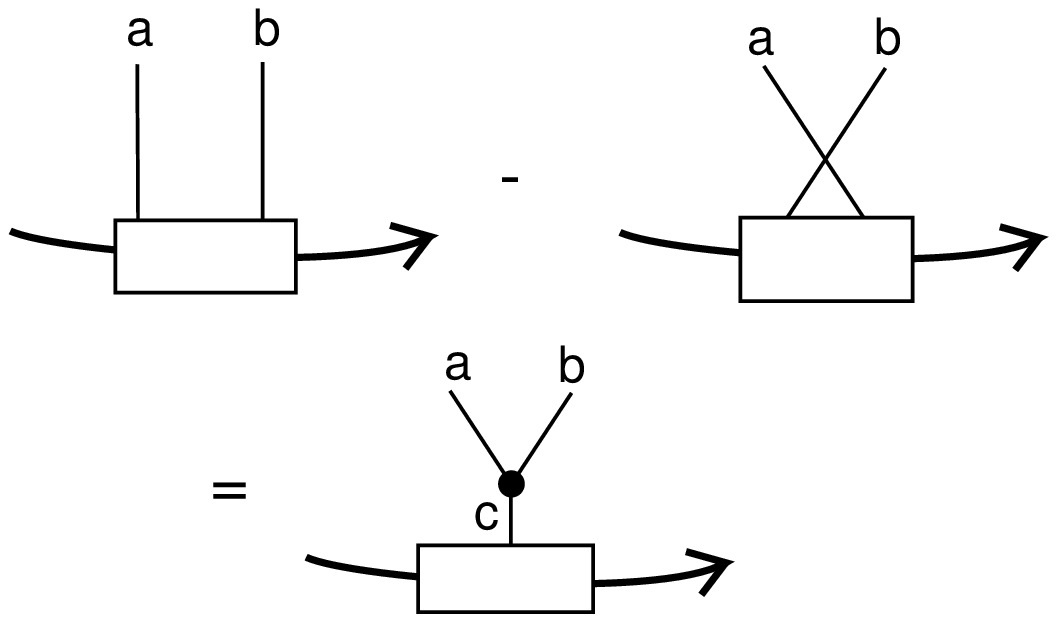}}
\vvvv

\begin{center}
{\bf Figure 23 - Algebraic Clothing}
\end{center}
\newpage

Here the structure tensor of the Lie algebra has been assumed (for simplicity) to
be invariant under cyclic permutation of the indices. This invariance means that
our last Lemma applies to this Lie algebraic interpretation of $STU$. The upshot
is that we can manufacture weight systems for graphs that satisfy the 4-term
relation by replacing paired points on the chord diagram by an insertion of
$T^{a}$ in one point of the pair and a corresponding insertion of $T^{a}$ at the
other point in the pair and summing over all $a$. The result of all such
insertions on a given chord diagram is a big sum of specific matrix products
along the circle of the diagram, each of which (being a circular product) is
interpreted as a trace. \vspace{3mm}

Let's say this last matter more precisely: Regard a graph with $k$ nodes as obtained
by identifying $k$ {\em pairs} of points on a circle.  Thus a code such as $1212$
taken in cyclic order specifies such a graph by regarding the points $1,2,1,2$ as
arrayed along a circle with the first and second $1$'s and $2$'s identified to form
the graph. Define,for a code $a_{1}a_{2}...a_{m}$
\vvv

$$wt(a_{1}a_{2}...a_{m}) =
trace(T^{a_{1}}T^{a_{2}}T^{a_{1}}...T^{a_{m}})$$
\vvvv

where the Einstein summation
convention is in place for the double appearances of indices on the right-hand
side. This gives the weight system. \vspace{3mm}

The weight system described by the above procedure satisfies the 4-term relation,
but does not necessarily satisfy the vanishing condition for isolated pairings.
This is because the framing compensation for converting an invariant of regular
isotopy to ambient isotopy has not yet been introduced. We will show how to do
this in the course of the discussion in the next paragraph. The main point to
make here is that by starting with the idea of extending an invariant of knots to
a Vassiliev invariant of embedded graphs and searching out the conditions on
graph evaluation demanded by the topology, we have inevitably entered the domain
of relations between Lie algebras and link invariants. Since the $STU$ relation does not
demand Lie algebras for its satisfaction we see that the landscape is wider than
the Lie algebra  context, but it is not yet understood how big is the class of
link invariants derived from Lie algebras. \vspace{3mm}

In fact, we can line up this weight system with the formalism related to the 
knot diagram by writing the Lie algebra insertions back on the 4-valent graph. We
then get a Casimir insertion at the node.  

To get the framing compensation, note that an isolated pairing corresponds to the
trace of the Casimir. Let $\gamma$ denote this trace.  See Figure 24.
\vvv

$$\gamma =
tr(\sum_{a} T^{a}T^{a})$$
\newpage

\centerline{\includegraphics[scale=1.2]{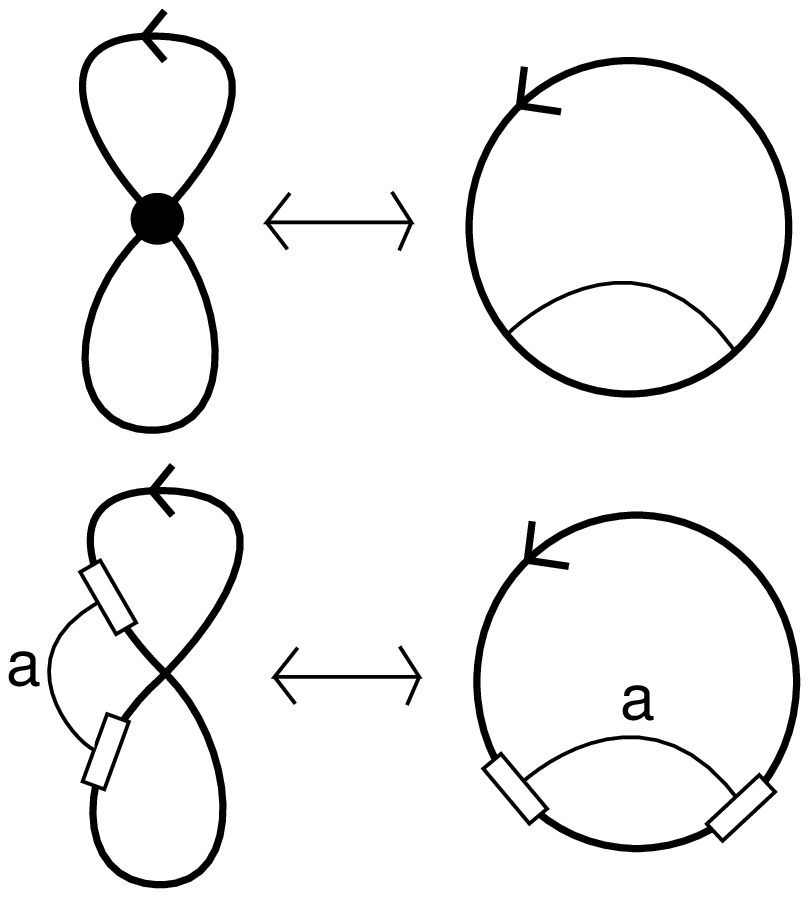}}
\vvvv\vvv

\begin{center}
{\bf  Figure 24 - Weight System and Casimir Insertion}
\end{center}
\vvvv\vvv

Let $D$ be the trace of the identity. Then it is easy to see that we must
compensate the given weight system by subtracting $(\gamma/D)$ multiplied by the
result of dropping the identification of the two given points. We can diagram
this by drawing two crossed arcs without a node drawn to bind them. Then the
modified recursion formula becomes as shown in Figure 25.
\vvv

For example, in the case of $SU(N)$ we have $D=N$, $\gamma = (N^{2} -1)/2$ so
that we get the transformation shown in Figure 25, including the use of the
Fierz identity.
\vvv

For $N=2$ the final formula of Figure 25  is,up to a multiple, exactly the top
row formula that we deduced for the Jones polynomial from its combinatorial
structure.
\newpage

\centerline{\includegraphics[scale=1.2]{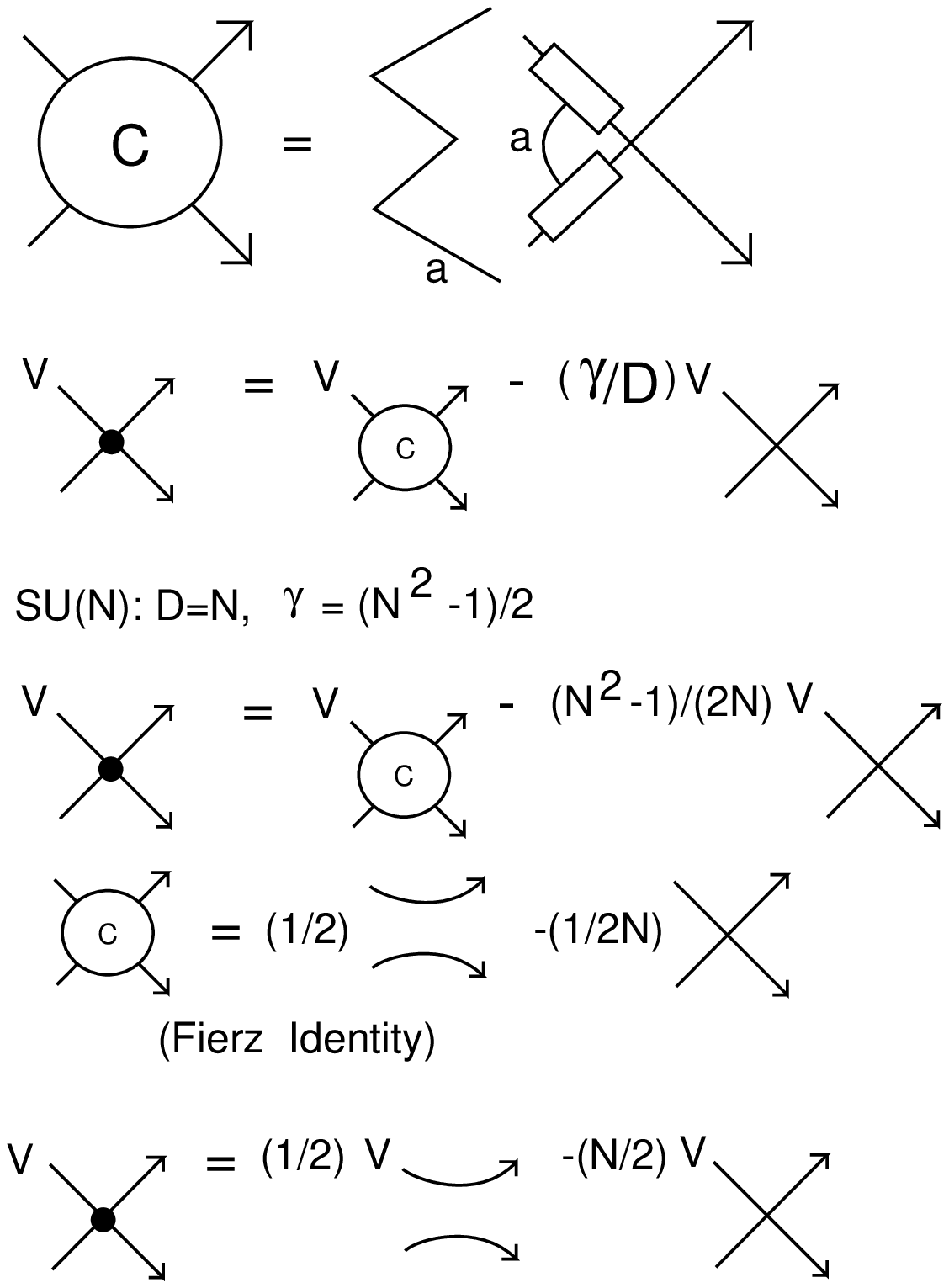}}
\vvvv\vvvv

\begin{center}
{\bf Figure 25 - Modified Recursion Formula}
\end{center}
\newpage

\subsection{Vassiliev Invariants and Witten's Functional Integral}

In \cite{Witten}  Edward Witten proposed a formulation of a class of 3-manifold
invariants as generalized Feynman integrals taking the form  $Z(M)$  where

$$Z(M) = \int dAexp[(ik/4\pi)S(M,A)].$$

Here  $M$ denotes a 3-manifold without boundary and $A$ is a gauge field  (also
called a gauge potential or gauge connection)  defined on $M$.  The gauge field
is a one-form on a trivial $G$-bundle over $M$ with values in a representation of
the  Lie algebra of $G$. The group $G$ corresponding to this Lie algebra is said to
be the gauge group. In this integral the ``action"   $S(M,A)$  is taken to be the
integral over $M$ of the trace of the Chern-Simons three-form    $CS = AdA +
(2/3)AAA$.  (The product is the wedge product of differential forms.)
\vspace{3mm}

$Z(M)$  integrates over all gauge fields modulo gauge equivalence (See
\cite{Atiyah:YM} for a discussion of the definition and meaning of gauge
equivalence.) \vspace{3mm}

The formalism  and   internal logic of Witten's integral supports  the existence
of a large class of topological invariants of 3-manifolds and  associated
invariants of knots and links in these manifolds. \vspace{3mm}

The invariants associated with this integral have been given rigorous
combinatorial descriptions
\cite{RT}, \cite{Turaev-Wenzl}, \cite{Kirby-Melvin}, \cite{Lickorish},
\cite{Walker}, \cite{TL}, but questions and conjectures arising from the integral
formulation are still outstanding. (See for example \cite{Atiyah},
\cite{Garoufalidis}, \cite{Gompf&Freed}, \cite{Jeffrey}, \cite{Rozansky}.)
Specific conjectures about this integral take the form of just how it involves
invariants of links and 3-manifolds, and how these invariants behave in certain
limits of the coupling constant $k$ in the integral. Many conjectures of this
sort can be verified through the combinatorial models. On the other hand, the
really outstanding conjecture about the integral is that it exists! At the
present time there is no measure theory or generalization of measure theory that
supports it. It is a fascinating exercise to take the speculation seriously,
suppose that it does really work like an integral and explore the formal
consequences.   Here is a formal structure of great beauty. It is also a
structure whose consequences can be verified by a remarkable variety of
alternative means. Perhaps in the course of the exploration there will appear a
hint of the true nature of this form of integration. \vspace{3mm}

We now  look at the formalism of the Witten integral in more detail and see how
it involves invariants of knots and links corresponding to each classical Lie
algebra.   In order to accomplish this task, we need to introduce the Wilson
loop.  The Wilson loop is an exponentiated version of integrating the gauge field
along a loop  $K$  in three-space that we take to be an embedding (knot) or a
curve with transversal self-intersections.  For this discussion, the Wilson loop
will be denoted by the notation  $ W_{K}(A) = <K|A>$ to denote the dependence on
the loop $K$ and the field $A$.   It is usually indicated by the symbolism  
$tr(Pexp(\int_{K} A))$ .   Thus $$W_{K}(A) = <K|A>  = tr(Pexp(\int_{K} A)).$$   
Here the $P$  denotes  path ordered integration - we are integrating and
exponentiating matrix valued functions, and so must keep track of the order of
the operations.  The  symbol  $tr$  denotes the trace of the resulting matrix.
\vspace{3mm}

With the help of the Wilson loop functional on knots and links,  Witten writes
down a functional integral for link invariants in a 3-manifold  $M$:

$$Z(M,K) = \int dAexp[(ik/4 \pi)S(M,A)] tr(Pexp(\int_{K} A)) $$

$$= \int dAexp[(ik/4 \pi)S] <K|A>.$$

Here $S(M,A)$ is the Chern-Simons Lagrangian, as in the previous discussion.

We abbreviate  $S(M,A)$  as $S$ and write  $<K|A>$  for the Wilson loop. Unless
otherwise mentioned, the manifold  $M$  will be the three-dimensional sphere 
$S^{3}$ \vspace{3mm}

An analysis of the formalism of this functional integral reveals quite a bit
about its role in knot theory.   This analysis depends upon key facts relating
the curvature of the gauge field to both the Wilson loop and the Chern-Simons
Lagrangian. The idea for using the curvature in this way is due to Lee Smolin
\cite{Smolin1, Smolin2} (See also \cite{Ramusino}). To this end, let us recall the local
coordinate structure of the gauge field  $A(x)$,  where  $x$  is a point in
three-space.   We can write   $A(x)  =  A_{a}^{k}(x)T^{a}dx_{k}$  where  the
index  $a$ ranges from $1$ to $m$ with the Lie algebra basis $\{T^{1}, T^{2},
T^{3}, ..., T^{m}\}$.  The index $k$   goes from $1$  to  $3$.     For each
choice of $a$  and  $k$,  $A_{a}^{k}(x)$   is a smooth function defined on
three-space. In  $A(x)$  we sum over the values of repeated indices.  The Lie
algebra generators $T^{a}$  are  matrices  corresponding to a given
representation of the Lie algebra of the gauge group $G.$   We assume some
properties of these matrices as follows: \vspace{3mm}

\noindent 1.  $[T^{a} , T^{b}] = i f_{abc}T^{c}$  where  $[x ,y] = xy - yx$ , and
$f_{abc}$ (the matrix of structure constants)  is totally antisymmetric.  There
is summation over repeated indices. \vspace{3mm}

\noindent 2.  $tr(T^{a}T^{b}) = \delta^{ab}/2$ where  $\delta^{ab}$ is the
Kronecker delta  ($\delta^{ab} = 1$ if $a=b$ and zero otherwise). \vspace{6mm}

We also assume some facts about curvature. (The reader may enjoy comparing with
the exposition in \cite{KNOTS}.  But note the difference of conventions on the
use of   $i$ in the Wilson loops and curvature definitions.)   The first fact  is
the relation of Wilson loops and curvature for small loops: \vspace{3mm}

\noindent {\bf Fact 1.} The result of evaluating a Wilson loop about a very small
planar circle around a point $x$ is proportional to the area enclosed by this
circle times the corresponding value of the curvature tensor of the gauge field
evaluated at $x$. The curvature tensor is  written 
$$F_{a}^{rs}(x)T^{a}dx_{r}dy_{s}.$$ It is the local coordinate expression of 
$AdA +AA.$ \vspace{3mm}

\noindent {\bf Application of Fact 1.}  Consider a given Wilson line  $<K|S>$.
Ask how its value will change if it is deformed infinitesimally in the
neighborhood of a point $x$ on the line.  Approximate the change according to
Fact 1, and regard the point $x$ as the place of curvature evaluation.  Let 
$\delta<K|A>$  denote the change in the value of the line.    $\delta <K|A>$  is
given by the formula $$\delta <K|A> = dx_{r}dx_{s}F_{a}^{rs}(x)T^{a}<K|A>.$$ This
 is the first order approximation to the change in the Wilson line. \vspace{3mm}

In this formula it  is understood that the Lie algebra matrices  $T^{a}$  are to
be inserted into the Wilson line at the point $x$,  and that we are summing over
repeated indices. This means that  each  $T^{a}<K|A>$ is  a new Wilson line
obtained from the original  line  $<K|A>$  by leaving the form of the loop
unchanged,  but inserting the matrix  $T^{a}$ into that loop at the point  $x$. 
A  Lie algebra generator is diagrammed by a little box with a single index line
and two input/output lines which correspond to its role as a matrix (hence as
mappings of a vector space to itself). See Figure 26. \vspace{3mm}

\centerline{\includegraphics[scale=1.0]{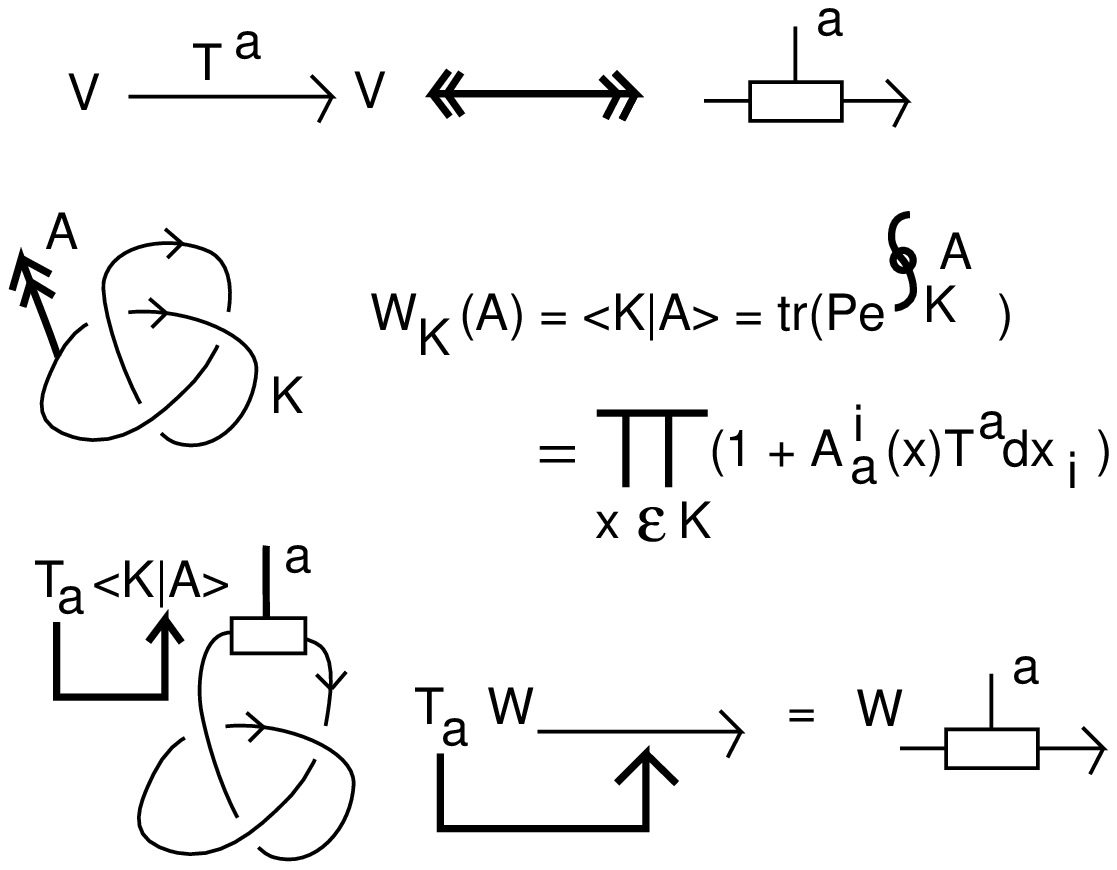}}
\vvvv

\begin{center}
{\bf Figure 26 - Wilson Loop Insertion}
\end{center}
\newpage

\noindent {\bf Remark.}  In thinking about the Wilson line $<K|A> =
tr(Pexp(\int_{K} A))$,  it is helpful to recall Euler's formula for the
exponential:

$$e^{x} = lim_{n \rightarrow \infty}(1+x/n)^{n}.$$
\vv

The  Wilson line  is  the limit, over partitions of the loop $K$,  of products of
the matrices  $(1 + A(x))$  where $x$ runs over the partition. Thus we can write
symbolically,
\vv

$$<K|A> =  \prod_{x \in K}(1 +A(x))  =  \prod_{x \in K}(1 +
A_{a}^{k}(x)T^{a}dx_{k}).$$
\vvv

It is understood that a product of matrices around a closed loop connotes the
trace of the product.  The ordering is forced by the one-dimensional nature of
the loop.   Insertion of a given matrix into this product at a point on the loop
is then a well-defined concept.   If  $T$  is a given matrix then it is
understood that   $T<K|A>$  denotes the insertion of $T$ into some point of the
loop. In the case above, it is understood from context in the formula
$$dx_{r}dx_{s}F_{a}^{rs}(x)T^{a}<K|A>$$ that the insertion is to be performed at
the point $x$  indicated in the argument of the curvature. \vspace{3mm}
\vvv

\noindent {\bf Remark.}  The   previous remark implies the following formula for
the variation of the Wilson loop with respect to the gauge field:

$$\delta <K|A>/\delta (A_{a}^{k}(x))  =  dx_{k}T^{a}<K|A>.$$
\vvv

Varying the Wilson loop with respect to the gauge field results in the insertion
of an infinitesimal Lie algebra element into the loop.
\vvvv

\noindent {\bf Proof.}
\v

$$\delta <K|A>/\delta (A_{a}^{k}(x))$$

$$= \delta \prod_{y \in K}(1 + A_{a}^{k}(y)T^{a}dy_{k})/\delta (A_{a}^{k}(x))$$

$$= \prod_{y<x \in K}(1 + A_{a}^{k}(y)T^{a}dy_{k}) [T^{a}dx_{k}] \prod_{y>x \in
K}(1 + A_{a}^{k}(y)T^{a}dy_{k})$$

$$= dx_{k}T^{a}<K|A>.$$
\newpage

\noindent{\bf Fact 2.}  The variation of the Chern-Simons Lagrangian  $S$  with
respect to the gauge potential at a given point in three-space is related to the
values of the curvature tensor at that point by the following formula:

$$F_{a}^{rs}(x)  =  \epsilon_{rst} \delta S/\delta (A_{a}^{t}(x)).$$
\vvv

Here 
$\epsilon_{abc}$ is the epsilon symbol for three indices, i.e. it is $+1$ for
positive permutations of $123$ and $-1$ for negative permutations of $123$ and
zero if any two indices are repeated.
\vvv

With these facts at hand we are prepared to determine how the Witten integral
behaves under a small deformation of the loop $K.$
\vvv

In accord with the theme of this paper, we shall use a system of abstract tensor
diagrams to look at the differential algebra related to the functional integral. 
The translation to diagrams is accomplished with the aid of Figure 27 and Figure
28.  In Figure 27 we give diagrammatic equivalents for the component parts of our
machinery. Tensors become labelled boxes. Indices become lines emanating from the
boxes. Repeated indices that we intend to sum over become lines from one box to
another. (The eye can immediately apprehend the repeated indices and the tensors
where they are repeated.)  Note that we use a capital $D$ with lines extending
from the top and the bottom for the partial derivative with respect to the gauge
field, a capital $W$ with a link diagrammatic subscript for the Wilson loop,
a cubic vertex for the three-index epsilon, little triangles with emanating arcs
for the differentials of the space variables.
\vvv

The Lie algebra generators are little boxes with single index lines and two
input/output lines which correspond to their roles as matrices (hence as mappings
of a vector space to itself). The Lie algebra generators are, in all cases of our
calculation, inserted into the Wilson line either through the curvature tensor or
through insertions related to differentiating the Wilson line.
\vvv

In Figure 28 we give the diagrammatic calculation of the change of the functional
integral corresponding to a tiny change in the Wilson loop. The result is a
double insertion of Lie Algebra generators into the line, coupled with the
presence of a volume form that will vanish if the deformation does not twist in
three independent directions. This shows that the functional integral is formally
invariant under regular isotopy since the regular isotopy moves are changes in
the Wilson line that happen entirely in a plane. One does not expect the integral
to be invariant under a Reidemeister move of type one, and it is not. This
framing compensation can be determined by the methods that we are discussing
\cite{Func}, but we will not go into the details of those calculations here.
\newpage

\centerline{\includegraphics[scale=1.2]{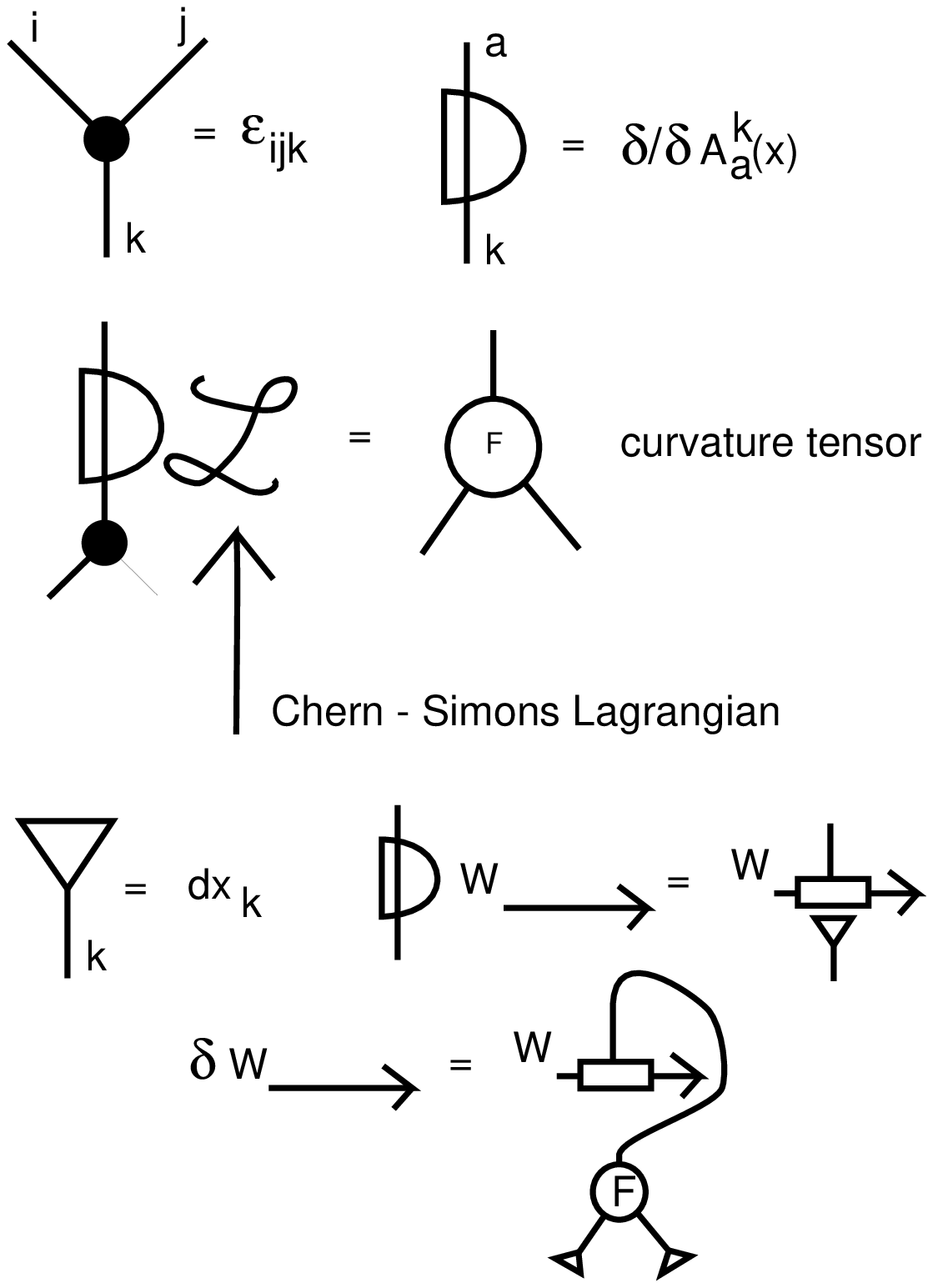}}
\vvvv\vvvv

\begin{center}
{\bf Figure 27 - Notation}
\end{center}
\newpage

\centerline{\includegraphics[scale=1.0]{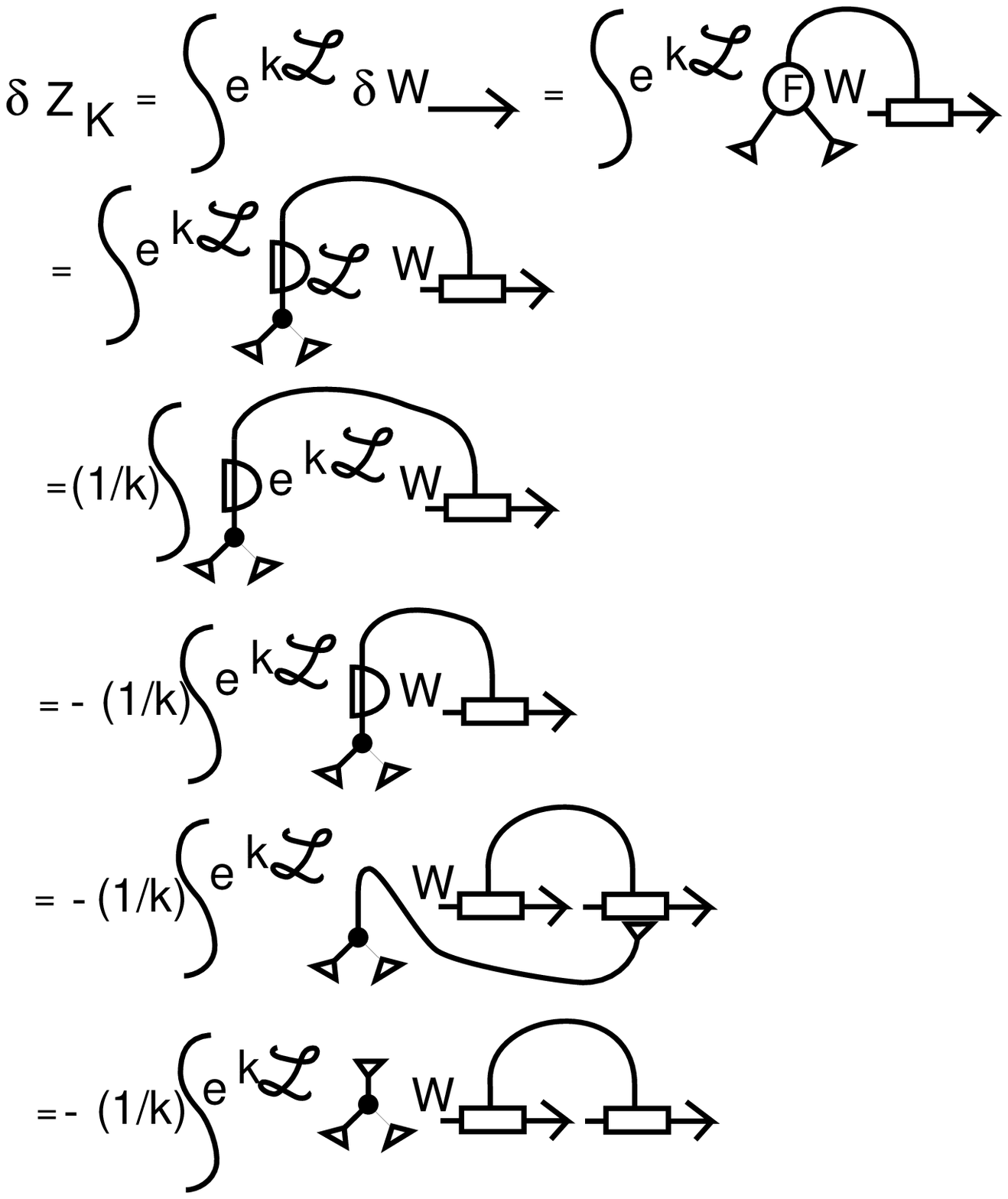}}
\vvv

\begin{center}
{\bf Figure 28 - Derivation}
\end{center}
\vvv

In Figure 29 we show the application of the calculation in Figure 28 to the case
of switching a crossing. The same formula applies, with a different
interpretation,  to the case where  $x$ is  a double point of transversal self-intersection
of a loop $K$,  and the deformation consists in
shifting one of the crossing segments perpendicularly to the plane of intersection  so that the
self-intersection point disappears.  In this  case,  one $T^{a}$  is inserted
into each of the transversal crossing segments so that $T^{a}T^{a}<K|A>$ denotes
a Wilson loop with a self-intersection  at  $x$   and insertions of $T^{a}$  at
$x + \epsilon_{1}$ and  $x + \epsilon_{2}$  where $\epsilon_{1}$ and
$\epsilon_{2}$ denote small displacements along the two arcs of $K$ that
intersect at $x.$  In this case, the volume form is nonzero, with two directions
coming from the plane of movement of one arc, and the perpendicular direction is
the direction of the other arc.  The reason for the insertion into the two lines
is a direct  consequence of the calculational form of Figure 28: The first insertion
is in the moving line, due to curvature. The second insertion is the consequence
of differentiating the self-touching Wilson line. Since this line can be regarded
as a product, the differentiation occurs twice at the point of intersection, and
it is the second direction that produces the non-vanishing volume form.
\vvvv\vvv

\centerline{\includegraphics[scale=1.0]{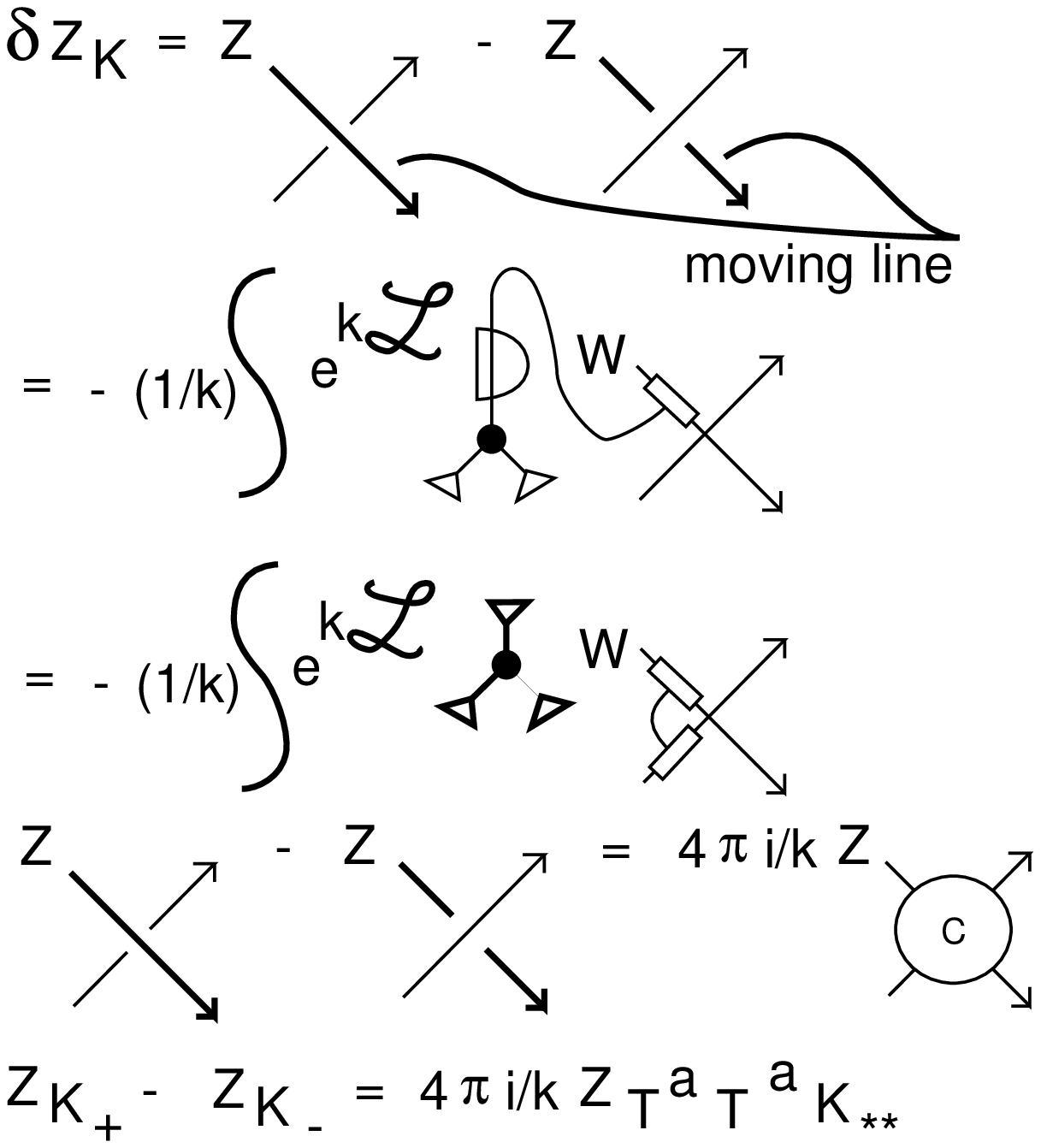}}
\vvvv\vv

\begin{center}
{\bf Figure 29 - Crossing Switch}
\end{center}
\newpage

Up to the choice of our conventions for constants, the switching formula is, as
shown in Figure 29,

$$Z(K_{+}) -  Z(K_{-}) =  (4 \pi i/k)\int dA  exp[(ik/4\pi)S]
T^{a}T^{a}<K_{**}|A>$$ $$= (4 \pi i/k) Z(T^{a}T^{a}K_{**}).$$

The key point is to notice that the Lie algebra insertion for this difference is
exactly what we did to make the weight systems for Vassiliev invariants (without
the framing compensation). Thus the formalism of the Witten functional integral
takes us directly to these weight systems in the case of the classical Lie
algebras. The functional integral  is central to the structure of the Vassiliev
invariants. \vspace{3mm}

\subsection{Combinatorial Constructions for Vassiliev Invariants} Perhaps the
most remarkable thing about this story of the structure of the Vassiliev
invariants is the way that Lie algebras are so naturally involved in the
structure of the weight systems. This shows the remarkably close nature of the
combinatorial structure of Lie algebras and the combinatorics of knots and links
via the Reidemeister moves.  A really complete story about the Vassiliev
invariants at this combinatorial level would produce their existence on the basis
of the weight systems with entirely elementary arguments. \vspace{3mm}

As we have already mentioned, one can prove that a given set of weights for the
top row, satisfying the abstract four-term relation does imply that there exists
a Vassiliev invariant of finite type $n$ realizing  these weights for graphs with
$n$ nodes.  Proofs of this result either use analysis \cite{Bar-Natan},
\cite{Altschuler-Friedel} or non-trivial algebra \cite{Cart}, \cite{Bar-Natan}.
There is no known elementary combinatorial proof of the existence of Vassiliev
invariants for given top rows. \vspace{3mm}

Of course quantum link invariants (See section 4 of these lectures.) do give
combinatorial constructions for large classes of link invariants.  These
constructions rest on solutions to the Yang-Baxter equations, and it is not known
how to describe the subset of finite type Vassiliev invariants that are so
produced. \vspace{3mm}

It is certainly helpful to look at the structure of Vassiliev invariants that
arise from already-defined knot invariants.  If $V(K)$ is an already defined
invariant of knots (and possibly links), then its extension to a Vassiliev
invariant is calculated on embedded graphs $G$ by expanding each graphical vertex
into a difference by resolving the vertex into a positive crossing and a negative
crossing. If we know that $V(K)$ is of finite type $n$ and $G$ has $n$ nodes then
we can take any embedding of $G$ that is convenient, and calculate $V(G)$ in
terms of all the knots that arise in resolving the nodes of this chosen
embedding. This is a finite collection of knots. Since there is a finite
collection of 4-valent graphs with $n$ nodes, it follows that the top row
evaluation for the invariant $V(K)$ is determined by the values of $V(K)$ on a
finite collection of knots. Instead of asking for the values of the Vassiliev
invariant on a top row, we can ask for this set of knots and the values of the
invariant on this set of knots. A minimal set of knots that can be used to
generate a given Vassiliev invariant will be called a {\em knots basis} for the
invariant. Thus we have shown that the set consisting of the unknot, the
right-handed trefoil and the left handed-trefoil is a knots basis for a Vassiliev
invariant of type 3. See \cite{Mathias} for more information about this point of
view.
\vv

A tantalizing combinatorial approach to Vassiliev invariants is due to Michael
Polyak and Oleg Viro \cite{Polyak-Viro}.  They give explicit formulas for the
second, third and fourth Vassiliev invariants and conjecture that their method
will work for Vassiliev invariants of all orders. The method is as follows.
\vspace{3mm}

First one makes a new representation for oriented knots by taking {\em Gauss
diagrams}. A Gauss diagram is a diagrammatic representation of the classical {\em
Gauss code} of the knot. The Gauss code is obtained from the oriented knot
diagram by first labelling each crossing with a naming label (such as
$1$,$2$,...) and also indicating the crossing type ($+1$ or $-1$).  Then choose a
basepoint on the knot diagram and begin walking along the diagram, recording the
name of the crossings encountered, their sign and whether the walk takes you over
or under that crossing.  For example, if you go under crossing $1$ whose sign is
$+$ then you will record $o+1$. Thus the Gauss code of the positive trefoil
diagram is $$(o1+)(u2+)(o3+)(u1+)(o2+)(u3+).$$ For prime knots the Gauss code is
sufficient information to reconstruct the knot diagram. See  \cite{GAUSS}
for a sketch of the proof of this result and for other references.
\vvv

\centerline{\includegraphics[scale=0.75]{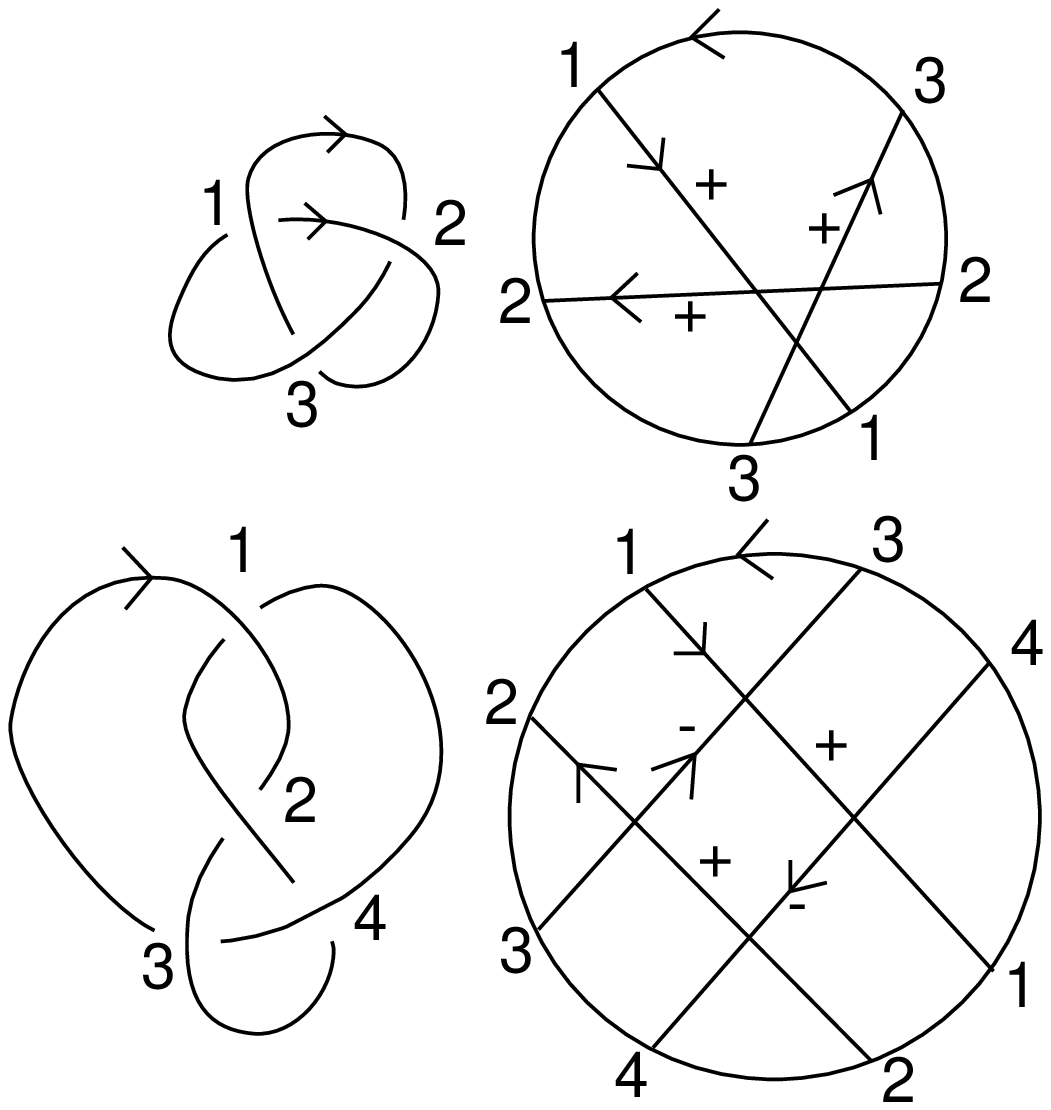}}
\vvvv

\begin{center}
{\bf Figure 30 - Gauss Diagrams}
\end{center}
\newpage
\vspace*{-.2in}

To form a Gauss diagram from a Gauss code, take an oriented circle with a
basepoint chosen on the circle. Walk along the circle marking it with the labels
for the crossings in the order of the Gauss code. Now draw chords between the
points on the circle that have the same label. Orient each chord from
overcrossing site to undercrossing site. Mark each chord with $+1$ or $-1$
according to the sign of the corresponding crossing in the Gauss code. The
resulting labelled and basepointed graph is the Gauss diagram for the knot.  See
Figure 30 for examples. \vspace{3mm}

The Gauss diagram is deliberately formulated to have the structure of a chord
diagram (as we have discussed for the weight systems for Vassiliev invariants).
If $G(K)$ is the Gauss diagram for a knot $K$, and $D$ is an oriented (i.e. the
chords as well as the circle in the diagram are oriented) chord diagram, let
$|G(K)|$ denote the number of chords in $G(K)$ and $|D|$ denote the number of
chords in $D.$ If $|D| \leq |G(K)|$ then we may consider oriented embeddings of
$D$ in $G(K)$. For a given embedding $i: D \longrightarrow G(K)$ define
$$<i(D) | G(K)> = sign(i)$$
where $sign(i)$ denotes the product of the signs of the chords
in $G(K) \cap i(D).$ Now suppose that $C$ is a collection of {\em oriented} chord
diagrams, each with $n$ chords, and that
$$eval: C \longrightarrow R$$
is an
evaluation mapping on these diagrams that satisfies the four-term relation at
level $n$. Then we can define
$$<D|K> = \sum_{i:D \longrightarrow G(K)} <i(D)|G(K)> $$
and
$$v(K) = \sum_{D \in C} <D|K> eval(D).$$

For appropriate oriented chord subsets this definition can produce Vassiliev
invariants $v(K)$ of type $n$.  For example, in the case of the Vassiliev
invariant of type three taking value $0$ on the unknot and value $1$ on the
right-handed trefoil, $-1$ on the left-handed trefoil,  Polyak and Viro give the
specific formula
$$v_{3}(K) = <A|K> + (1/2)<B|K>$$

\noindent where $A$ denotes the trefoil chord diagram as we described it in
section 3 and $B$ denotes the three-chord diagram consisting of two parallel
chords pierced by a third chord. In Figure 31 we show the specific orientations
for the chord diagrams $A$ and $B$.  The key to this construction is in the
choice of orientations for the chord diagrams in $C = \{ A, B \}.$  It is a nice
exercise in translation of the  Reidemeister moves to Gauss diagrams to see that 
$v_{3}(K)$ is indeed a knot invariant. \vspace{3mm}

It is possible that all Vassiliev invariants can be constructed by a method
similar to the formula $v(K) = \sum_{D \in C} <D|K> eval(D).$ This remains to be
seen.
\newpage
\vspace*{-.3in}

\centerline{\includegraphics[scale=0.85]{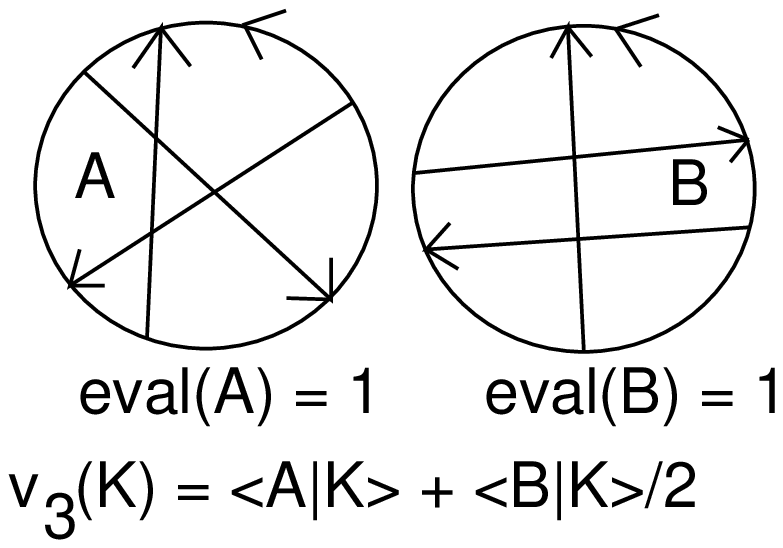}}
\vvvv

\begin{center}
{\bf Figure 31 - Oriented Chord Diagrams for $\mathbf{v_{3}}.$}
\end{center}

\subsection{$8_{17}$}

It is an open problem whether there are Vassiliev invariants that can detect the
difference between a knot and its reverse (The reverse of an oriented knot is
obtained by flipping the orientation.). The smallest instance of a non-invertible
knot is the knot $8_{17}$ depicted in Figure 32.

\centerline{\includegraphics[scale=0.8]{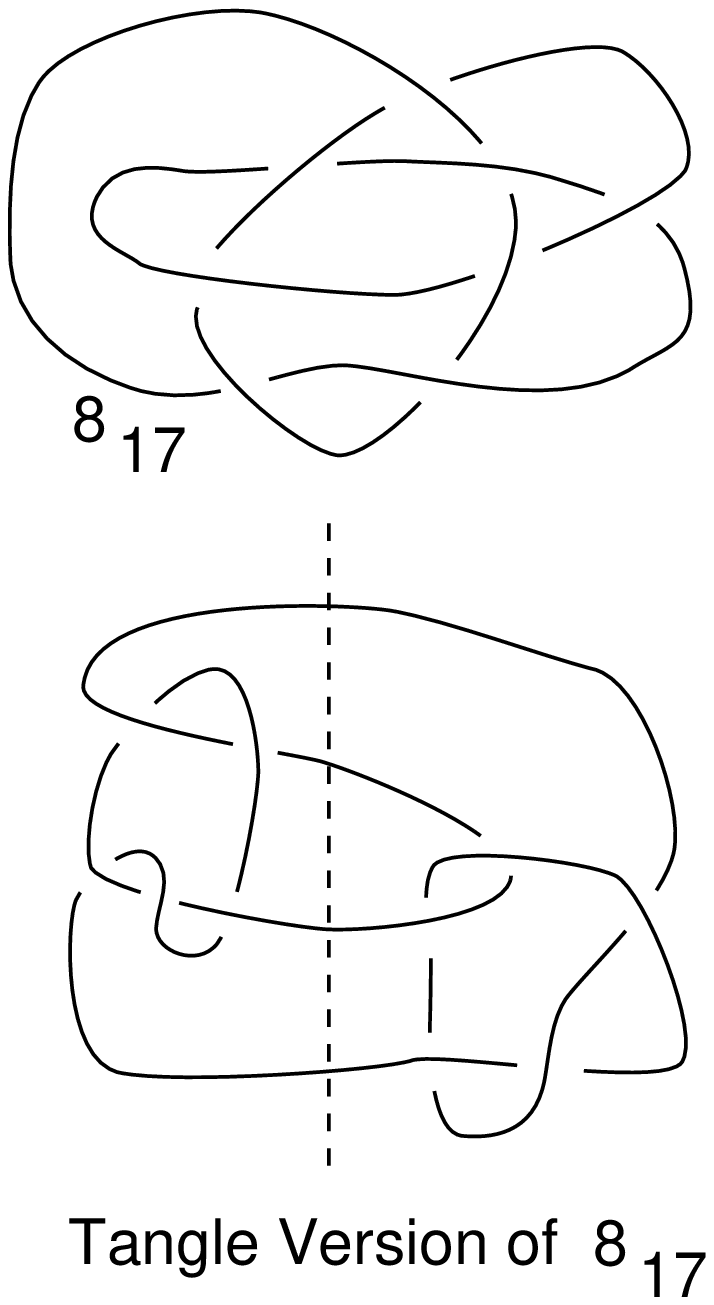}}
\vvvv

\begin{center}
{\bf Figure 32 - Tangle Decomposition of $8_{17}.$}
\end{center}

Thus, at the time of this
writing there is no known Vassiliev invariant that can detect the
non-invertibility of $8_{17}.$  On the other hand, the tangle decomposition shown
in Figure 32 can be used in conjunction with the results of Siebenmann and
Bonahon \cite{S and B} and the formulations of John Conway \cite{Con} to show
this non-invertibility.  These tangle decomposition methods use higher level
information about the diagrams than is easy to encode in Vassiliev invariants. 
The purpose of this section is to underline this discrepancy between different
levels in the combinatorial topology. \vspace{3mm}

\section{Quantum Link Invariants} In this section we describe the construction of
quantum link invariants from knot and link diagrams that are arranged with
respect to a given direction in the plane. This special direction will be called
``time". Arrangement with respect to the special direction means that
perpendiculars to this direction meet the diagram transversely (at edges or at
crossings)  or tangentially (at maxima and minima).  The designation of the
special direction as time allows the interpretation of the consequent evaluation
of the diagram as a generalized scattering amplitude. \vspace{3mm}

In the course of this discussion we find the need to reformulate the Reidemeister
moves for knot and link diagrams that are arranged to be transverse (except for a
finite collection of standard critical points) to the specific special direction
introduced in the previous paragraph. This brings us back to our theme of
diagrams and related structures. This particular reformulation of the
Reidemeister moves is quite far-reaching. It encompasses the relationship of link
invariants with solutions to the Yang-Baxter equation and the relationship with
Hopf algebras (to be dealt with in Section 5). \vspace{3mm}

\subsection{Knot Amplitudes} Consider first a circle in a spacetime plane with
time represented vertically and space horizontally as in Figure 33. \vspace{3mm}

The circle represents a vacuum to vacuum process that includes the creation of
two ``particles" and their subsequent annihilation. We could divide the circle
into these two parts (creation ``cup"  and annihilation ``cap") and consider the
amplitude   $<cap|cup>$. Since the diagram for the creation of the two particles
ends in two separate points, it is natural to take a vector space of the form  $V
\otimes V$  as the target for the bra and as the domain of the ket.  We imagine
at least one particle property being catalogued by each factor of the tensor 
product.  For example, a basis of  $V$  could enumerate the spins of the created
particles.

\centerline{\includegraphics[scale=1.0]{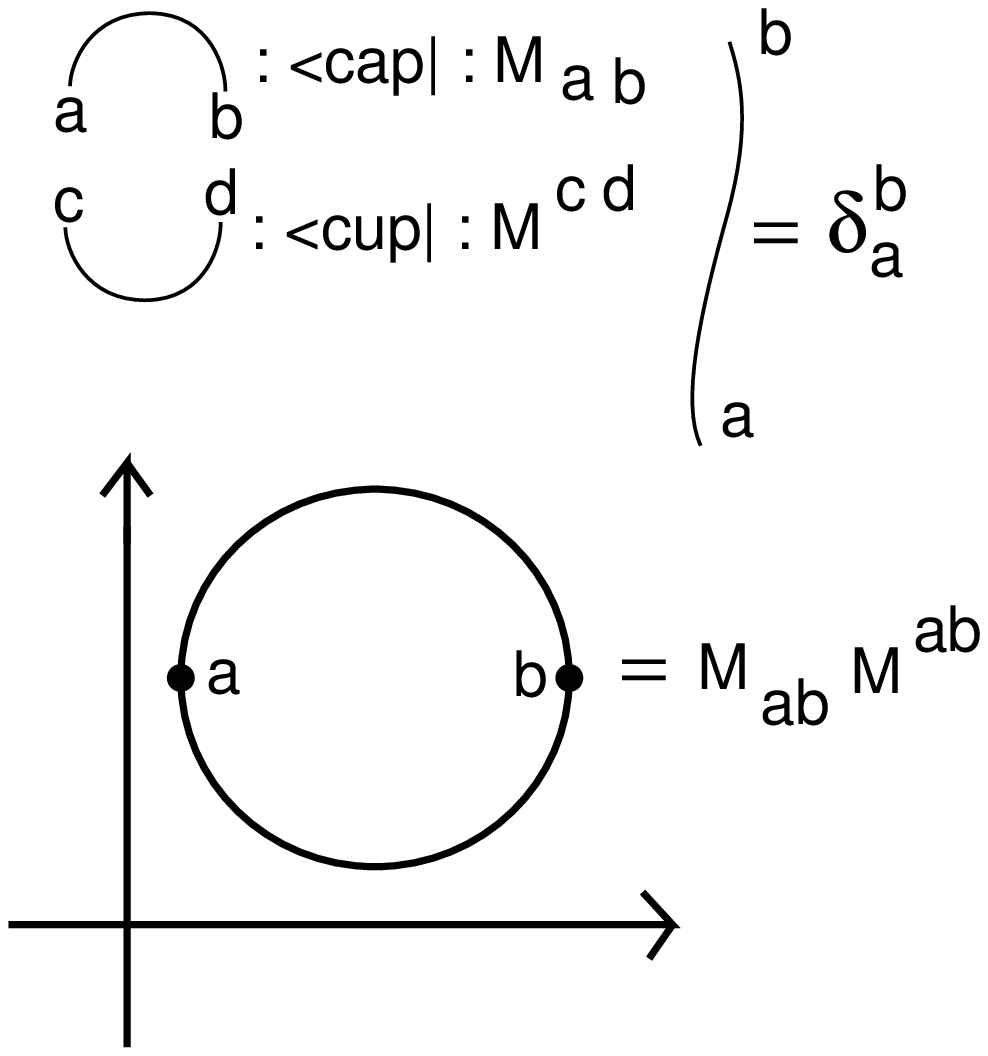}}
\vvv

\begin{center}
{\bf  Figure 33 - Spacetime Circle}
\end{center}
\vvvv

Any non-self-intersecting differentiable curve can be rigidly rotated until it is
in general position with respect to the vertical.  It will then be seen to be
decomposed into an interconnection of minima and maxima.  We can evaluate an
amplitude for any curve in general position with respect to a vertical direction.
Any simple closed curve in the plane is isotopic to a circle, by the Jordan Curve
Theorem.  If these are topological amplitudes,  then the value for any simple
closed curve should be equal to the original amplitude for the circle. What
condition on creation (cup) and annihilation (cap) will insure topological
amplitudes?  The answer derives from the fact that isotopies of the simple closed
curves are generated by the cancellation of adjacent maxima and minima as
illustrated in Figure 34. \vspace{3mm}

This condition is articulated by taking a matrix representation for the
corresponding operators. Specifically,  let   $\{ e_{1}, e_{2}, ..., e_{n} \}$ 
be a basis for  $V.$ Let $e_{ab} = e_{a} \otimes  e_{b}$ denote the elements of
the tensor basis for  $V \otimes  V$.  Then there are matrices  $M_{ab}$   and
$M^{ab}$  such that $$|cup>(1)  =  \sum M^{ab}e_{ab}$$ with the summation taken
over all values of $a$ and $b$ from $1$ to $n.$ Similarly,  $<cap|$  is described
by    $$<cap|(e_{ab}) =  M_{ab}.$$ Thus the amplitude for the circle is
$$<cap|cup>(1)  =  <cap| \sum M^{ab}e_{ab}$$ $$= \sum M^{ab} <cap|(e_{ab}) = \sum
M^{ab}M_{ab}.$$ In general, the value of  the amplitude on a simple closed curve
is obtained by translating it into an ``abstract tensor expression"  in the
$M^{ab}$ and $M_{ab}$,  and  then summing over these products for all cases of
repeated indices. Note that here the value ``1" corresponds to the vacuum.

\centerline{\includegraphics[scale=0.9]{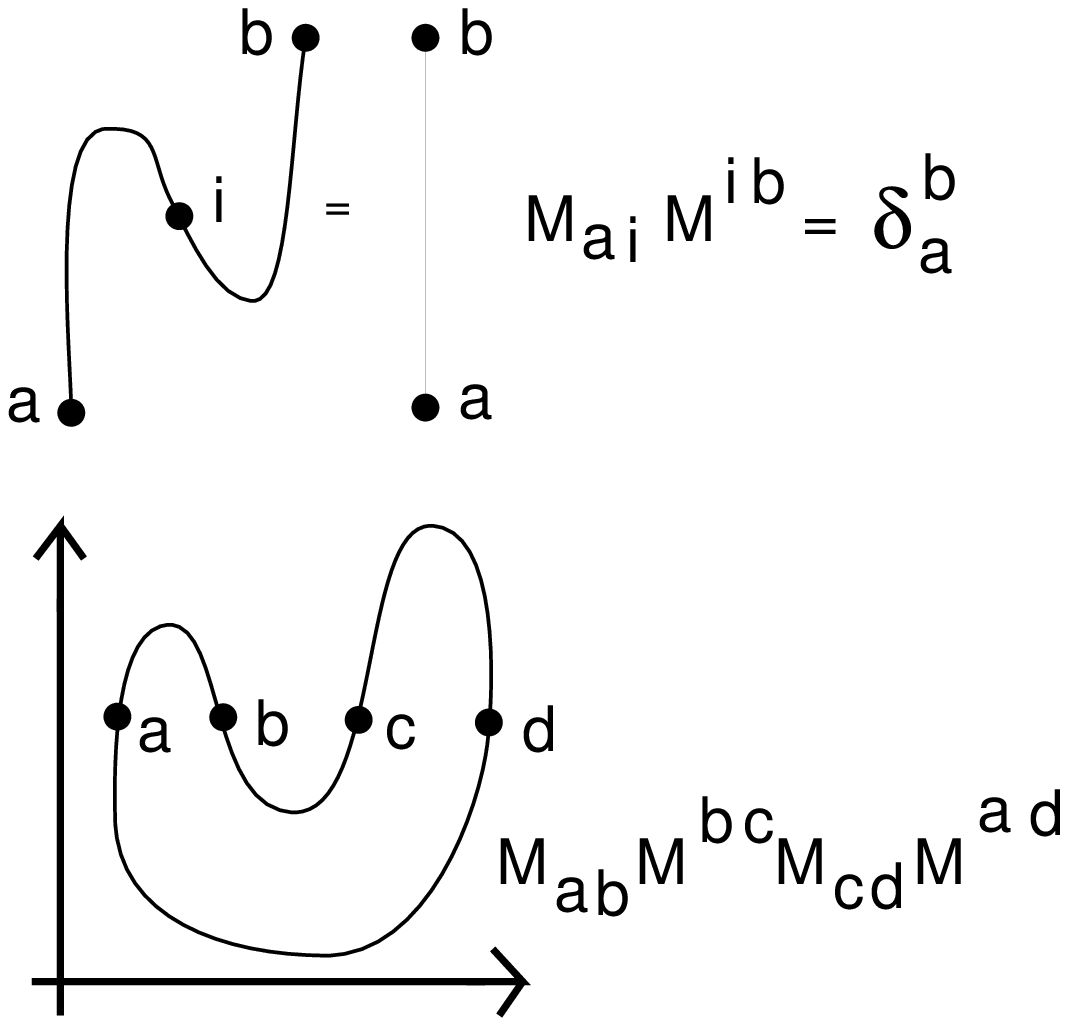}}
\vvvv

\begin{center}
{\bf Figure 34- Spacetime Jordan Curve}
\end{center}
\vvv

Returning to the topological conditions we see that they are just that the
matrices $M^{ab}$  and  $M_{cd}$  are inverses in the sense that $$\sum_{i}
M^{ai}M_{ib}  =\delta^{a}_{b}$$ where   $\delta^{a}_{b}$ denotes the identity
matrix. See Figure 34. \vspace{3mm}

One of the simplest choices is to take a $2$ x $2$ matrix $M$ such that $M^{2} =
I$ where $I$ is the identity matrix.  Then the entries of $M$ can be used for
both the cup and the cap.  The value for a loop is then equal to the sum of the
squares of the entries of $M$:

$$<cap|cup>(1) = \sum M^{ab}M_{ab} =\sum M_{ab}M_{ab}= \sum M_{ab}^{2}.$$

In particular, consider the following choice for $M.$  It has square equal to the
identity matrix and yields a loop value of $d = -A^{2} -A^{-2}$, just the right
loop value for the bracket polynomial model for the Jones polynomial \cite{stat},
\cite{state}.

$$M =  
\left[ \begin{array}{cc} 
0& iA  \\ 
-iA^{-1} & 0 
\end{array} 
\right] 
$$

Any knot or link can be represented by a picture that is configured with respect
to a vertical direction in the plane. The picture will decompose into minima
(creations), maxima (annihilations) and crossings of the two types shown in
Figure 35.  Here the knots and links are unoriented. These models generalize
easily to include orientation.
\vvvv

\centerline{\includegraphics[scale=1.0]{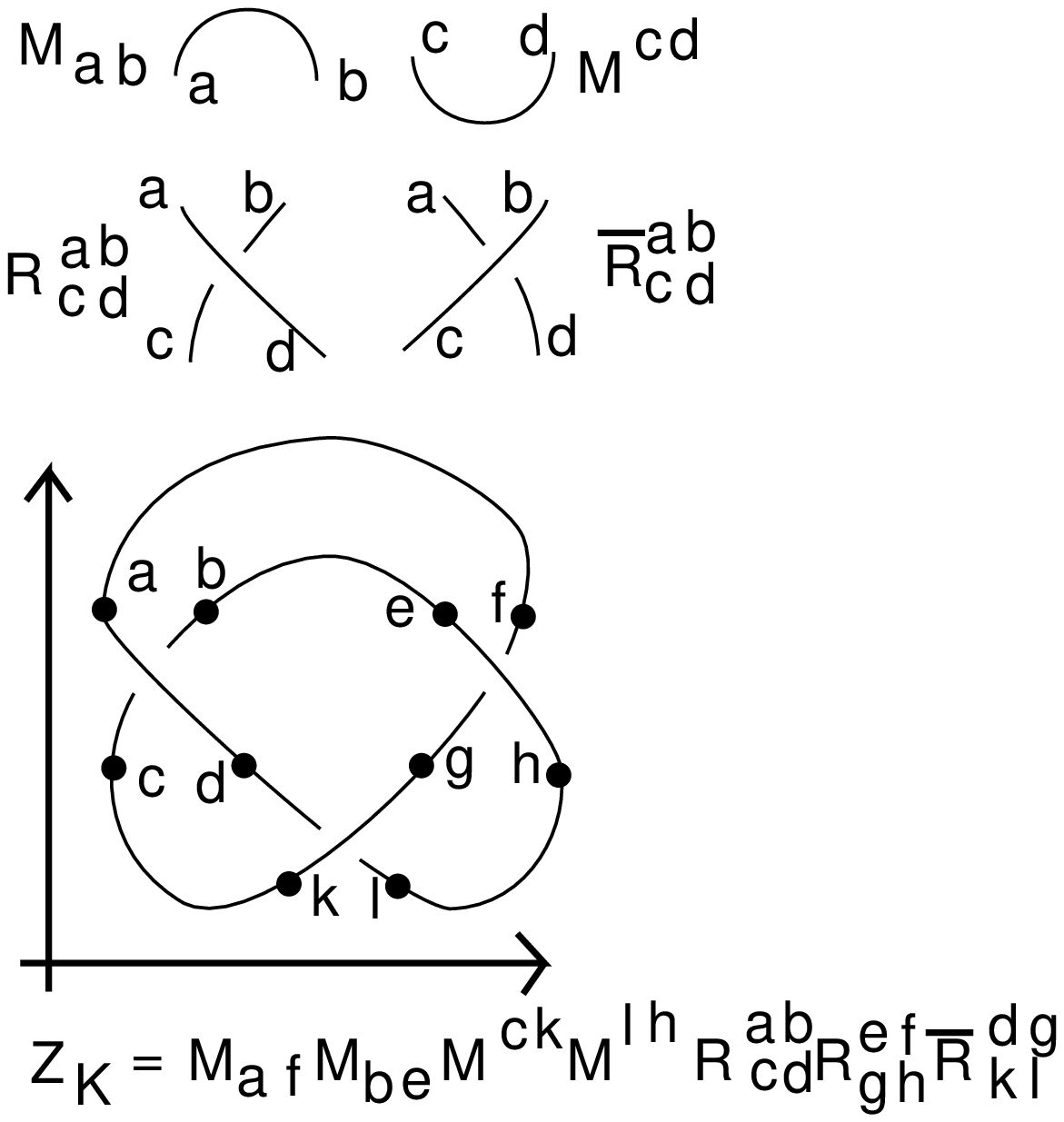}}
\vvvv

\begin{center}
{\bf Figure 35 - Cups, Caps and Crossings}
\end{center}
\newpage

Next to each of the crossings we have indicated mappings of $V \otimes V$  to
itself , called   $R$  and  $R^{-1}$  respectively.  These mappings represent the
transitions corresponding to elementary braiding.   We now have the vocabulary of
$cup$, $cap$,  $R$  and $R^{-1}$.   Any knot or link can be written as a
composition of these fragments, and consequently a choice of such mappings
determines an amplitude for knots and links.  In order for such an amplitude to
be topological  (i.e. an invariant of regular isotopy the equivalence relation
generated by the second and third of the classical Reidemeister moves)  we want
it to be invariant under a list of local moves on the diagrams as shown in Figure
36.  These moves are an augmented list of  Reidemeister moves, adjusted to take
care of the fact that the diagrams are arranged with respect to a given direction
in the plane.  The proof that these moves generate regular isotopy is composed in
exact parallel to the proof that we gave for the classical Reidemeister moves in
section 2. In the piecewise linear setting, maxima and minima are replaced by
upward and downward pointing angles. The fact that the triangle, in the
Reidemeister piecewise linear triangle move, must be projected so that it is
generically transverse to the vertical direction in the plane introduces the
extra restriction that expands the move set.
\vvv

In this context, the algebraic translation of  Move $III$  is the Yang-Baxter
equation  that occurred for the first time in problems of exactly solved models
in statistical mechanics \cite{Baxter}.

\centerline{\includegraphics[scale=1.0]{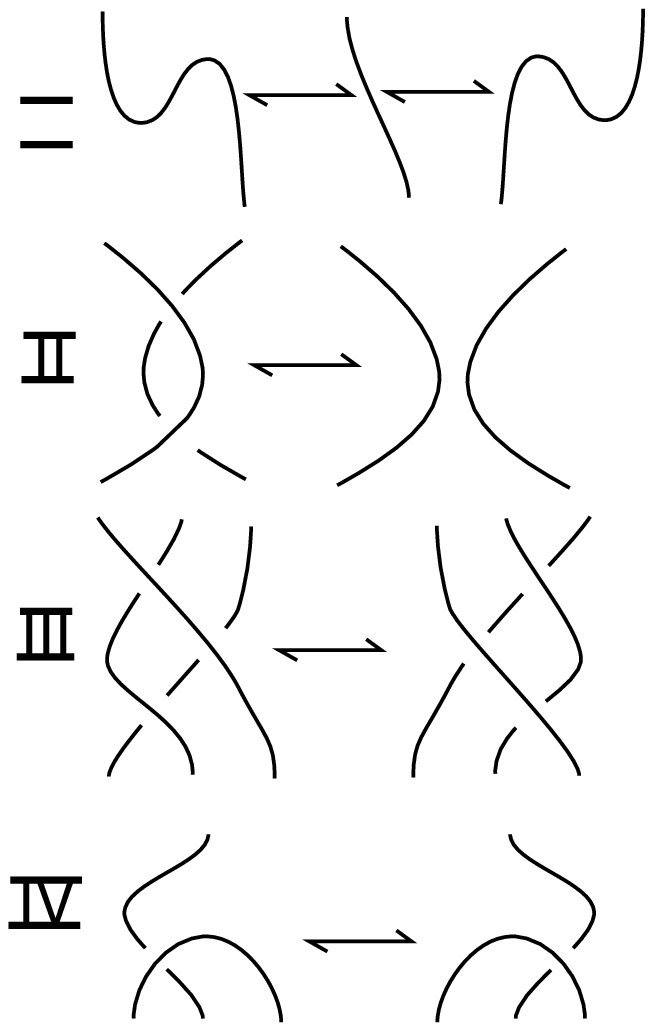}}
\vvvv

\begin{center}
{\bf Figure 36 - Regular Isotopy with respect to a Vertical Direction}
\end{center}
\newpage

All the moves taken together  are directly related to the axioms for a
quasi-triangular Hopf algebra  (aka quantum group). Many seeds of the structure
of Hopf algebras are prefigured in the patterns of link diagrams and the
structure of the category of tangles.    The interested reader can consult
\cite{RT}, \cite{Yetter}, \cite{GAUSS} and \cite{KNOTS},
\cite{KR}, \cite{K:alg} and section 5 of this paper for more information on this
point. \vspace{3mm}

Here is the list of the algebraic versions of the topological moves. Move $0$ is
the cancellation of maxima and minima. Move $II$  corresponds to the second
Reidemeister move. Move $III$ is the Yang-Baxter equation.  Move $IV$ expresses
the relationship of switching a line across a maximum. (There is a corresponding
version of $IV$ where the line is switched across a minimum.)

$$0.\hspace{.1in} M^{ai}M_{ib}  =\delta^{a}_{b}$$
$$II.\hspace{.1in} R^{ab}_{ij} \overline{R^{ij}_{cd}} = \delta^{a}_{c}
\delta^{b}_{d}$$
$$III.\hspace{.1in} R^{ab}_{ij}R^{jc}_{kf}R^{ik}_{de} =
R^{bc}_{ij}R^{ai}_{dk}R^{kj}_{ef}$$
$$IV. \hspace{.1in}R^{ai}_{bc}M_{id} = M_{bi}R^{ia}_{cd}$$

In the case of the Jones polynomial we have all the algebra present to make the
model. It is easiest  to indicate the model for the bracket polynomial: Let $cup$
and $cap$ be given by the $2$ x $2$ matrix $M$, described above so that $M_{ij}
=M^{ij}$. Let $R$ and $R^{-1}$ be given by the equations

$$R^{ab}_{cd} = AM^{ab}M_{cd} + A^{-1} \delta^{a}_{c} \delta^{b}_{d},$$
$$(R^{-1})^{ab}_{cd} = A^{-1}M^{ab}M_{cd} + A \delta^{a}_{c} \delta^{b}_{d}.$$

This definition of the R-matrices exactly parallels the diagrammatic expansion of
the bracket, and it is not hard to see, either by algebra or diagrams, that all
the conditions of the model are met. \vspace{3mm}

\subsection{Oriented Amplitudes} Slight but significant modifications are needed
to write the oriented version of the models we have discussed in the previous
section. See \cite{KNOTS}, \cite{Turaev}, \cite{LOMI}, \cite{Hennings}. In this
section we sketch the construction of oriented topological amplitudes.

The generalization  to oriented link diagrams naturally involves the introduction
of right and left oriented caps and cups. These are drawn as shown in Figure 37
below. \vspace{3mm}

\centerline{\includegraphics[scale=1.0]{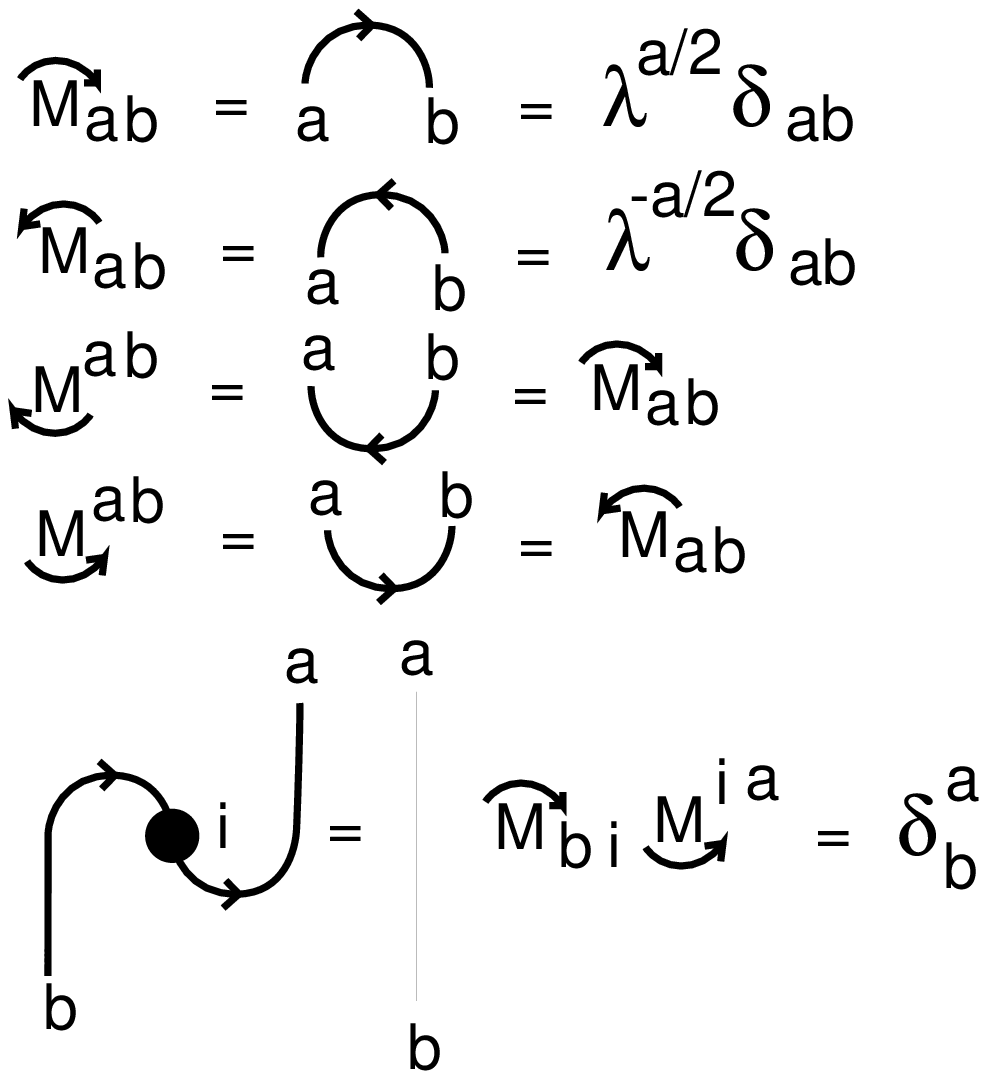}}
\vvvv

\begin{center}
{\bf Figure 37 - Right and Left Cups and Caps}
\end{center}

A right cup cancels  with a right cap to produce an upward pointing identity
line.  A left cup cancels with a left cap to produce a downward pointing identity
line. \vspace{3mm}

Just as we considered the simplifications that occur in the unoriented model by
taking the cup and cap matrices to be identical, lets assume here that right caps
are identical with left cups and that consequently left caps are identical with
right cups.  In fact, let us assume that the right cap and left cup are given by
the matrix $$M_{ab} = \lambda^{a/2} \delta_{ab}$$ where $\lambda$ is a constant
to be determined by the situation, and $\delta_{ab}$ denotes the Kronecker delta.
 Then the left cap and right cup are given by the inverse of $M$:
$$M_{ab}^{-1} =
\lambda^{-a/2} \delta_{ab}.$$

We assume that along with $M$ we are given a solution $R$ to the Yang-Baxter
equation, and that in an oriented diagram the specific choice of $R^{ab}_{cd}$ is
governed by the local orientation of the crossing in the diagram. Thus $a$ and
$b$ are the labels on the lines going into the crossing and $c$ and $d$ are the
labels on the lines emanating from the crossing.
\newpage

Note that with respect to the vertical direction for the amplitude, the crossings
can assume the aspects: both lines pointing upward, both lines pointing downward,
one line up and one line down (two cases). See Figure 38. \vspace{3mm}

\centerline{\includegraphics[scale=0.85]{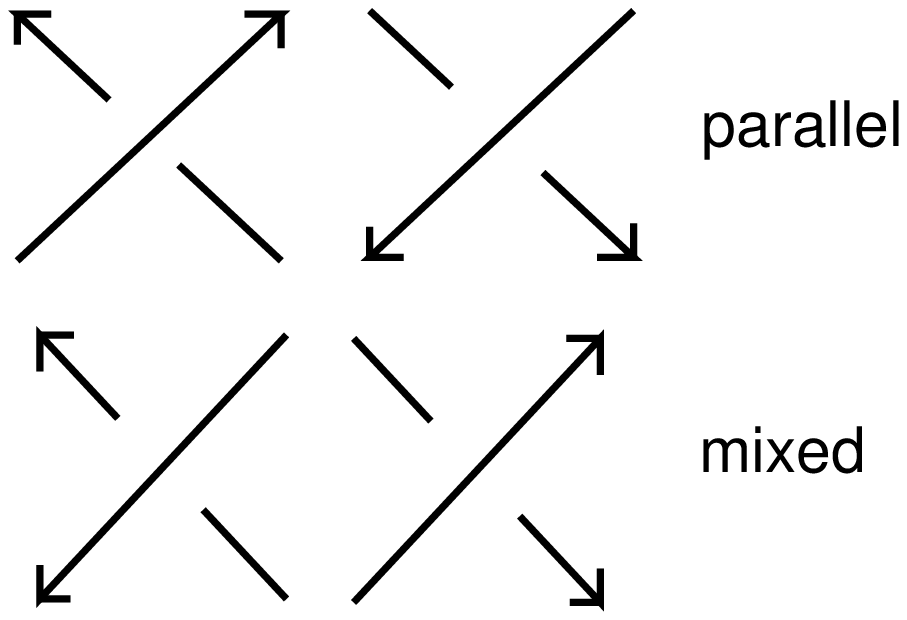}}
\vvvv

\begin{center}
{\bf Figure 38 - Oriented Crossings}
\end{center}
\vvv\vv

Call the cases of one line up and one line down the {\em mixed} cases and the
upward and downward cases the {\em parallel} cases. A given mixed crossing can be
converted ,in two ways, into a combination of a parallel crossing of the same
sign plus a cup and a cap. See Figure 39. \vspace{3mm}

\centerline{\includegraphics[scale=1.0]{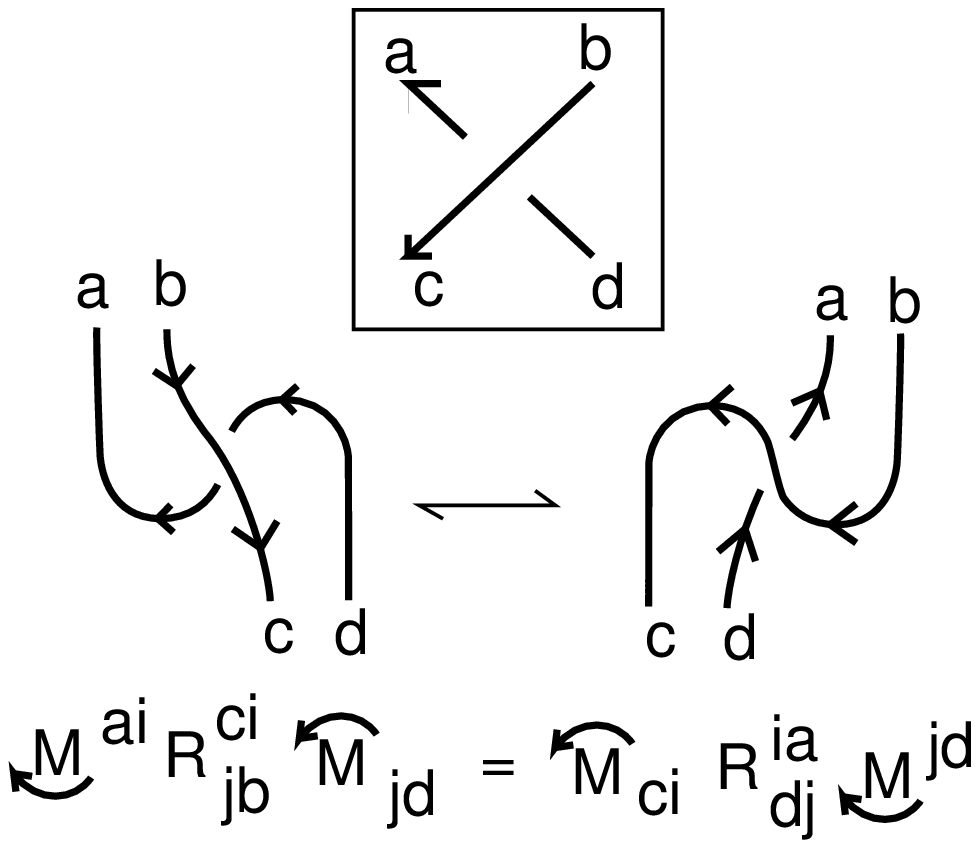}}
\vvvv

\begin{center}
{\bf Figure 39 - Conversion}
\end{center}
\newpage

This leads to an equation that must  be satisfied by the $R$ matrix in relation
to powers of $\lambda$ (again we use the Einstein summation convention):

$$\lambda^{a/2} \delta^{ai} R^{ci}_{jb} \lambda^{-d/2} \delta_{jd} =
\lambda^{-c/2} \delta_{ic} R^{ia}_{dj} \lambda^{b/2} \delta^{jb}.$$

This simplifies to the equation

$$\lambda^{a/2}   R^{ca}_{db} \lambda^{-d/2}  = \lambda^{-c/2} R^{ca}_{db}
\lambda^{b/2} ,$$
\vvvv

\centerline{\includegraphics[scale=1.0]{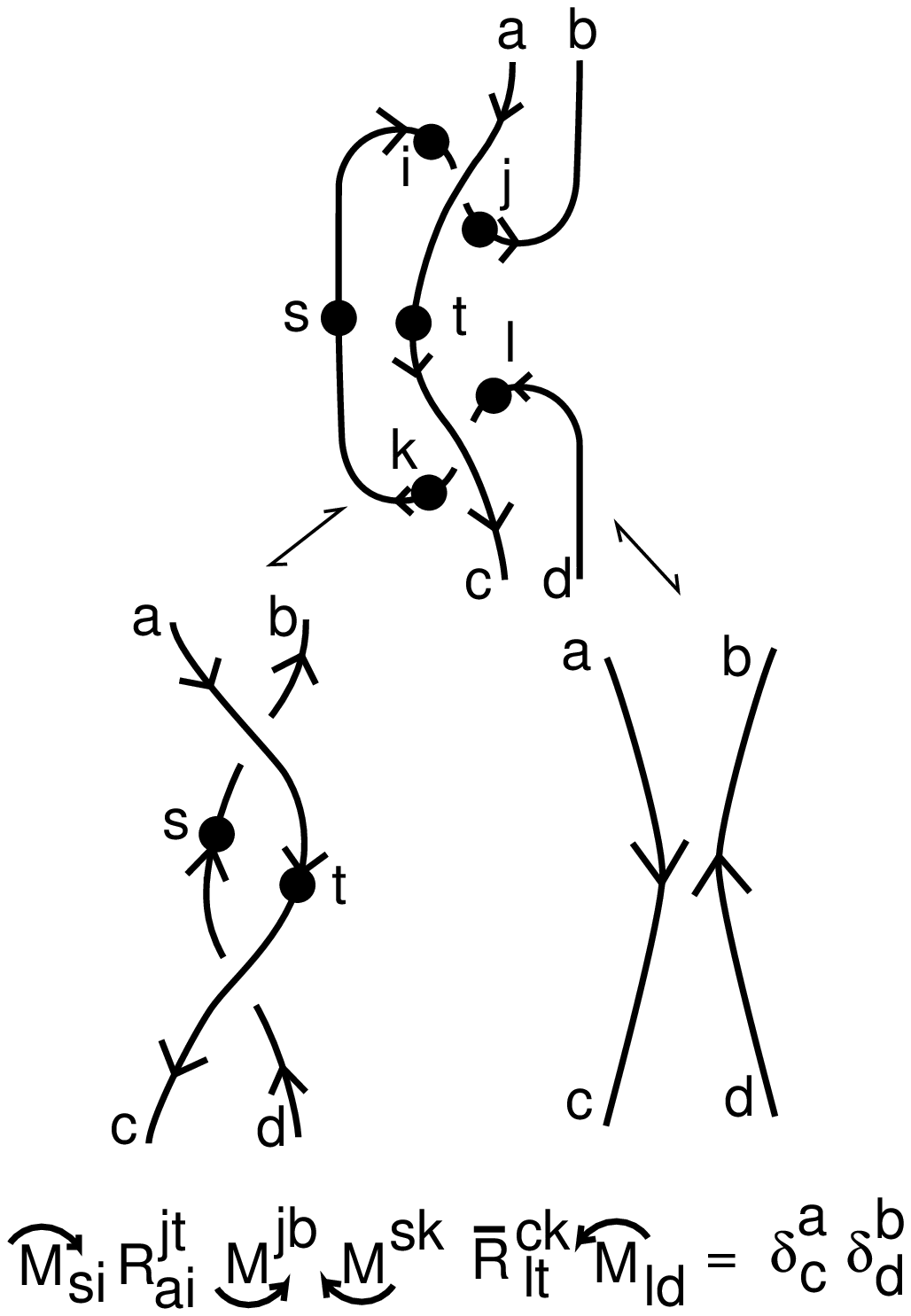}}
\vvvv

\begin{center}
{\bf Figure 40 - Antiparallel Second Move}
\end{center}
\newpage

\noindent from which we see that $R^{ca}_{db}$ is necessarily equal to zero unless $b+d =
a+c$.  We say that the $R$ matrix is {\em spin preserving} when it satisfies this
condition.  Assuming that the $R$ matrix is spin preserving,  the model will be
invariant under all orientations of the second and third Reidemeister moves just
so long as it is invariant under the anti-parallel version of the second
Reidemeister move as shown in Figure 40. \vspace{3mm}

This antiparallel version of the second Reidemeister move places the following
demand on the relation between $\lambda$ and $R$:

$$\sum_{st} \lambda^{(s-b)/2} \lambda^{(t-c)/2}R^{bt}_{as}
\overline{R^{cs}_{dt}}= \delta^{a}_{c} \delta^{b}_{d}.$$

Call this the {\em $R-\lambda$ equation.} The reader familiar with \cite{Jones}
or with the piecewise linear version as described in \cite{KNOTS} will
recognise this equation as the requirement for regular homotopy invariance in
these models. \vspace {3mm}

\subsection {Quantum Link Invariants and Vassiliev Invariants} Vassiliev
invariants can be used as building blocks for all the presently known quantum
link invariants. \vspace{3mm}

It is this  result that we can now make clear in the context of the models given
in our section on quantum link invariants.  Suppose that  $\lambda$ is written as
a power series in a variable $h$ , say $\lambda = exp(h)$ to be specific. 
Suppose also, that the R-matrices can be written as power series in $h$ with
matrix coefficients so that $PR = I + r_{+}h + O(h^2)$ and $PR^{-1} =  I + r_{-}h
+ O(h^2)$ where $P$ denotes the map of $V \otimes V$ that interchanges the tensor
factors. Let $Z(K)$ denote the value of the oriented amplitude described by this
choice of $\lambda$ and $R$. Then we can write $$Z(K) = Z_{0}(K) + Z_{1}(K)h +
Z_{2}(K)h^2 + ...$$ where each $Z_{n}(K)$ is an invariant of regular isotopy of
the link $K$. Furthermore, we see at once that {\em $h$ divides the series for
$Z(K_{+}) - Z(K_{-})$}.  By the definition of the Vassiliev invariants this
implies that $h^{k}$ divides $Z(G)$ if $G$ is a graph with $k$ nodes.  Therefore
$Z_{n}(G)$ vanishes if $n$ is less than the number of nodes of $G$. Therefore
$Z_{n}$ is a Vassiliev invariant of finite type $n$. Hence the quantum link
invariant is built from an infinite sequence of interlocked Vassiliev invariants.
\vspace{3mm}

It is an open problem whether the class of finite type Vassiliev invariants is
greater than those generated from quantum link invariants.  It is also possible
that there are quantum link invariants that cannot be generated by Vassiliev
invariants. \vspace{3mm}

\subsection {Vassiliev Invariants and Infinitesimal Braiding}

Kontsevich \cite{Kontsevich}, \cite{Bar-Natan} proved that a weight assignment
for a Vassiliev top row that satisfies the 4-term relation and the framing
condition (that the weights vanish for graphs with isolated double points)
actually extends to a Vassiliev invariant defined on all knots. His method is
motivated by the perturbative expansion of the Witten integral and by Witten's
interpretation of the integral in terms of conformal field theory. This section
will give a brief description of the Kontsevich approach and the questions that
it raises about the functional integral itself. \vspace{3mm}

The key to this approach is to see that the 4-term relations are a kind of
``infinitesimal braid relations". That is, we can re-write the 4-term relations in
the form of tangle operators as shown in Figure 41. \vspace{3mm}
\vvv

\centerline{\includegraphics[scale=0.85]{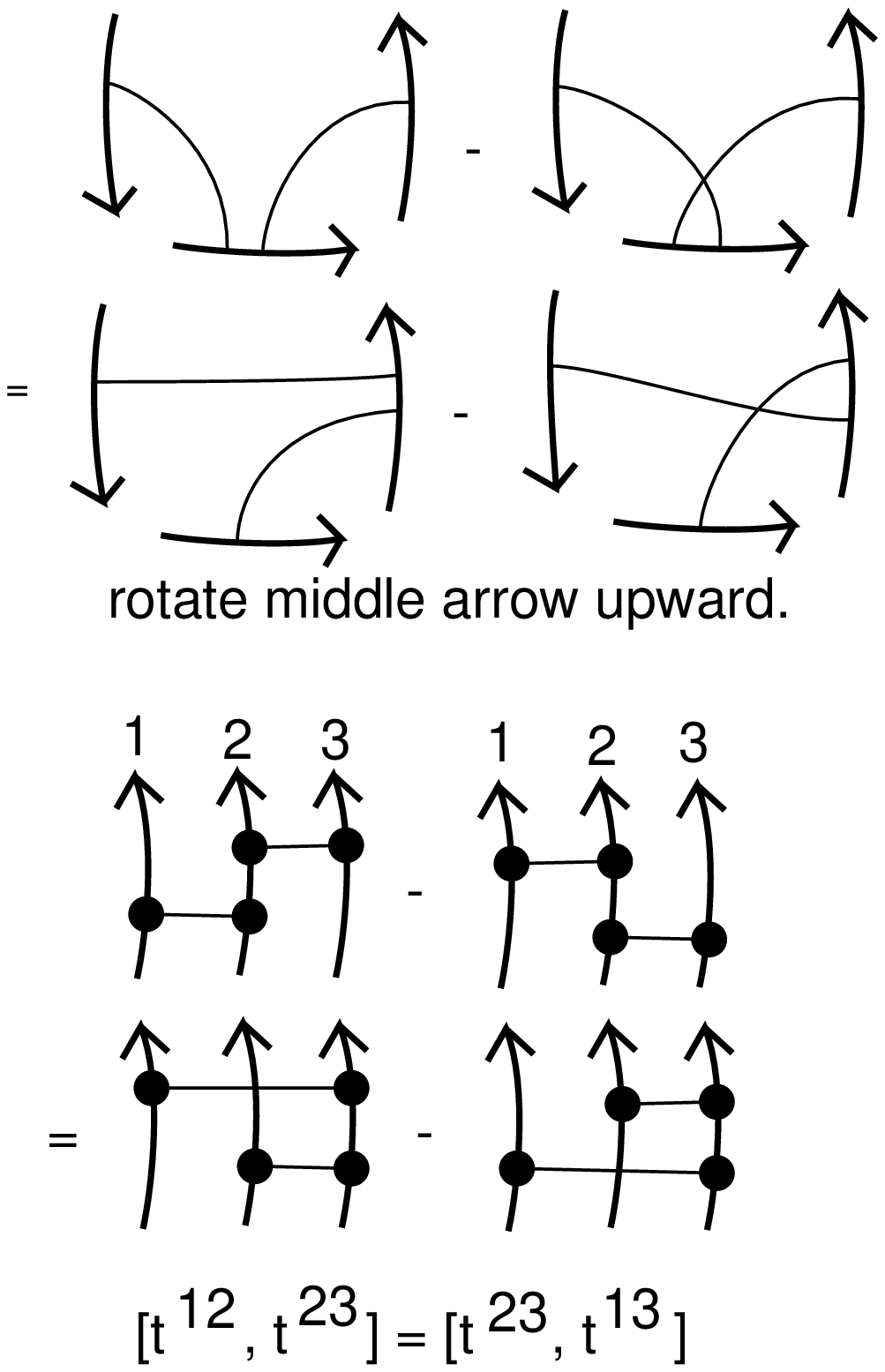}}
\vvvv

\begin{center}
{\bf Figure 41 - Infinitesimal Braiding}
\end{center}
\newpage

This shows that the commutator equation

$$[t^{12},t^{13}] + [t^{13}, t^{23}] = 0$$

\noindent is an algebraic form of the 4-term relation. The 4-term relations
translate exactly into these infinitesimal braid relations studied by Kohno
\cite{Kohno}. Kohno showed that his version of infinitesimal braid relations
corresponded to a flatness condition for a certain connection (the
Knizhnik-Zamolodchikov connection) and that this meant that these relations
constituted an integrability condition for making representations of the braid
group via monodromy.   Others have verified that the braid group representations
related to the Chern-Simons form and the Witten integral arise in this same way
from the Knizhnik-Zamolodchikov equations. In the case of Chern-Simons theory the
weights in the K-Z equations come from the Casimir of a classical Lie algebra,
just as we have discussed. Kontsevich observed that since the arbitrary 4-term
relations could also be regarded as an integrability condition it was possible to
use them in a generalization of Kohno's ideas to produce braid group
representations via iterated integration. He then generalized the process of
producing these braid group representations to the production of knot invariants
and these become realizations of Vassiliev invariants that have given admissible
weight systems for their top rows. \vspace{3mm}

The upshot of the Kontsevich work is a very specific integral formula for the
Vassiliev invariants. See \cite{Bar-Natan} for the specifics. It is clear from
the nature of the construction that the Kontsevich formula captures  the various
orders of perturbative terms in the Witten integral.  At this writing there is no
complete published description of this correspondence. \vspace{3mm}

\subsection{Weight Systems and the Classical Yang Baxter Equation}

Lets return momentarily to the series form of the solution to the Yang-Baxter
equation, as we had indicated it in the previous subsection.

$$PR = I + r_{+}h + O(h^2)$$ $$PR^{-1} =  I + r_{-}h + O(h^2).$$

Since $RR^{-1} = P$, it follows that $r_{-} = -r_{+}^{'}$ where $a^{'}$ denotes
the transpose of $a.$  In the case that $R$ satisfies the $R-\lambda$ equation,
it follows that

$$t = r_{+} - r_{-} = r + r^{'}$$

(letting $r$ denote $r_{+}$) satisfies the infinitesimal braiding relations

$$[t^{12},t^{13}] + [t^{13}, t^{23}] = 0.$$

It is interesting to contemplate this fact, since $r$, being the coefficient of
$h$ in the series for $R$,  necessarily satisfies the {\em classical Yang Baxter
Equation}

$$[r^{13},r^{23}] + [r^{12},r^{23}] +[r^{12},r^{13}]=0 .$$

\noindent (The classical Yang-Baxter equation for $r$ is a direct consequence of
the fact that $R$ is a solution of the (quantum) Yang-Baxter equation.)
\vspace{3mm}

Via the quantum link invariants, we have provided a special condition (the
assumption that $r$ is the coefficient of $h$ in a power series solution of the
quantum Yang-Baxter equation $R$, and that $R$ satisfies the $R-\lambda$
equation)  ensuring that a solution $r$ of the classical Yang Baxter equation
will produce a solution $t = r + r^{'}$ of the infinitesimal braiding relation,
whence a weight system for Vassiliev invariants.  More work needs to be done  to
fully understand the relationship between solutions of the classical Yang-Baxter
equation and the construction of Vassiliev invariants. \vspace{3mm}

\section{Hopf Algebras and Invariants of Three-Manifolds}

This section is a rapid sketch of the relationship between the description of
regular isotopy with respect to a vertical direction (as described in our
discussion of quantum link invariants) and the way that this formulation of the
Reidemeister moves is related to Hopf algebras and to the construction of link
invariants and invariants of three-manifolds via Hopf algebras. More detailed
presentations of this material can be found in \cite{Hennings}, \cite{GAUSS},
\cite{KR}, \cite{KRS}, \cite{K:alg}. \vspace{3mm}

Let's begin by recalling the Kirby calculus \cite{KCalc}. In the context of link
diagrams the Kirby calculus has an elegant formulation in terms of (blackboard)
framed links represented by link diagrams up to {\em ribbon equivalence}. Ribbon
equivalence consists in diagrams up to regular isotopy coupled with the
equivalence of a positive (negative) curl of Whitney degree $1$ with a positive
(negative) curl of Whitney degree $-1. $  See Figure 42. \vspace{3mm}

Here we refer informally to the Whitney degree of a plane curve. The Whitney
degree is the total turn of the tangent vector. If the curve is not closed, then
it is assumed that the tangent direction of the initial point is the same as the
tangent vector of the endpoint. In Figure 42 we illustrate how curls encode
framings and how ribbon equivalent curls correspond to identical framings.  A
link is said to be {\em framed} if it is endowed with a smooth choice of normal
vector field. Framing a link is equivalent to specifying an embedded band(s) of
which it is the core.  The {\em core} of a band is the center curve. Thus $S^{1}
\times \{.5 \}$ is the core of $S^{1} \times [0,1].$ \vspace{3mm}

\centerline{\includegraphics[scale=1.0]{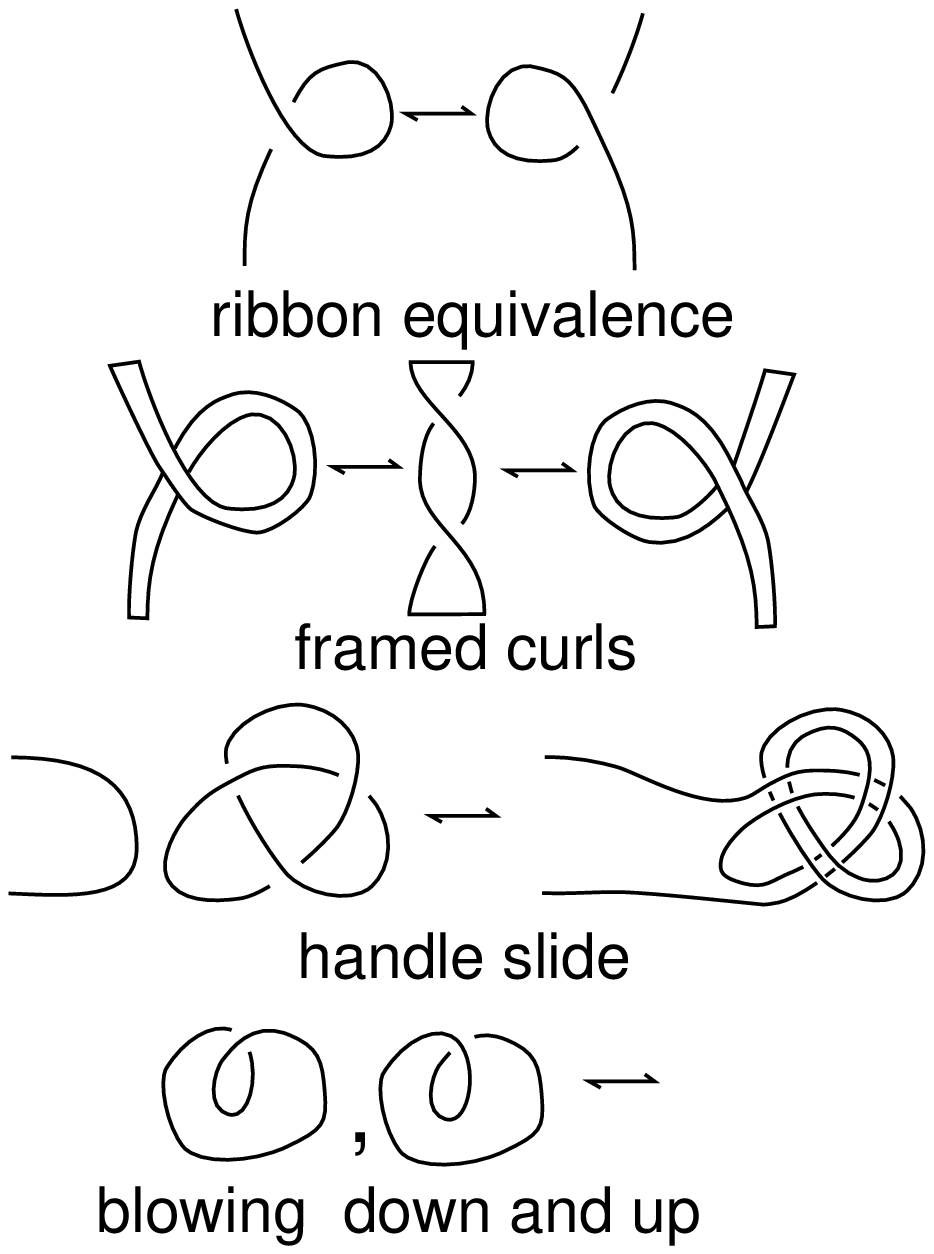}}
\vvvv

\begin{center}
{\bf Figure 42 - Framing and Kirby Calculus}
\end{center}

Now introduce two new moves on link diagrams called {\em handle sliding} and {\em
blowing up and down}.  These moves are illustrated in Figure 42.  Handle sliding
consists in duplicating a parallel copy of one link component and then band
connect summing it with another component. Blowing up consists in adding an
isolated unknotted component with a single curl. Blowing down consists in
deleting such a component.  These are the basic moves of the Kirby Calculus. Two
link diagrams  are said to be {\em KC-equivalent} if there is a combination of
ribbon equivalence, handle-sliding and blowing up and blowing down that takes one
diagram to the other. \vspace{3mm}

The invariants of three-manifolds described herein are based on the
representation of closed three-manifolds via surgery on framed links. Let
$M^{3}(K)$ denote the three-manifold obtained by surgery on the blackboard framed
link corresponding to the diagram $K$. In $M^{3}(K)$ the longitude associated
with the diagram, as shown in Figure 43, bounds the meridian disk of the solid
torus attached via the surgery. \vspace{3mm}

\centerline{\includegraphics[scale=1.0]{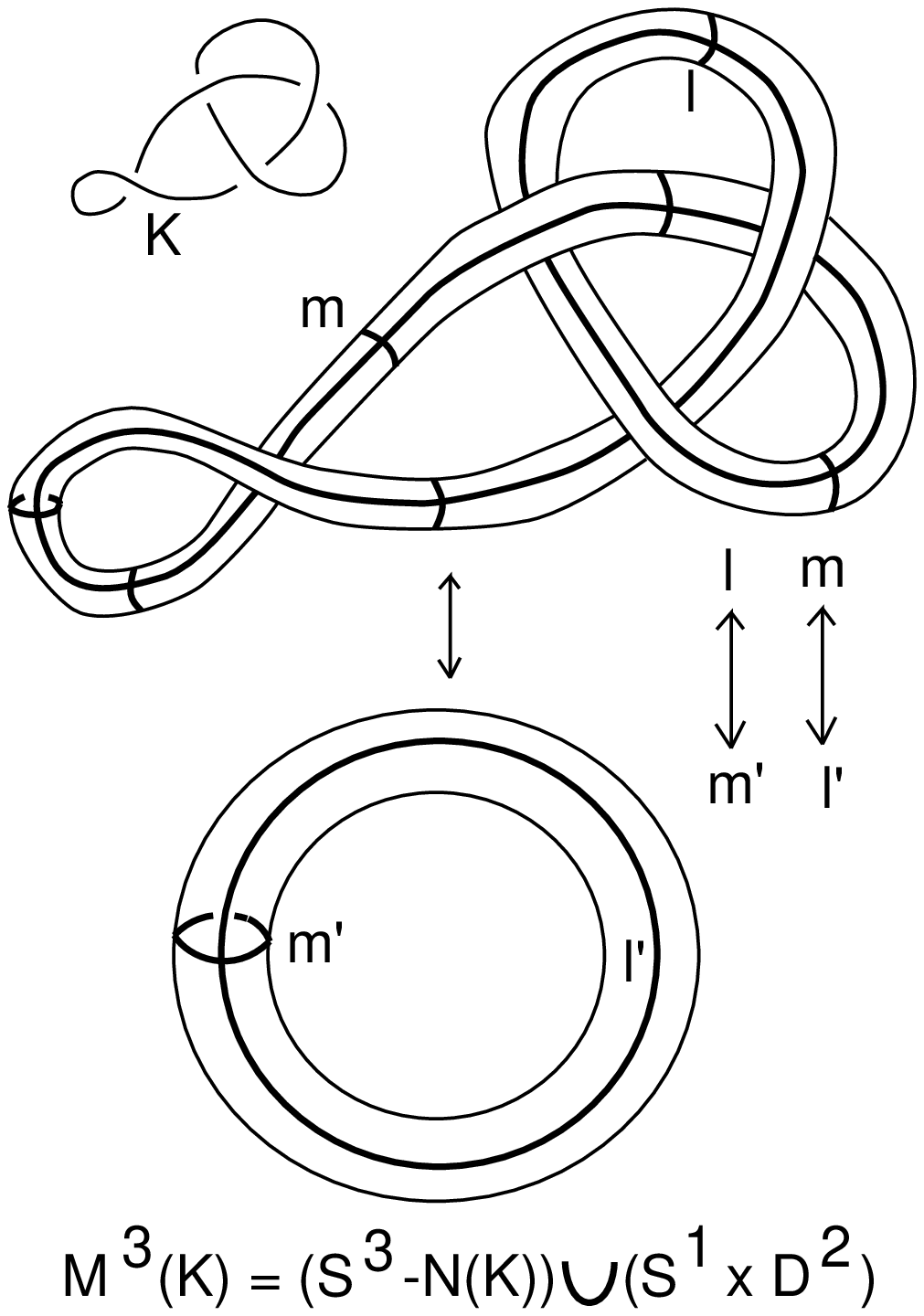}}
\vvvv

\begin{center}
{\bf Figure 43 - Surgery on a Blackboard Framed Link}
\end{center}

The basic result about Kirby Calculus is that  $M^{3}(K)$ is homeomorphic to
$M^{3}(L)$ if and only if $K$ and $L$ are KC equivalent. Thus invariants of links
that are also invariant under Kirby moves will produce invariants of
three-manifolds. It is the purpose of this section to sketch one on the
approaches to constructing such invariants.
\vv

The ideas behind this approach are quite simple. We are given a finite
dimensional quasitriangular Hopf algebra $A$.  We associate to A a tensor
category $Cat(A).$ The objects in this category are the base field $k$ of the
Hopf algebra, and tensor powers of a formal object $V$.  It is assumed that the
tensor powers of $V$ are canonically associative and that the tensor product of V
with k on either side is canonically isomorphic to V.  The  morphisms in $Cat(A)$
are represented by Hopf algebra decorated immersed curves  arranged with respect
to a vertical direction.
\vv

An {\em immersed curve diagram} is a link diagram where there is no distinction
between undercrossings and overcrossings. Segments of the diagram can cross one
another transversely as in a standard link diagram, and we can arrange such a
diagram with respect to a vertical direction just as we did for link diagrams.  A
{\em vertical place} on such a diagram is a point that is not critical with
respect to the vertical direction, and is not a crossing. A {\em decoration} of
such an immersed curve diagram consists in a subset of vertical places labelled
by elements of the Hopf algebra $A.$  The diagrams can have endpoints and these
are either at the bottom of the diagram or at the top (with respect to the
vertical).  The simplest decorated diagram is a vertical line segment with a
label $a$ (corresponding to a element $a$ of the Hopf algebra) in its interior.
In the category $Cat(A)$ this segment  is regarded as a morphism $[a]: V
\longrightarrow V$ where $V$ is the formal object alluded to above. Composition
of these morphisms corresponds to multiplication in the algebra: $[a][b] = [ab].$
By convention, we take the order of multiplication from bottom to top with
respect to the vertical direction.
\vv

A tensor product $a \otimes b$  in  $A \otimes A$  is represented by two parallel
segments, one decorated by $a$, the other decorated by $b$. It is our custom to
place the decorations for $a$ and for $b$ at the same level in the diagram.  In
the Hopf algebra we have the coproduct $\Delta: A \longrightarrow A\otimes A.$ We
shall write $$\Delta(a) =\Sigma a_{1} \otimes a_{2}$$ where it is understood 
that  this means that the coproduct of $a$ is a sum over elements of the {\em
form}  $a_{1} \otimes a_{2}.$   It is also useful to use a version of the
Einstein summation convention and just write $$\Delta(a) = a_{1} \otimes a_{2}$$ 
where it is understood that the right hand side is a summation.   In diagrams,
application of the antipode makes parallel lines with doubled decorations
according to the two factors of the coproduct.  See Figure 44.
\vv

A crossing of two undecorated segments is regarded as a morphism $P: V \otimes V
\longrightarrow  V \otimes V.$ Since the lines interchange, we expect $P$ to
behave as the permutation of the two tensor factors.  That is, we take the
following formula to be axiomatic:
$$P \circ ([a] \otimes [b]) = ([b] \otimes [a]) \circ P.$$

A cap  (see Figure 44) is regarded as a morphism from $V \otimes V$ to $k$, while
a cup is regarded as a morphism form $k$ to $V \otimes V.$   As in the case of
the crossing the relevance of these morphisms to the category is entirely encoded
in their properties. The basic property of the cup and the cap is that {\em if
you ``slide" a decoration across the maximum or minimum in a counterclockwise
turn, then the antipode $S$ of the Hopf algebra is applied to the decoration.} 
In categorical terms this property says
\newpage

$$Cup \circ ([a] \otimes 1) = Cup \circ (1 \otimes [Sa])$$
and
$$([Sa] \otimes 1) \circ Cap  =  (1 \otimes [a]) \circ Cap.$$
These properties and some other naturality properties of the cups and the
caps are illustrated in Figure 44.  These naturality properties of the flat
diagrams include regular homotopy of immersions, as illustrated in Figure 44.
\vvvv\vvvv

\centerline{\includegraphics[scale=1.0]{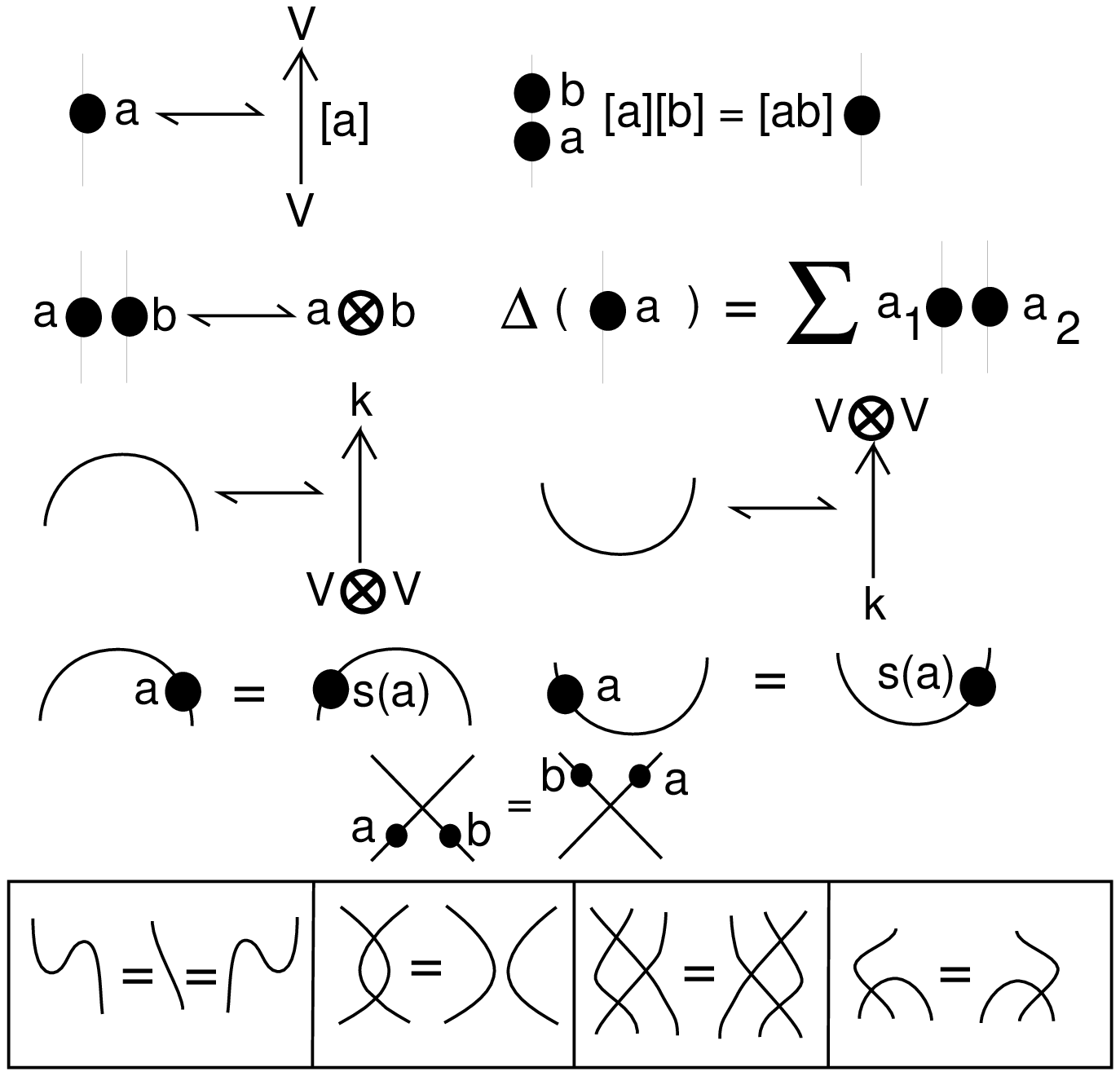}}
\vvvv

\begin{center}
{\bf Figure 44 - Morphisms in $Cat(A)$}
\end{center}
\newpage

\centerline{\includegraphics[scale=1.0]{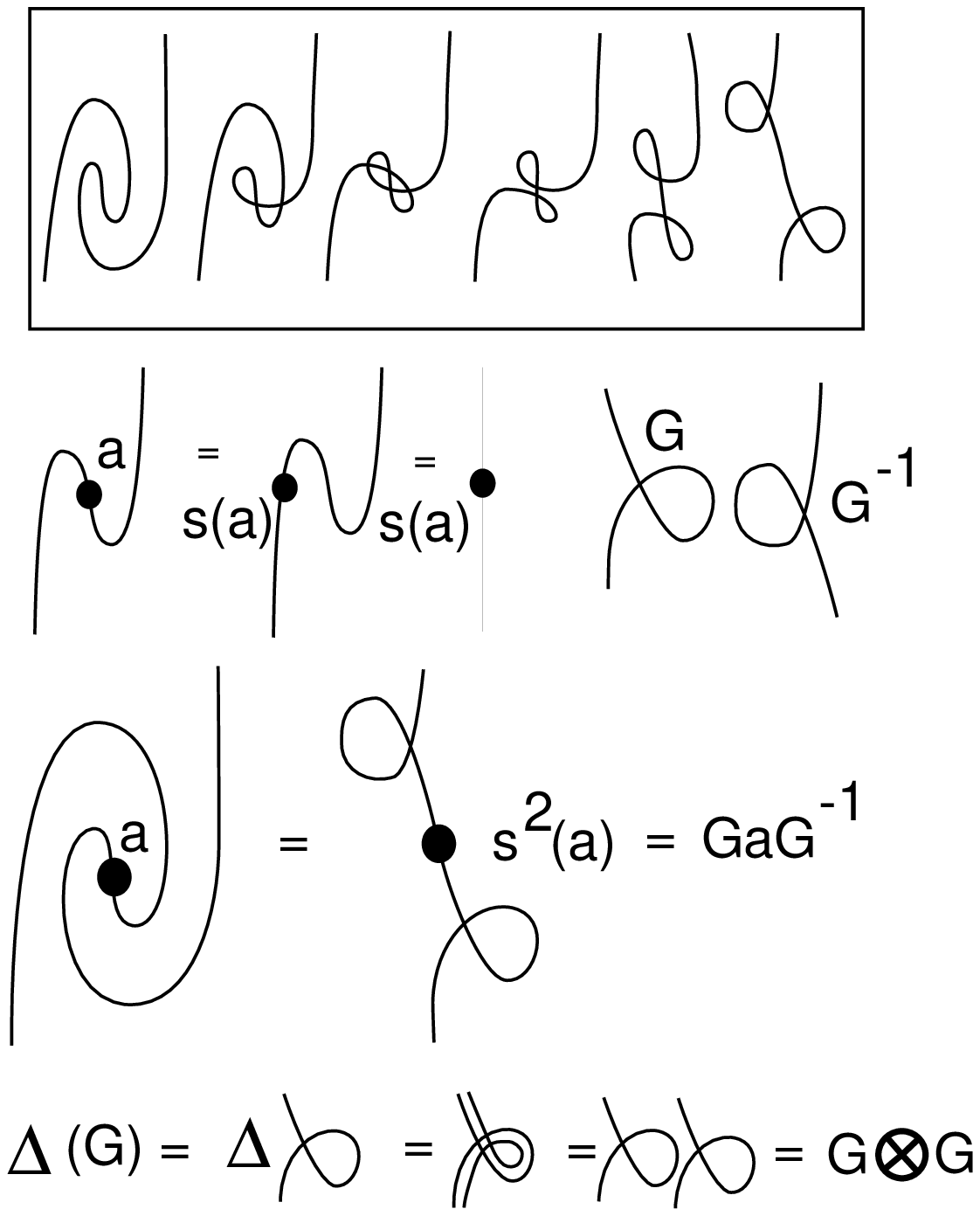}}
\vvvv

\begin{center}
{\bf Figure 45 - Diagrammatics of the Antipode}
\end{center}
\vvv

In Figure 45 we see how this property of the cups and the caps leads to a
diagrammatic interpretation of the antipode. This, in turn, leads to the
interpretation of the flat curl as a grouplike element $G$ in $A$ such that
$S^{2}(a) = GaG^{-1}$ for all $a$ in $A.$  $G$ is a flat curl diagram interpreted
as a morphism in the category. We see that formally it is natural to interpret
$G$ as an element of $A$ and that $\Delta (G) = G \otimes G$ is a direct
consequence of the diagrams for $Cat(A).$ In a so-called {\em ribbon Hopf
algebra} there is such a grouplike already in the algebra. In the general case it
is natural to extend the algebra to contain such an element. \vspace{3mm}

We are now in a position to describe a functor $F$ from the  tangle category $T$
to $Cat(A).$  (The tangle category is defined for link diagrams without
decorations. It has the same objects as $Cat(A).$ The morphisms in the tangle
category have relations corresponding to the augmented Reidemeister moves
described in the section on quantum link invariants.) $F$ simply decorates each
positive (with respect to the vertical - see Figure 45)  crossing of the tangle
with the Yang-Baxter element (given by the quasi-triangular Hopf algebra $A$)
$\rho = \Sigma e \otimes e^{'}$ and each negative crossing (with respect to the
vertical) with $\rho^{-1} = \Sigma S(e) \otimes e^{'}$. The form of the
decoration is indicated in Figure 46.
\vvvv

\centerline{\includegraphics[scale=0.9]{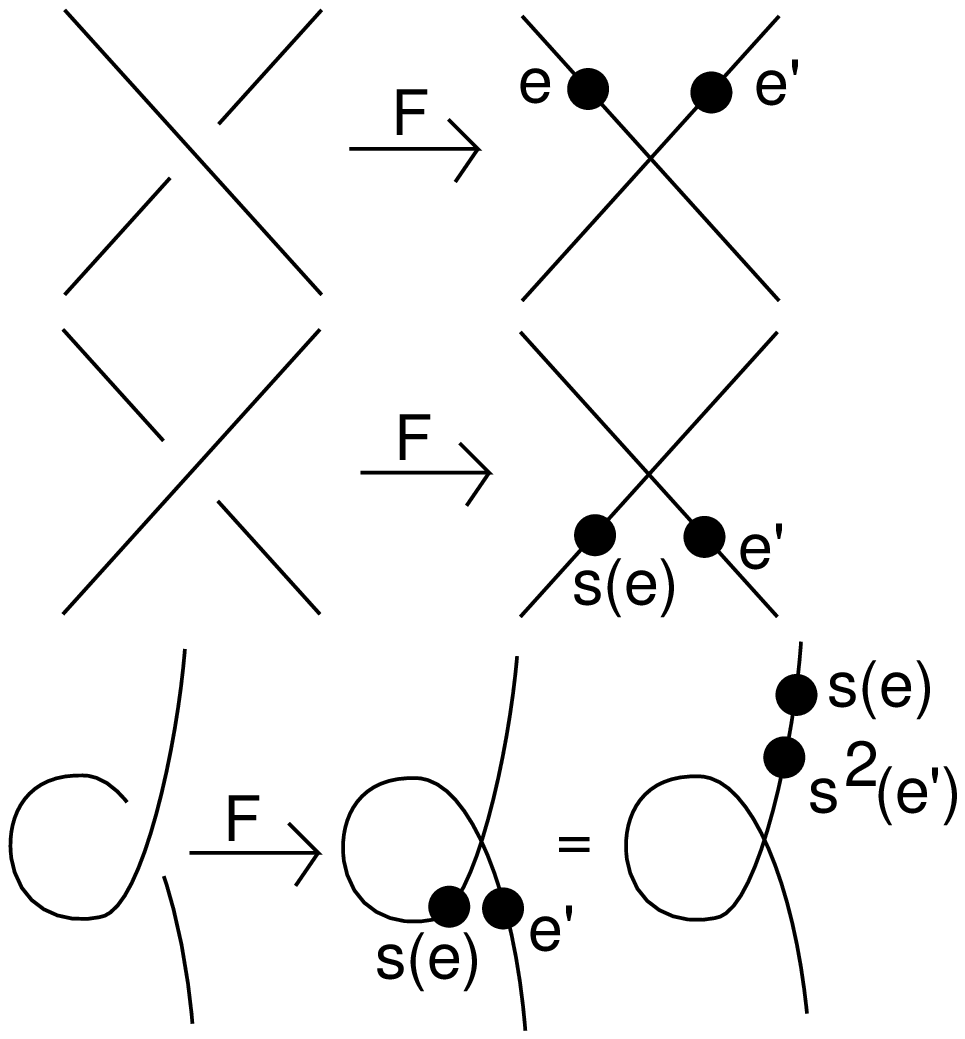}}
\vvvv

\begin{center}
{\bf Figure 46 - The Functor $F: T \longrightarrow Cat(A).$}
\end{center}
\vvv

The key point about this category is that because Hopf algebra elements can be
moved around the diagram, we can concentrate all the algebra in one place.
Because the flat curls are identified with either $G$ or $G^{-1}$, we can use
regular homotopy of  immersions to bring each component of a link diagram to the
form of a circle with a single concentrated decoration (involving a sum over many
products). An example is shown in Figure 46.  Let us denote by $\lambda(a): k
\longrightarrow k$ the morphism that corresponds to decorating the right hand
side of a standard circle with $a$. That is, $\lambda(a) = Cap \circ (1 \otimes
[a]) \circ Cup.$ We can regard $\lambda$ as a linear functional defined on $A$ as
a vector space over $k.$ \vspace{3mm}

We wish to find out what properties of $\lambda$ will be appropriate for
constructing invariants of three-manifolds. View Figure 47.  Handle sliding is
accomplished by doubling a component and then band summing. The doubling
corresponds to applying the antipode. As a result, we have that in order for
$\lambda$ to be invariant under handle-sliding it is sufficient that it have the
property  $\lambda(x)1 = \Sigma \lambda ( x_{1}) x_{2}.$  This is the formal
defining property of a {\em right integral} on the Hopf algebra $A.$  Finite
dimensional Hopf algebras have such functionals and suitable normalizations lead
to well-defined three-manifold invariants.  For more information see the
references cited at the beginning of this section. This completes our capsule
summary of Hopf algebras and invariants of three-manifolds. \vspace{3mm}
\vvv

\centerline{\includegraphics[scale=1.0]{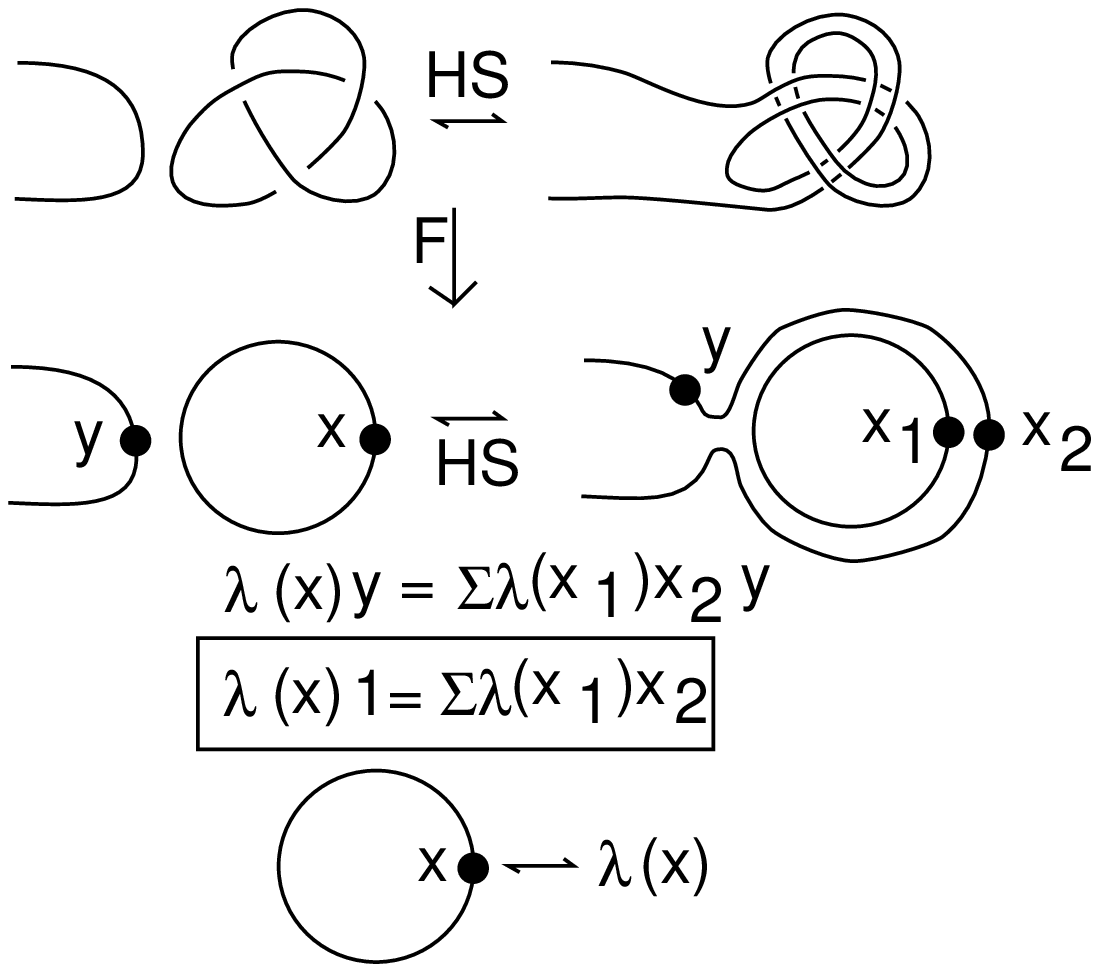}}
\vvvv

\begin{center}
{\bf Figure 47 - Handle Sliding and Right Integral}
\end{center}
\vvvv

There are a number of problems related to this formulation of invariants of
three-manifolds. First of all, while it is the case that the invariants that come
from integrals can be  different from invariants defined through representations
of Hopf algebras as in \cite{RT} it is quite difficult to compute them and
consequently little is known. Another beautiful problem is related to the work of
Greg Kuperberg \cite{Kup1}, \cite{Kup2}. Kuperberg defines invariants of
three-manifolds associated via Hopf diagrams associated with a Heegard splitting
of the three-manifold. Does our invariant on the Drinfeld double of a Hopf
algebra $H$ give the same result as Kuperberg's invariant for $H$?  This
conjecture is verified in the (easy) case where $H$ is the group ring of a finite
group. Finally, it should be mentioned that the way in which handle-sliding
invariance is proven for the universal three-manifold invariant of finite type of
Le and Murakami \cite{LeMur} is directly analogous to our method of relating
handle sliding, coproduct and right integral. It remains to be seen what is the
relationship between three-manifold invariants of finite type and the
formulations discussed here.
\newpage

\section{Temperley-Lieb Algebra}

This section is devoted to the structure of the Temperley-Lieb algebra as
revealed by its diagrammatic interpretation. We begin with a combinatorial
description of this algebra. It is customary, in referring to the Temperley-Lieb
algebra to refer to a certain free algebra over an appropriate ring. This free
algebra is the analog of the group ring of the symmetric group $S_{n}$ on $n$
letters. It is natural therefore to first describe that multiplicative structure
that is analogous to $S_{n}.$ We shall call this structure the {\em Temperley-Lieb
 Monoid $M_{n}$}. We shall describe the Temperley-Lieb
algebra itself after first defining this monoid. \vspace{3mm}

There is one Temperley-Lieb monoid, $M_{n}$, for each natural number $n$. The
{\em connection elements} of  $M_{n}$ consist in diagrams in the plane that make
connections involving two rows of $n$ points. These rows will be referred to as
the {\em top} and {\em bottom} rows.  Each point in a row is paired with a unique
point different from itself in either the top or the bottom row (it can be paired
with a point in its own row). These pairings are made by arcs drawn in the space
between the two rows. {\em No two arcs  are allowed to intersect one another}. 
Such a connection element will be denoted by $U$, with subscripts to indicate
specific elements. If the top row is the set $Top=\{1,2,3,...,n\}$ and the bottom
row is $Bot=\{1',2',...,n'\},$ then we can regard $U$ as a function from $Top
\bigcup Bot$ to itself such that $U(x)$ is never equal to $x$, $U(U(x)) = x$ for
all $x$, and satisfying the planar non-intersection property described above. In
topological terms $U$ is an $n$-tangle with no crossings, taken up to regular
isotopy of tangles in the plane. \vspace{3mm}

If $U$ and $V$ are two elements of $M_{n}$ as described above, then their product
$UV$ is the tangle product obtained by attaching the bottom row of $U$ to the top
row of $V.$ Note that the result of taking such a product will produce a new
connection structure plus some loops in the plane. Each loop is regarded as an
instance of the {\em loop element $\delta$} of the Temperley-Lieb monoid $M_{n}.$
The loop element commutes with all other elements of the monoid and has no other
relations with these elements. Thus $UV = \delta^{k} W$ for some non-negative
integer $k$,  and  some connection element $W$ of the monoid. \vspace{3mm}

The Temperley-Lieb algebra $T_{n}$ is  the free additive module on $M_{n}$ modulo
the identification $$\delta = -A^{2} - A^{-2},$$ over the ring $Z[A,A^{-1}]$ of
Laurent polynomials in the variable $A.$  Products are defined on the connection
elements and extended linearly to the algebra.  The reason for this loop
identification is the application of the Temperley-Lieb algebra for the bracket
polynomial and for representations of the braid group \cite{Jones}, \cite{state},
\cite{TL}. \vspace{3mm}

The Temperley-Lieb monoid  $M_{n}$ is generated by the elements $1, U_{1}, U_{2},
...,U_{n-1}$ where the identity element $1$ connects each $i$ in the top row with
its corresponding member $i'$ in the bottom row. Here $U_{k}$ connects $i$ to
$i'$ for $i$ not equal to $k,k+1, k',(k+1)'.$  $U_{k}$ connects $k$ to $k+1$ and
$k'$ to $(k+1)'$.   It is easy to see that

$$U_{k}^{2} = \delta U_{k},$$ $$U_{k}U_{k \pm 1}U_{k} = U_{k},$$ $$U_{i}U_{j} =
U_{j}U_{i} , \hspace{.1in}  | i-j | > 1.$$

See Figure 48. \vspace{3mm}

\centerline{\includegraphics[scale=1.0]{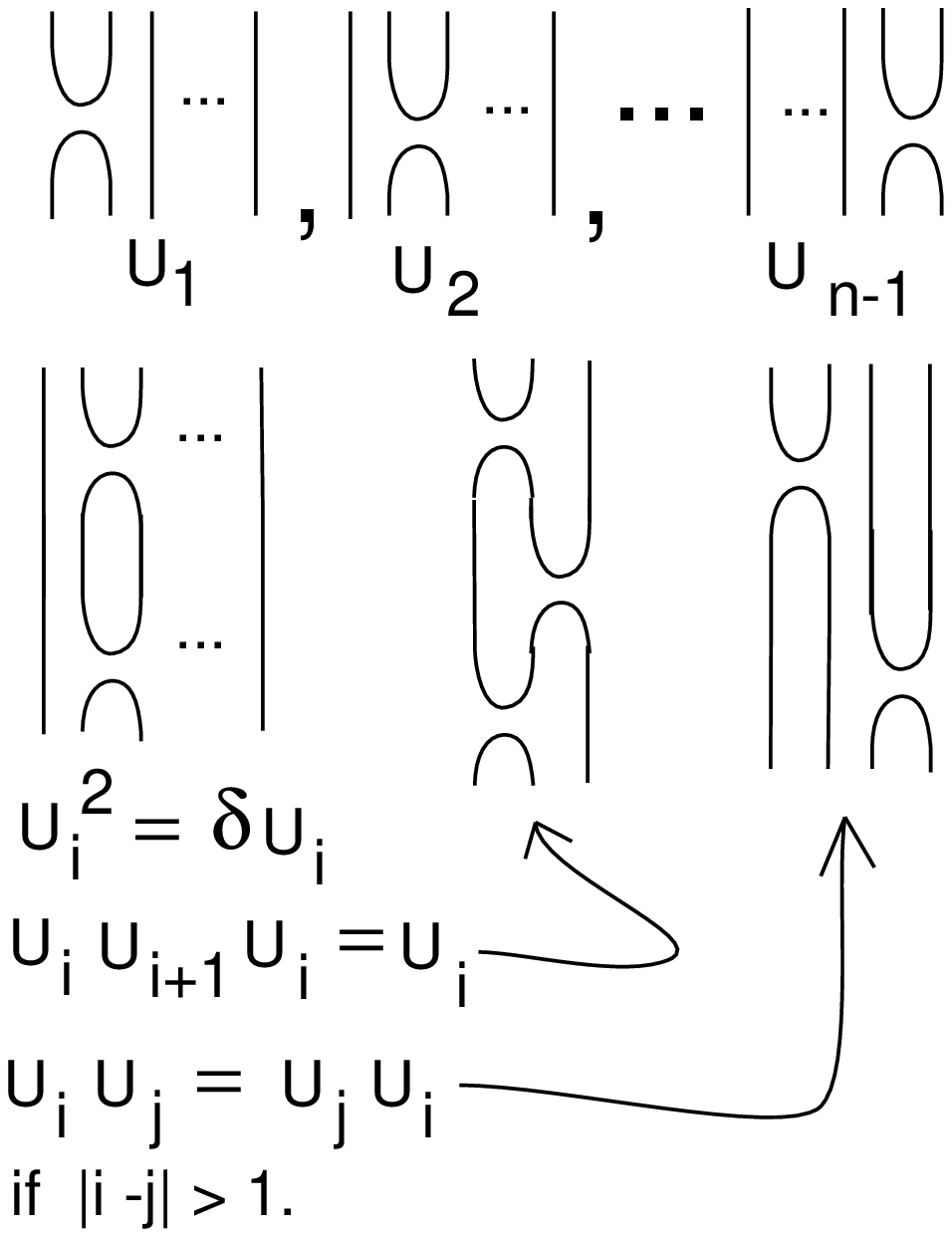}}
\vvvv

\begin{center}
{\bf Figure 48 - Relations in the Temperley-Lieb Monoid}
\end{center}

We shall prove that the Temperley-Lieb Monoid is the universal monoid on $ G_{n}
= \{ 1, U_{1}, U_{2}, ...,U_{n-1} \}$ modulo these relations. In order to
accomplish this end we give a direct diagrammatic method for writing any
connection element of the monoid as a certain canonical product of elements of
$G_{n}.$  This method is illustrated in Figure 49. \vspace{3mm}

\centerline{\includegraphics[scale=1.0]{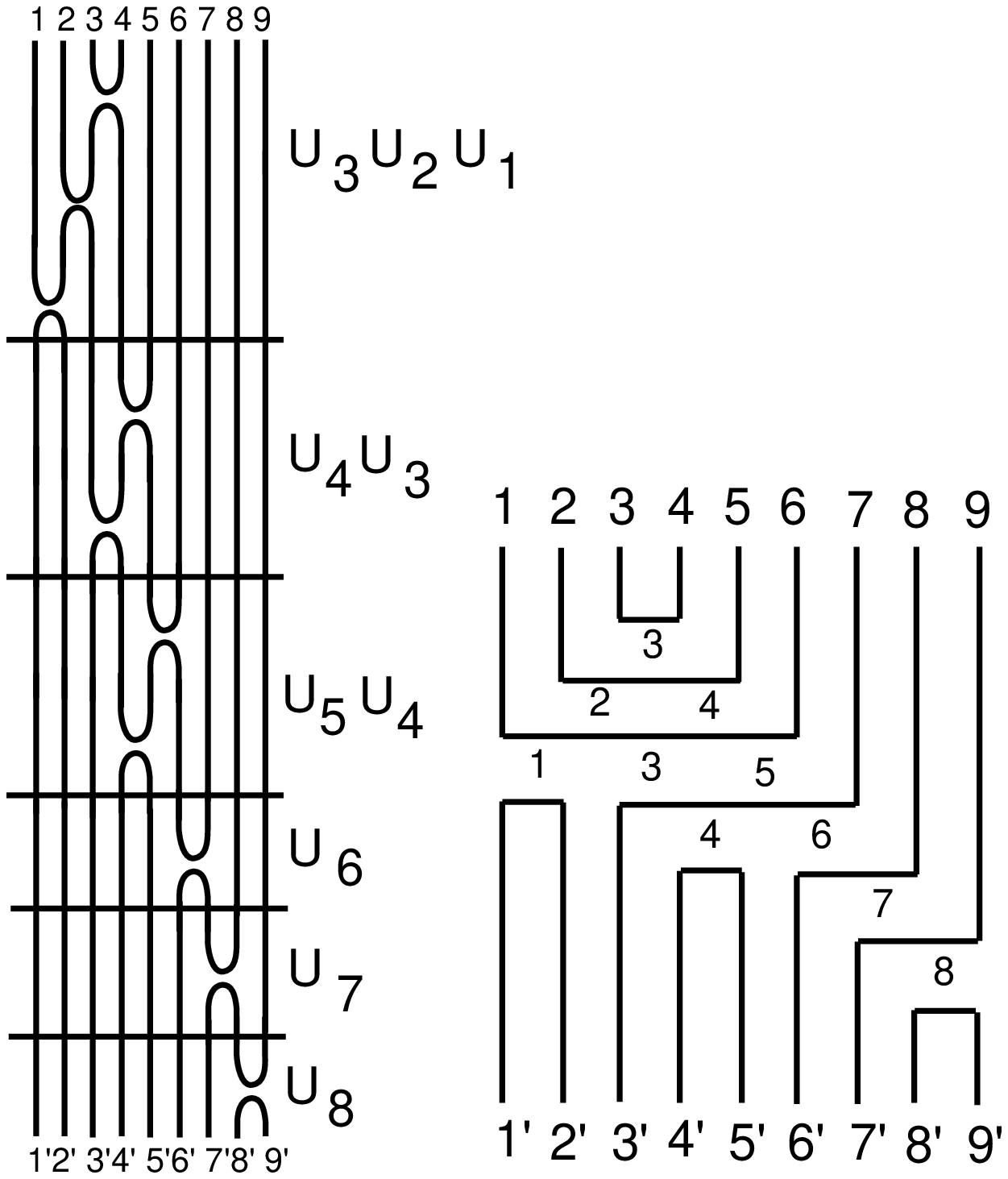}}
\vvvv

\begin{center}
{\bf Figure 1 - Figure 49 - Canonical Factorization in the Temperley-Lieb Monoid }
\end{center}
\vvv

As shown in Figure 49, we represent the connection diagram with vertical and
horizontal straight arcs such that except for the height of the straight arcs,
the form of the connection between any two points is unique - consisting in two
vertical arcs and one horizontal arc. The horizontal arc has its endpoints on
the vertical lines that go through the row points that are being connected.
(Diagram is drawn so that each pair in the set $\{(i,i'): 1 \leq n \}$ determines
such a vertical line. The vertical arcs in the connection are chosen as segments
from these vertical lines. All connections are chosen so that the connections do
not intersect. It is from this diagram that we shall read out a factorization
into a product of elements of $G_{n}$. \vspace{3mm}

The factorization is achieved via a decoration of the straight arc diagram by
dotted vertical arcs, as shown in Figure 49.  Each dotted arc connects midpoints
of  the restrictions of horizontal arcs to the {\em columns} of the diagram, where a
column of the diagram is the space between two consecutive vertical lines
(vertical lines described as in the previous paragraph).  The {\em index} of a
column is the row number associated to the left vertical boundary of the column.
The dotted lines in a given column are uniquely determined by starting at the
bottom of the column and pairing up the horizontal arcs in that column in
vertical succession. Each dotted arc is labelled by the index of the column in
which it stands. \vspace{3mm}

In a given diagram a  {\em sequence} of dotted arcs  is a maximal set of dotted
arcs (with consecutive indices) that are interconnected by horizontal segments
such that one can begin at the top of the dotted arc (in that sequence) of
highest index, go down the arc and left to the top of the next arc along a
horizontal segment, continuing in this manner until the whole sequence is
traversed.  It is clear from the construction of the diagram that the dotted
segments in the diagram collect into a disjoint union of sequences $\{ s^{1},
s^{2},...,s^{k} \}$  where each $s^{i}$ denotes the corresponding descending
sequence of of consecutive indices:

$$s^{i}=(m_{i},m_{i}-1,m_{i}-2,..., n_{i}+1, n_{i}).$$

\noindent These indices satisfy the inequalities:

$$m_{1} < m_{2} < m_{3} < ... < m_{k}$$

\noindent and

$$n_{1} < n_{2} < n_{3} < ... < n_{k}.$$

The sequences  $\{ s^{1}, s^{2},...,s^{k} \}$  occur in that order on the diagram
read from left to right. Of course the descent of each sequence goes from right
to left. If $D$ is a diagram with sequence structure $s(D)$ as we have just
described, let $U(s(D))$ be the following product of generators of the Temperley-Lieb monoid:

$$U(D) = U(s^{1})U(s^{2})U(s^{3})...U(s^{k})$$ where

$$U(s^{i}) = U_{m_{i}}U_{m_{i}-1} ...U_{n_{i}+1}U_{n_{i}}.$$

By looking carefully at the combinatorics of these diagrams, as illustrated in
Figure 49, one sees that $D$ and $U(D)$ represent the same connection structure
in the Temperley-Lieb monoid. Furthermore, any sequence structure $s$ satisfying
the inequalities given above (call these the canonical inequalities)  will
produce a standard diagram from the product $U(s).$  Thus the sequence structure
of a Temperley-Lieb diagram completely classifies this diagram as a connection
structure in the monoid. \vspace{3mm}

We must now prove that any product of elements of  $ G_{n} = \{ 1, U_{1}, U_{2},
...,U_{n-1} \}$  can be written, up to a loop factor, as $U(s(D))$ for some
diagram $D$, or equivalently as $U(s)$ for a sequence structure satisfying the
canonical inequalities. This is a simple exercise in using the relations we have
already given for the products of elements of $G_{n}.$  We leave the details to
the reader. This completes the proof that the relations in Figure 48 are a
complete set of relations for the Temperley-Lieb monoid. \vspace{3mm}

\subsection {Parentheses}

Elements of the Temperley-Lieb monoid $M_{n}$ are in one to one correspondence
with well-formed parenthesis expressions using $n$ pairs of parentheses. The
proof of this statement is shown in Figure 50. \vspace{3mm}

\centerline{\includegraphics[scale=1.0]{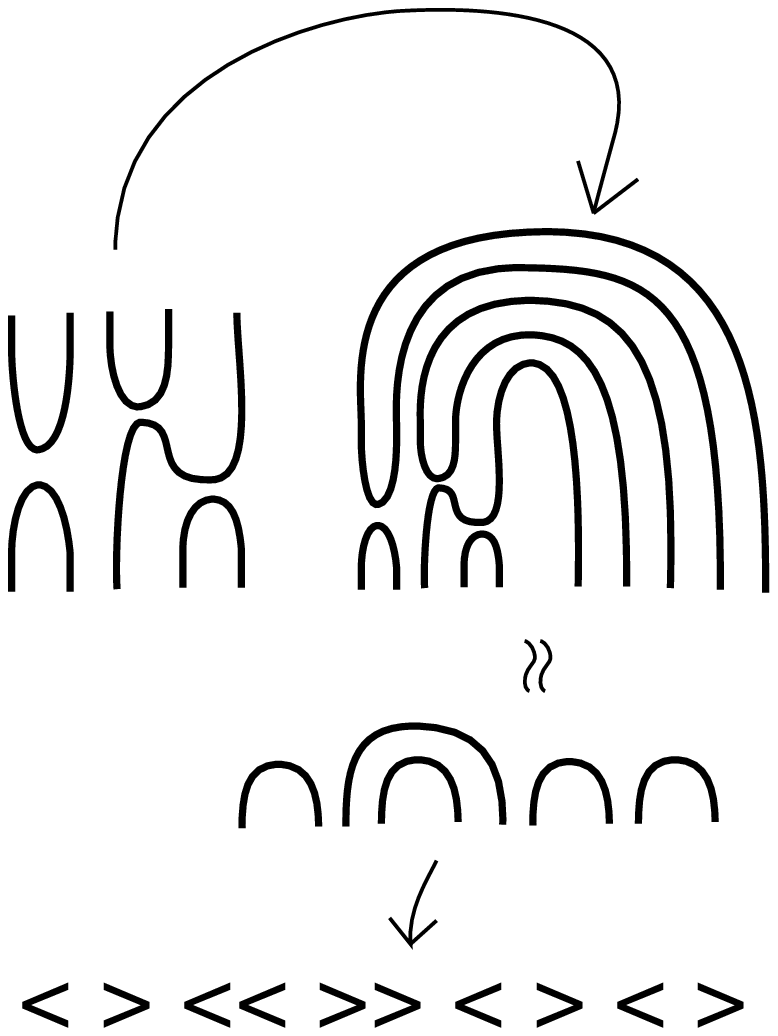}}
\vvvv

\begin{center}
{\bf Figure 50 - Creating Parentheses}
\end{center}
\vvv

One result of this reformulation of the Temperley-Lieb monoid  is that one can
rewrite the product structure in terms of operations on parentheses, getting an
interesting formal algebra that encodes properties of the topology of plane
curves \cite{KSpin}, \cite{spin}. In this section we will take the relationship
with parentheses as an excuse to indicate a further relationship with
non-associative products. \vspace{3mm}

First consider the abstract structure of non-associative products. They are
usually written in the forms:

$$(a*b)*c$$ $$a*(b*c)$$ $$((a*b)*c)*d$$ $$(a*(b*c))*d$$ $$a*((b*c)*d)$$
$$a*(b*(c*d))$$ $$(a*b)*(c*d)$$

\noindent In writing products in this standard manner, one sees the same
structure of parentheses occurring in different products, as in $((a*b)*c)*d$   and
 $(a*(b*c))*d.$  In each of these cases if we eliminated the algebraic literals and
the operation symbol $*$, we would be left with the parenthetical:  $(()).$
\vspace{3mm}

There is another way to write the products so that different products correspond
to different arrangements of parentheses.  To do this, rewrite the products in
the operator notation shown below:

$$a*b = a(b) \leadsto () $$ $$(a*b)*c = a(b)(c) \leadsto ()()$$ $$a*(b*c) = 
a(b(c)) \leadsto (())$$ $$((a*b)*c)*d =a(b)(c)(d)  \leadsto  ()()()$$
$$(a*(b*c))*d =a(b(c))(d)  \leadsto (())()$$ $$a*((b*c)*d) =a(b(c)(d))  \leadsto
(()())$$ $$a*(b*(c*d)) =a(b(c(d)))  \leadsto ((()))$$ $$(a*b)*(c*d) =a(b)(c(d)) 
\leadsto ()(())$$

In the operator notation, each product of $n+1$ terms is uniquely associated with
an expression using $n$ pairs of parentheses. \vspace{3mm}

Of course we can replace the expressions in parentheses in terms of nested caps
as we did in Figure 50. Once this is done, we notice the very interesting fact
that {\em re-associations} can be visualized in terms of ``handle-sliding" at the
cap level.  The meaning of this remark should be apparent to the reader from
Figure 51. \vspace{3mm}

\centerline{\includegraphics[scale=1.0]{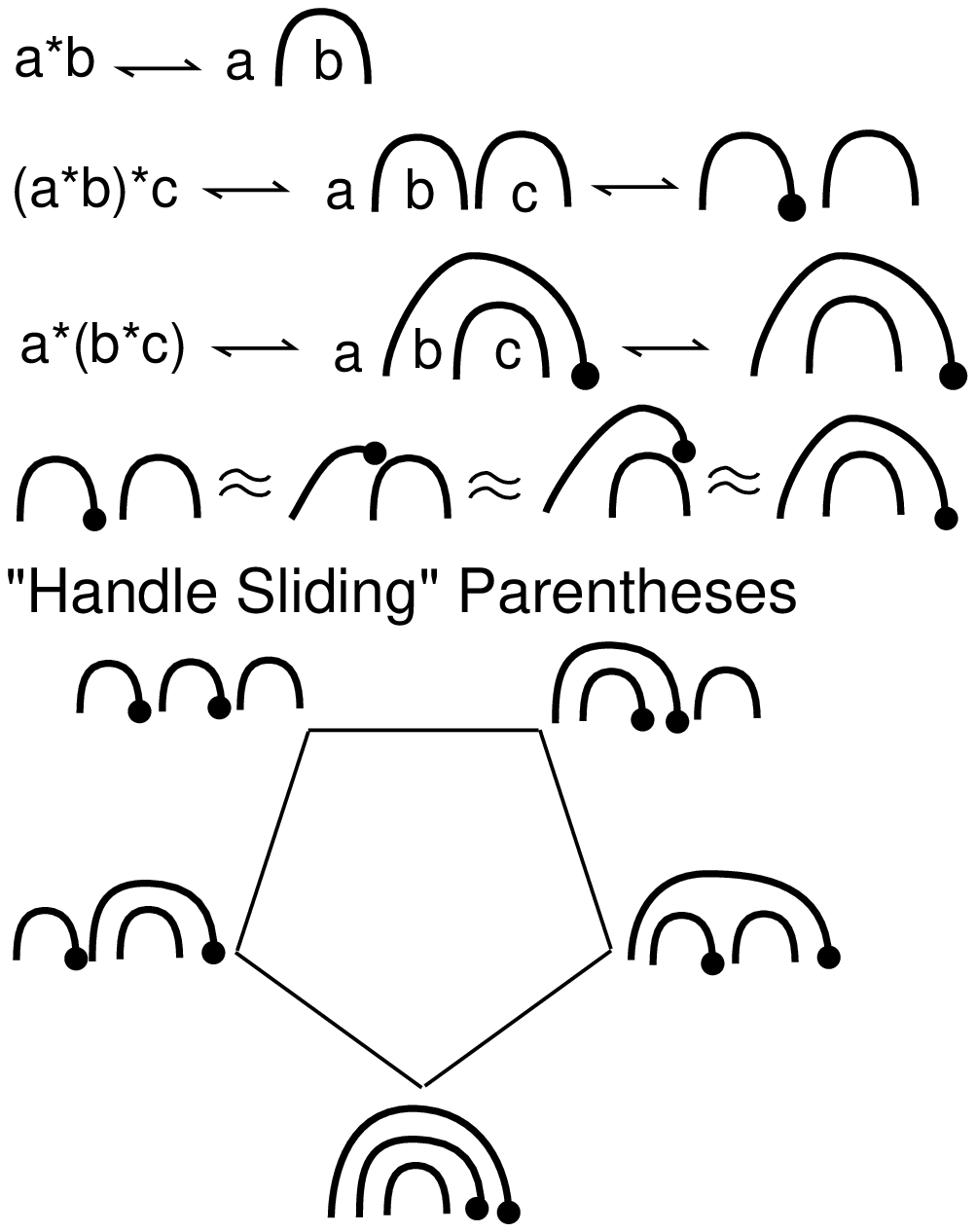}}
\vvvv

\begin{center}
{\bf Figure 51 - Pentagon}
\end{center}

In Figure 51 we illustrate the pentagon of re-associations of a product of four
terms in terms of sliding caps. Notice that we can do this sliding without
writing any algebraic literals by labelling the ``active" right feet of the caps
that do the sliding. In Figure 52 we show  the structure of the Stasheff
polyhedron with re-associations of five literals in the cap sliding formalism.
\vspace{3mm}

Cap sliding is a new formulation of a continuous background for these
re-association moves. There is a related continuous background in terms of
recoupling formalisms for trees. This formalism is intimately related to many
topics in topological quantum field theory. The interested reader will find more
about these points of view in \cite{spin}, \cite{CF}, \cite{CKY},  in the author's
lecture notes with Scott Carter and Masahico Saito \cite{CKS} and in his book
with Sostenes Lins \cite{TL}.
\newpage
\vspace*{-.2in}

\centerline{\includegraphics[scale=1.0]{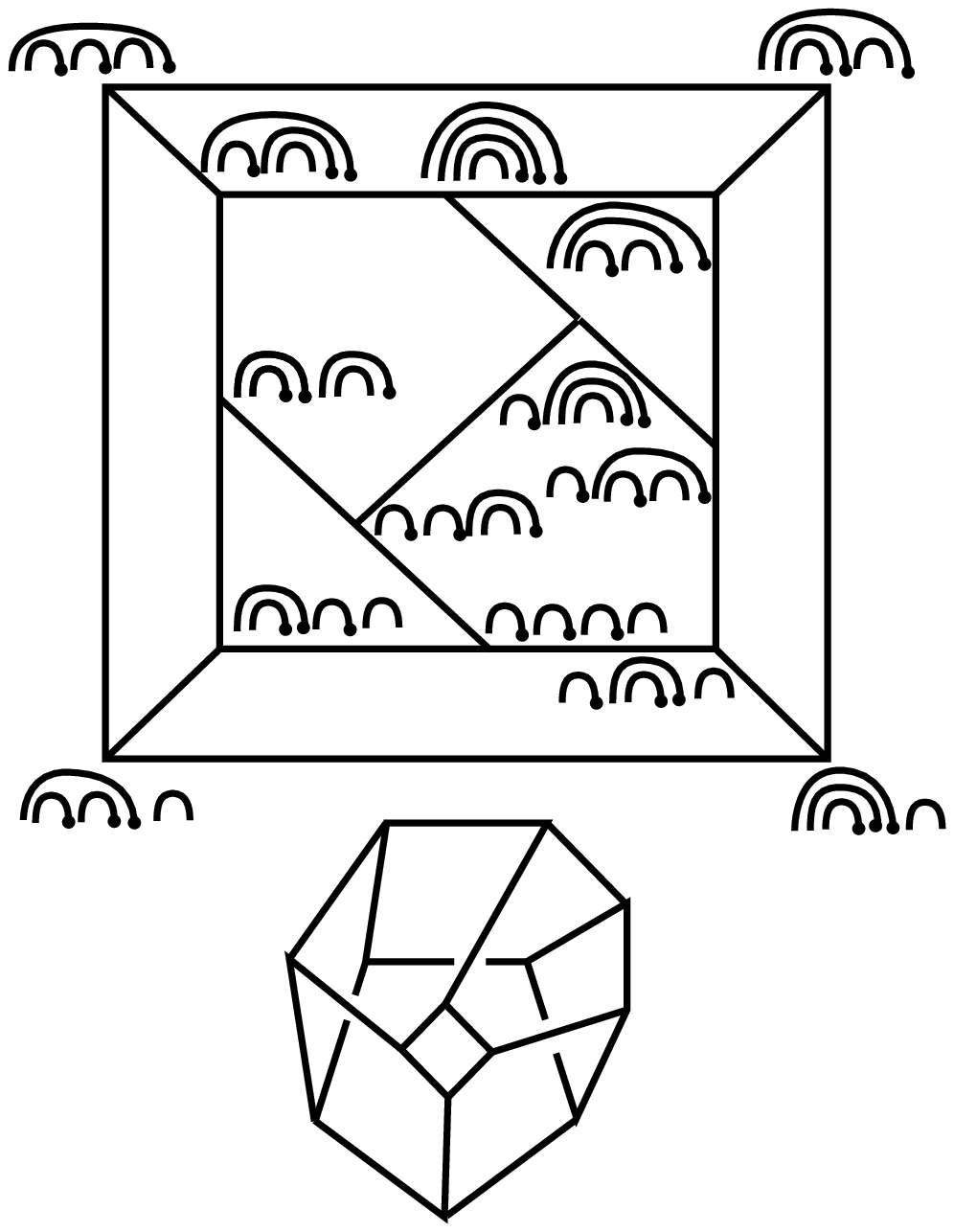}}
\vvvv

\begin{center}
{\bf Figure 52 - The Stasheff Polyhedron}
\end{center}
\vvv

\section{Virtual Knot Theory}

Knot theory
studies the embeddings of curves in three-dimensional space. Virtual knot theory studies the  embeddings of curves in thickened
surfaces of arbitrary genus, up to the addition and removal of empty handles from the surface.  Virtual knots
have a special  diagrammatic theory that makes handling them very similar to the handling of classical knot diagrams. In fact,
this diagrammatic theory simply involves adding a new type of crossing to the knot diagrams, a {\em virtual crossing} that is
neither under nor over. From a combinatorial point of view, the virtual crossings are artifacts of the representation of the
virtual knot or link in the plane. The extension of the Reidemeister moves that takes care of them respects this viewpoint.
A virtual crossing  (See Figure 53) is represented by two crossing arcs with a small circle
placed around the crossing point.
\bigbreak

Moves on virtual diagrams generalize the Reidemeister moves for classical knot and link
diagrams.  See Figure 53.  One can summarize the moves on virtual diagrams by saying that the classical crossings
interact with one another according to the usual Reidemeister moves. One adds the detour moves for consecutive sequences
of virtual crossings and this completes the description of the moves on virtual diagrams. It is a consequence of  moves
(B) and (C) in Figure 53 that an arc going through any consecutive sequence of virtual crossings
can be moved anywhere in the diagram keeping the endpoints fixed;
the places  where the moved arc crosses the diagram become new virtual crossings. This replacement
is the {\em detour move}. See Figure 53.1.
\bigbreak

One can generalize many structures in classical knot theory to the virtual domain, and use the virtual knots to test
the limits of classical problems such as  the question whether the Jones polynomial detects knots and the classical Poincar\'{e}
conjecture. Counterexamples to these conjectures exist in the virtual domain, and it is an open problem whether any of these
counterexamples are equivalent (by addition and subtraction of empty handles) to  classical knots and links. Virtual knot theory
is a significant domain to be investigated for its own sake and for a deeper understanding of classical knot theory.
\bigbreak

Another way to understand the meaning of virtual diagrams is to regard them as representatives for oriented Gauss codes
(Gauss diagrams) \cite{VKT,GPV}. Such codes do not always have planar realizations and an attempt to embed such a code in the plane
leads to the production of the virtual crossings. The detour move makes the particular choice of virtual crossings
irrelevant. Virtual equivalence is the same as the equivalence relation generated on the collection
of oriented Gauss codes modulo an abstract set of Reidemeister moves on the codes.
\bigbreak

One can consider {\em virtual braids}, generalizing the classical Artin Braid group. We shall not discuss this topic here, but refer the
reader to \cite{SVKT,DVK,KL3,Kamada,Ma0}.
\bigbreak

One intuition for virtual knot theory is the idea of a particle
moving in three-dimensional space in a trajectory that
occasionally disappears, and then reappears elsewhere. By
connecting the disappearance points and the reappearance points
with detour lines in the ambient space we get a picture of the
motion, but the detours, being artificial, must be treated as
subject to replacements.
\newpage

{\tt    \setlength{\unitlength}{0.92pt}
\begin{picture}(389,343)
\thicklines   \put(77,1){\framebox(252,95){}}
              \put(199,97){\framebox(189,245){}}
              \put(1,97){\framebox(191,240){}}
              \put(197,55){\vector(1,0){19}}
              \put(211,55){\vector(-1,0){17}}
              \put(278,131){\vector(1,0){19}}
              \put(292,131){\vector(-1,0){17}}
              \put(263,237){\vector(1,0){19}}
              \put(277,237){\vector(-1,0){17}}
              \put(265,305){\vector(1,0){19}}
              \put(279,305){\vector(-1,0){17}}
              \put(87,166){\vector(1,0){19}}
              \put(101,166){\vector(-1,0){17}}
              \put(64,240){\vector(1,0){19}}
              \put(78,240){\vector(-1,0){17}}
              \put(66,310){\vector(1,0){19}}
              \put(80,310){\vector(-1,0){17}}
              \put(332,2){\makebox(24,25){$C$}}
              \put(355,71){\makebox(23,26){$B$}}
              \put(31,75){\makebox(20,21){$A$}}
              \put(303,30){\circle{18}}
              \put(262,31){\circle{18}}
              \put(142,69){\circle{18}}
              \put(103,69){\circle{18}}
              \put(282,11){\line(1,1){39}}
              \put(241,52){\line(1,-1){41}}
              \put(121,91){\line(1,-1){40}}
              \put(83,51){\line(1,1){38}}
              \put(286,60){\line(1,1){34}}
              \put(242,12){\line(1,1){35}}
              \put(242,91){\line(1,-1){79}}
              \put(82,91){\line(1,-1){79}}
              \put(82,12){\line(1,1){35}}
              \put(127,57){\line(1,1){34}}
              \put(361,124){\circle{18}}
              \put(322,123){\circle{18}}
              \put(342,143){\circle{18}}
              \put(244,145){\circle{18}}
              \put(263,164){\circle{18}}
              \put(225,164){\circle{18}}
              \put(232,220){\circle{18}}
              \put(232,252){\circle{18}}
              \put(230,308){\circle{18}}
              \put(340,104){\line(1,1){39}}
              \put(299,146){\line(1,-1){41}}
              \put(302,104){\line(1,1){77}}
              \put(204,106){\line(1,1){79}}
              \put(241,183){\line(6,-5){43}}
              \put(205,146){\line(1,1){36}}
              \put(212,236){\line(5,-4){48}}
              \put(261,278){\line(-6,-5){49}}
              \put(251,328){\line(-6,-5){43}}
              \put(209,329){\line(1,-1){40}}
              \put(250,289){\line(0,1){40}}
              \put(296,327){\line(1,0){40}}
              \put(336,326){\line(0,-1){40}}
              \put(336,287){\line(-1,0){40}}
              \put(207,278){\line(1,-1){40}}
              \put(247,236){\line(-1,-1){39}}
              \put(338,236){\line(1,2){20}}
              \put(339,236){\line(1,-2){19}}
              \put(299,276){\line(1,-2){19}}
              \put(318,236){\line(-1,-2){20}}
              \put(204,186){\line(1,-1){79}}
              \put(302,183){\line(1,-1){79}}
              \put(184,147){\line(-1,-1){16}}
              \put(144,107){\line(1,1){16}}
              \put(144,107){\line(-1,1){15}}
              \put(104,147){\line(1,-1){18}}
              \put(106,185){\line(1,-1){79}}
              \put(106,106){\line(1,1){35}}
              \put(150,154){\line(1,1){34}}
              \put(87,149){\line(-1,1){15}}
              \put(46,189){\line(1,-1){17}}
              \put(31,174){\line(1,1){15}}
              \put(7,149){\line(1,1){16}}
              \put(52,155){\line(1,1){34}}
              \put(7,110){\line(1,1){35}}
              \put(7,189){\line(1,-1){79}}
              \put(124,241){\line(-1,-2){20}}
              \put(105,281){\line(1,-2){19}}
              \put(145,241){\line(1,-2){19}}
              \put(144,241){\line(1,2){20}}
              \put(36,220){\line(4,-3){24}}
              \put(11,242){\line(6,-5){16}}
              \put(27,254){\line(-4,-3){16}}
              \put(60,282){\line(-6,-5){24}}
              \put(47,240){\line(-1,-1){39}}
              \put(7,282){\line(1,-1){40}}
              \put(143,292){\line(-1,0){40}}
              \put(143,331){\line(0,-1){40}}
              \put(103,332){\line(1,0){40}}
              \put(20,310){\line(-1,-1){15}}
              \put(46,333){\line(-1,-1){16}}
              \put(46,293){\line(0,1){40}}
              \put(5,333){\line(1,-1){40}}
\end{picture}}

\begin{center}
{ \bf Figure 53 -- Generalized Reidemeister Moves for Virtuals} \end{center}
\vspace{3mm}

\begin{center}
{\tt    \setlength{\unitlength}{0.92pt}
\begin{picture}(482,223)
\thicklines   \put(10,172){\line(1,0){196}}
              \put(64,212){\line(0,-1){81}}
              \put(104,213){\line(0,-1){80}}
              \put(143,213){\line(0,-1){81}}
              \put(85,92){\line(0,-1){79}}
              \put(125,92){\line(0,-1){80}}
              \put(371,90){\line(0,-1){80}}
              \put(334,90){\line(0,-1){79}}
              \put(351,211){\line(0,-1){81}}
              \put(391,210){\line(0,-1){80}}
              \put(311,209){\line(0,-1){81}}
              \put(230,171){\line(1,0){41}}
              \put(271,171){\line(0,-1){121}}
              \put(271,50){\line(1,0){160}}
              \put(431,51){\line(0,1){120}}
              \put(431,171){\line(1,0){41}}
              \put(64,172){\circle{20}}
              \put(105,172){\circle{20}}
              \put(144,172){\circle{20}}
              \put(335,49){\circle{18}}
              \put(371,49){\circle{18}}
              \put(164,109){\vector(1,0){79}}
              \put(243,109){\vector(-1,0){80}}
\end{picture}}

{ \bf Figure 53.1 -- Detour Move}
\end{center}

\subsection{Flat Virtual Knots and Links}
Every classical knot or link diagram can be regarded as a $4$-regular plane graph with extra structure at the
nodes. This extra structure is usually indicated by the overcrossing and undercrossing conventions that give
instructions for constructing an embedding of the link in three-dimensional space from the diagram.  If we take the diagram
without this extra structure, it is the shadow of some link in three-dimensional space, but the weaving of that link is not
specified. It is well known that if one is allowed to apply the Reidemeister moves to such a shadow (without regard to the types
of crossing since they are not specified) then the shadow can be reduced to a disjoint union of circles. This reduction is
no longer true for virtual links. More precisely, let a {\em flat virtual diagram} be a diagram with virtual crossings as we have
described them and {\em flat crossings} consisting in undecorated nodes of the $4$-regular plane graph. Virtual crossings
are flat crossings that have been decorated by a small circle. Two flat virtual diagrams are {\em equivalent} if there is a
sequence of generalized flat Reidemeister moves (as illustrated in Figure 53) taking one to the other. A generalized
flat Reidemeister move is any move as shown in Figure 53, but one can ignore the overcrossing or undercrossing structure.
Note that in studying flat virtuals the rules for changing virtual crossings among themselves and the rules for changing
flat crossings among themselves are identical. However, detour moves as in Figure 53C are available for virtual crossings
with respect to flat crossings and not the other way around.
\v

We shall say that a virtual diagram {\em overlies} a flat diagram if the virtual diagram is obtained from the flat diagram by
choosing a crossing type for each flat crossing in the virtual diagram. To each virtual diagram $K$ there is an associated
flat diagram $F(K)$ that is obtained by forgetting the extra structure at the classical crossings in $K.$ Note that if $K$
is equivalent to $K'$ as virtual diagrams, then $F(K)$ is equivalent to $F(K')$ as flat virtual diagrams. Thus, if we can
show that $F(K)$ is not reducible to a disjoint union of circles, then it will follow that $K$ is a non-trivial virtual link.
\v

\begin{center}
{\tt    \setlength{\unitlength}{0.72pt}
\begin{picture}(507,286)
\thicklines   \put(364,1){\makebox(22,23){K}}
              \put(84,1){\makebox(24,23){D}}
              \put(73,242){\makebox(20,20){H}}
              \put(502,160){\line(-1,-1){38}}
              \put(383,40){\line(1,1){70}}
              \put(438,42){\line(-4,5){51}}
              \put(343,160){\line(5,-6){34}}
              \put(300,41){\line(0,1){34}}
              \put(301,180){\line(0,-1){93}}
              \put(262,160){\line(1,0){80}}
              \put(420,180){\line(-1,0){120}}
              \put(299,40){\line(1,0){82}}
              \put(501,200){\line(-1,0){120}}
              \put(381,200){\line(0,-1){120}}
              \put(382,80){\line(-1,0){122}}
              \put(261,82){\line(0,1){78}}
              \put(440,42){\line(1,0){65}}
              \put(421,179){\line(3,-5){83}}
              \put(301,158){\circle{20}}
              \put(382,178){\circle{20}}
              \put(502,200){\line(0,-1){40}}
              \put(416,71){\circle{20}}
              \put(158,73){\circle{20}}
              \put(125,42){\line(1,1){120}}
              \put(244,202){\line(0,-1){40}}
              \put(124,180){\circle{20}}
              \put(43,160){\circle{20}}
              \put(86,162){\line(4,-5){95}}
              \put(163,181){\line(3,-5){83}}
              \put(182,44){\line(1,0){65}}
              \put(3,84){\line(0,1){78}}
              \put(124,82){\line(-1,0){122}}
              \put(123,202){\line(0,-1){120}}
              \put(243,202){\line(-1,0){120}}
              \put(41,42){\line(1,0){82}}
              \put(42,180){\line(0,-1){138}}
              \put(162,182){\line(-1,0){120}}
              \put(4,162){\line(1,0){80}}
              \put(25,241){\oval(16,16)}
              \put(64,264){\line(-1,0){39}}
              \put(64,225){\line(0,1){39}}
              \put(24,225){\line(1,0){40}}
              \put(24,263){\line(0,-1){38}}
              \put(4,283){\line(0,-1){41}}
              \put(43,242){\line(-1,0){39}}
              \put(44,283){\line(0,-1){41}}
              \put(4,283){\line(1,0){40}}
\end{picture}}

{\bf Figure 54 -- Flats $H$ and $D$, and the knot $K.$} \end{center}
\newpage

Figure 54 illustrates an example of a flat virtual link $H.$ This link cannot be undone in the flat category because it
has an odd number of virtual crossings between its two components and each generalized Reidemeister move preserves the
parity of the number of virtual crossings between components.  Also illustrated in Figure 54 is a flat diagram $D$ and a
virtual knot $K$ that overlies it. This example is given in \cite{VKT}. The knot shown is undetectable by many
invariants (fundamental group, Jones polynomial) but it is knotted and this can be seen either by using a generalization
of the Alexander polynomial that we describe below, or by showing that the underlying diagram $D$ is a non-trivial flat
virtual knot using the filamentation invariant that is introduced in \cite{HR}. The filamentation invariant is a
combinatorial method that is sometimes successful in indentifying irreducible flat virtuals. At this writing we know
very few invariants of flat virtuals. The flat virtual diagrams present a strong challenge for the construction
of new invariants.  It is important to understand the structure of flat virtual knots and links. This structure
lies at the heart of the comparison of classical and virtual links. We wish to be able to determine when a given virtual
link is equivalent to a classcal link. The reducibility or irreducibility of the underlying flat diagram is the first
obstruction to such an equivalence.
\bigbreak

\subsection{Interpretation of Virtuals as Stable Classes of Links in  Thickened Surfaces}
There is a useful topological interpretation for this virtual
theory in terms of embeddings of links in thickened surfaces. See
\cite{VKT,DVK,KUP}.  Regard each virtual crossing as a shorthand
for a detour of one of the arcs in the crossing through a 1-handle
that has been attached to the 2-sphere of the original diagram. By
interpreting each virtual crossing in this way, we obtain an
embedding of a collection of circles into a thickened surface
$S_{g} \times R$ where $g$ is the number of virtual crossings in
the original diagram $L$, $S_{g}$ is a compact oriented surface of
genus $g$ and $R$ denotes the real line.  We say that two such
surface embeddings are {\em stably equivalent} if one can be
obtained from another by isotopy in the thickened surfaces,
homeomorphisms of the surfaces and the addition or subtraction of
empty handles.  Then we have the \smallbreak \noindent {\bf
Theorem \cite{VKT,DKT,KUP}.} {\em Two virtual link diagrams are
equivalent if and only if their correspondent surface embeddings
are stably equivalent.} \smallbreak \noindent \bigbreak

\noindent{\bf Virtual knots and links give rise to a host of
problems.} As we saw in the previous section, there are
non-trivial virtual knots with unit Jones polynomial.  Moreover,
there are non-trivial virtual knots with integer fundamental group
and trivial Jones polynomial. (The fundamental group is defined
combinatorially by generalizing the Wirtinger presentation.) These
phenomena underline the question of how planarity is involved in
the way the Jones polynomial appears to detect classical knots,
and that the relationship of the fundamental group (and peripheral
system) is a much deeper one than the surface combinatorics for
classical knots. It is possible to take the connected sum of two
trivial virtual diagrams and obtain a non-trivial virtual knot
(the Kishino knot).

Here long knots (or, equivalently $1-1$ tangles) come into play.
Having a knot, we can break it at some point and take its ends to
infinity (say, in a way that they coincide with the horizontal
axis line in the plane). One can study isotopy classes of such
knots. A well-known theorem says that in the classical case, knot
theory coincides with long knot theory. However, this is not the
case for virtual knots. By breaking the same virtual knot at
different points, one can obtain non-isotopic long knots \cite{FJK}.
Furthermore, even if the initial knot is trivial, the resulting
long knot may not be trivial. The ``connected sum'' of two
trivial virtual diagrams may not be trivial in the compact case.
The phenomenon occurs because these two knot diagrams may  be
non-trivial in the long category. It is sometimes more convenient
to consider long virtual knots rather than compact virtual knots,
since connected sum is well-defined for long knots. It is
important to construct long virtual knot invariants to see whether
long knots are trivial and whether they are
classical. One approach is to regard long knots as $1-1$ tangles
and use extensions of standard invariants  (fundamental group,
quandle, biquandle, etc). Another approach is to distinguish two
types of crossings: those having early undercrossing and those
having later undercrossing with respect to the orientation of the
long knot. The latter technique is described in \cite{MALONG}.
\bigbreak

Unlike classical knots, the connected sum of long knots is not commutative.
Thus, if we show that two long knots $K_{1}$ and $K_{2}$ do not
commute, then we see that they are different and both non-classical.

A typical example of such knots is the two parts of the Kishino
knot, see Figure 54.1.
\vvv

\begin{center}
{\tt    \setlength{\unitlength}{0.96pt}
\begin{picture}(195,239)
\thicklines   \put(111,11){\line(0,1){47}}
              \put(114,99){\line(0,1){47}}
              \put(95,10){\line(0,1){47}}
              \put(92,99){\line(0,1){47}}
              \put(111,58){\line(1,0){29}}
              \put(114,98){\line(1,0){29}}
              \put(65,58){\line(1,0){29}}
              \put(62,99){\line(1,0){29}}
              \put(24,59){\line(0,-1){12}}
              \put(163,79){\circle{20}}
              \put(45,78){\circle{20}}
              \put(24,60){\line(1,1){39}}
              \put(183,45){\line(0,1){17}}
              \put(133,43){\line(1,0){50}}
              \put(133,51){\line(0,-1){8}}
              \put(133,95){\line(0,-1){34}}
              \put(133,114){\line(0,-1){10}}
              \put(184,114){\line(-1,0){49}}
              \put(183,100){\line(0,1){13}}
              \put(77,46){\line(0,1){7}}
              \put(24,45){\line(1,0){53}}
              \put(77,95){\line(0,-1){33}}
              \put(77,113){\line(0,-1){10}}
              \put(25,113){\line(1,0){52}}
              \put(25,99){\line(0,1){13}}
              \put(141,59){\line(1,1){40}}
              \put(145,98){\line(1,-1){37}}
              \put(24,99){\line(1,-1){40}}
              \put(12,212){\line(1,-1){40}}
              \put(53,213){\line(1,0){41}}
              \put(96,211){\line(1,-1){37}}
              \put(52,171){\line(1,0){40}}
              \put(92,172){\line(1,1){40}}
              \put(13,212){\line(0,1){13}}
              \put(13,226){\line(1,0){52}}
              \put(65,226){\line(0,-1){10}}
              \put(65,208){\line(0,-1){33}}
              \put(12,171){\line(0,-1){12}}
              \put(12,158){\line(1,0){53}}
              \put(65,159){\line(0,1){7}}
              \put(134,213){\line(0,1){13}}
              \put(135,227){\line(-1,0){49}}
              \put(84,227){\line(0,-1){10}}
              \put(84,208){\line(0,-1){34}}
              \put(84,164){\line(0,-1){8}}
              \put(84,156){\line(1,0){50}}
              \put(134,158){\line(0,1){17}}
              \put(12,173){\line(1,1){39}}
              \put(33,191){\circle{20}}
              \put(114,192){\circle{20}}
              \put(139,175){\makebox(32,31){$K$}}
\end{picture}}
\vvv

{\bf Figure 54.1 -- Kishino and Parts}
\end{center}
\newpage

We have a natural map

$$\langle\mbox{Long virtual knots }\rangle\to \langle\mbox{Oriented compact virtual
knots}\rangle,$$

\noindent obtained by taking two infinite ends of the long knots together to
make a compact knot. This map is obviously well defined
and allows one to construct (weak) long virtual knot invariants
from classical invariants (by regarding compact knot
invariants as long knot invariants). There is no well-defined
inverse for this map. But, if we were able to construct the map
from compact virtual knots to long virtual knots, we could apply
the long techniques for the compact case. This map does have an
inverse for \underline{classical} knots. Thus, the long techniques
are applicable to classical (long) knots. It would be interesting
to obtain new classical invariants from it. The long category can
also be applied in the case of flat virtuals, where all problems
formulated above occur as well. \bigbreak

There are examples of virtual knots that are very difficult to
prove knotted, and there are infinitely many flat virtual diagrams
that appear to be irreducible, but we have no techniques to prove
it. How can one tell whether a virtual knot is classical? One can ask: Are
there non-trivial virtual knots whose connected sum is trivial?
The latter question cannot be shown by classical techniques, but it can be analyzed by using the surface interpretation
for virtuals. See \cite{MAN}.
\bigbreak

In respect to virtual knots, we are in the same
position as the compilers of the original knot tables. We are, in fact, developing tables.
At Sussex, tables of virtual knots are being constructed, and tables will appear in a book being written by Kauffman and
Manturov. The website ``Knotilus" has tables as do websites develped by Zinn-Justin and Zuber. The theory
of invariants of virtual knots, needs more development. Flat virtuals (whose
study is a generalization of the classification of immersions) are
a nearly unknown territory (but see \cite{HR,TURAEV}). The flat virtuals
provide the deepest challenge since we have very few invariants to
detect them. \bigbreak

\subsection{Jones Polynomial of Virtual Knots}

We use a generalization of the bracket state summation model for the Jones polynomial to extend it to virtual knots
and links.  We call a diagram in the plane
{\em purely virtual} if the only crossings in the diagram are virtual crossings. Each purely virtual diagram is equivalent by the
virtual moves to a disjoint collection of circles in the plane.
\bigbreak

Given a link diagram $K$, a state $S$ of this diagram is obtained by
choosing a smoothing for each crossing in the diagram and labelling that smoothing with either $A$ or $A^{-1}$
according to the convention that a counterclockwise rotation of the overcrossing line sweeps two
regions labelled $A$, and that a smoothing that connects the $A$ regions is labelled by the letter $A$. Then, given
a state $S$, one has the evaluation $<K|S>$ equal to the product of the labels at the smoothings, and one has the
evaluation $||S||$ equal to the number of loops in the state (the smoothings produce purely virtual diagrams).  One then has
the formula
$$<K> = \Sigma_{S}<K|S>d^{||S||-1}$$
where the summation runs over the states $S$ of the diagram $K$, and $d = -A^{2} - A^{-2}.$
This state summation is invariant under all classical and virtual moves except the first Reidemeister move.
The bracket polynomial is normalized to an
invariant $f_{K}(A)$ of all the moves by the formula  $f_{K}(A) = (-A^{3})^{-w(K)}<K>$ where $w(K)$ is the
writhe of the (now) oriented diagram $K$. The writhe is the sum of the orientation signs ($\pm 1)$ of the
crossings of the diagram. The Jones polynomial, $V_{K}(t),$ is given in terms of this model by the formula
$$V_{K}(t) = f_{K}(t^{-1/4}).$$
\noindent The reader should note that this definition is a direct generalization to the virtual category of the
state sum model for the original Jones polynomial \cite{state}. It is straightforward to verify the invariances stated above.
In this way one has the Jones polynomial for virtual knots and links.
\bigbreak

In terms of the interpretation of virtual knots as stabilized classes of embeddings of circles into thickened surfaces,
our definition coincides with the simplest version of the Jones polynomial for links in thickened surfaces. In that
version one counts all the loops in a state the same way, with no regard for their isotopy class in the surface.
It is this equal treatment that makes the invariance under handle stabilization work. With this generalized version of the
Jones polynomial, one has again the problem of finding a geometric/topological interpretation of this invariant. There is
no fully satisfactory topological interpretation of the original Jones polynomial and the problem is inherited by this
generalization.
\bigbreak

\noindent We have \cite{DVK} the
\smallbreak
\noindent
{\bf Theorem.} {\em To each non-trivial
classical knot diagram of one component $K$ there is a corresponding  non-trivial virtual knot diagram $Virt(K)$ with unit
Jones polynomial.}
\bigbreak

This Theorem is a key ingredient in the problems involving virtual knots. Here is a sketch of its proof.
The proof uses two invariants of classical knots and links that generalize to arbitrary virtual knots and links.
These invariants are the {\em Jones polynomial} and the {\em involutory quandle} denoted by the notation
$IQ(K)$ for a knot or link $K.$
\bigbreak

Given a
crossing $i$ in a link diagram, we define $s(i)$ to be the result of {\em switching} that crossing so that the undercrossing arc
becomes an overcrossing arc and vice versa. We also define the {\em virtualization}
$v(i)$ of the crossing by the local replacement indicated in Figure 55. In this Figure we illustrate how in the
virtualization of the crossing the  original crossing is replaced by a crossing that is flanked by two virtual crossings.
\vv

Suppose that $K$ is a (virtual or classical) diagram with a classical crossing labeled $i.$  Let $K^{v(i)}$ be the diagram
obtained from $K$ by virtualizing the crossing $i$ while leaving the rest of the diagram just as before. Let $K^{s(i)}$ be
the diagram obtained from $K$ by switching the crossing $i$ while leaving the rest of the diagram just as before. Then it
follows directly from the definition of the Jones polynomial that $$V_{K^{s(i)}}(t) = V_{K^{v(i)}}(t).$$
\noindent As far as the Jones
polynomial is concerned, switching a crossing and virtualizing a crossing look the same.
\vv

The involutory quandle \cite{KNOTS} is an algebraic invariant
equivalent to the fundamental group of the double branched cover of a knot or link in the classical case. In this algebraic
system one associates a generator of the algebra $IQ(K)$ to each arc of the diagram $K$ and there is a relation of the form
$c = ab$ at each crossing, where $ab$ denotes the (non-associative) algebra product of $a$ and $b$ in $IQ(K).$ 
See Figure 56.
In this Figure we have illustrated through the local relations the fact that  $$IQ(K^{v(i)}) = IQ(K).$$
\noindent As far the involutory quandle is concerned, the original crossing and the virtualized crossing look the same.
\vv

If a classical knot is actually knotted, then its involutory quandle is non-trivial \cite{W}. Hence if we start
with a non-trivial classical knot, we can virtualize any subset of its crossings to obtain a virtual knot that is still
non-trivial. There is a subset $A$ of the crossings of a classical knot $K$ such that the knot $SK$ obtained by
switching these crossings is an unknot.  Let $Virt(K)$ denote the virtual diagram obtained from $A$ by virtualizing
the crossings in the subset $A.$  By the above discussion the Jones polynomial of $Virt(K)$ is the same
as the Jones polynomial of $SK$, and this is $1$ since $SK$ is unknotted. On the other hand, the $IQ$ of $Virt(K)$ is the
same as the $IQ$ of $K$, and hence if $K$ is knotted, then so is $Virt(K).$   We have shown that $Virt(K)$ is a non-trivial
virtual knot with unit Jones polynomial.  This completes the proof of the Theorem.
\vvvv

\begin{center}
{\tt    \setlength{\unitlength}{0.92pt}
\begin{picture}(326,162)
\thinlines    \put(190,141){\line(0,-1){82}}
              \put(242,34){\makebox(41,41){s(i)}}
              \put(170,10){\makebox(41,41){v(i)}}
              \put(1,35){\makebox(41,42){i}}
              \put(292,81){\line(1,0){33}}
              \put(246,81){\line(1,0){35}}
              \put(286,161){\line(0,-1){160}}
              \put(143,81){\circle{16}}
              \put(191,80){\circle{16}}
              \put(190,59){\line(-1,0){25}}
              \put(167,22){\line(0,-1){20}}
              \put(144,22){\line(1,0){22}}
              \put(143,102){\line(0,-1){80}}
              \put(165,103){\line(-1,0){22}}
              \put(165,87){\line(0,1){15}}
              \put(165,59){\line(0,1){17}}
              \put(167,141){\line(1,0){23}}
              \put(166,161){\line(0,-1){20}}
              \put(126,81){\line(1,0){80}}
              \put(46,74){\line(0,-1){73}}
              \put(46,161){\line(0,-1){74}}
              \put(6,81){\line(1,0){80}}
\end{picture}}
\vvvv

{\bf Figure 55 --  Switching and Virtualizing a Crossing} \end{center}
\newpage

\begin{center}
{\tt    \setlength{\unitlength}{0.92pt}
\begin{picture}(214,164)
\thinlines    \put(197,64){\makebox(16,16){b}}
              \put(136,144){\makebox(22,18){c =}}
              \put(17,143){\makebox(22,18){c =}}
              \put(165,143){\makebox(19,20){ab}}
              \put(113,64){\makebox(16,16){b}}
              \put(67,64){\makebox(16,16){b}}
              \put(43,1){\makebox(17,19){a}}
              \put(47,142){\makebox(19,20){ab}}
              \put(1,63){\makebox(16,16){b}}
              \put(165,1){\makebox(17,19){a}}
              \put(185,142){\line(0,-1){82}}
              \put(138,82){\circle{16}}
              \put(186,81){\circle{16}}
              \put(185,60){\line(-1,0){25}}
              \put(162,23){\line(0,-1){20}}
              \put(139,23){\line(1,0){22}}
              \put(138,103){\line(0,-1){80}}
              \put(160,104){\line(-1,0){22}}
              \put(160,88){\line(0,1){15}}
              \put(160,60){\line(0,1){17}}
              \put(162,142){\line(1,0){23}}
              \put(161,162){\line(0,-1){20}}
              \put(121,82){\line(1,0){80}}
              \put(41,75){\line(0,-1){73}}
              \put(41,162){\line(0,-1){74}}
              \put(1,82){\line(1,0){80}}
\end{picture}}
\vvvv

{\bf Figure 56 --   $IQ(Virt(K) = IQ(K)$} \end{center}
\vvv

If there exists a classical knot with unit Jones polynomial, then
one of the knots $Virt(K)$ produced by this theorem may be
equivalent to  a classical knot.  It is an intricate task to
verify that specific examples of $Virt(K)$ are not classical. This
has led to an investigation of new invariants for virtual knots.
In this investigation a number of issues appear. One can examine
the combinatorial generalization of the fundamental group (or
quandle) of the virtual knot and sometimes one can prove by pure
algebra that the resulting group is not classical. This is related
to observations by Silver and Williams \cite{SW}, Manturov
\cite{MAPOLY, MAPOLY'} and by Satoh \cite{SATOH} showing that the
fundamental group of a virtual knot can be interpreted as the
fundamental group of the complement of a torus embedded in four-dimensional
Euclidean space. A very fruitful line of new
invariants comes about by examining a generalization of the
fundamental group or quandle that we call the {\em biquandle} of
the virtual knot. The biquandle is discussed in the next Section.
Invariants of flat knots (when one has them) are useful in this
regard. If we can verify that the flat knot $F(Virt(K))$ is
non-trivial, then $Virt(K)$ is non-classical. In this way the
search for classical knots with unit Jones polynomial expands to
the exploration of the structure of the infinite collection of
virtual knots with unit Jones polynomial.
\v

Another way of putting this theorem is as follows: In the arena of
knots in thickened surfaces there are many examples of knots with
unit Jones polynomial. Might one of these be equivalent via handle
stabilization to a classical knot? In \cite{KUP} Kuperberg shows
the uniqueness of the  embedding of minimal genus in the stable
class for a given virtual link. The minimal embedding genus can be
strictly less than the number of virtual crossings in a diagram
for the link.  There are many problems associated with this
phenomenon.
\v

There is a generalization of the Jones polynomial that
involves surface representation of virtual knots. See \cite{Dye,DKMin, MACURVES,MACURVES'}.
These invariants essentially use the fact that the Jones polynomial can be extended to knots in thickened surfaces by keeping
track of the isotopy classes of the loops in the state summation  for this polynomial. In the approach of Dye and Kauffman, one
uses this generalized polynomial directly. In the approach of Manturov, a polynomial invariant is defined using the stabilization
description of the virtual knots.
\newpage

\subsection{Biquandles}
In this section we give a sketch of some recent approaches to invariants of virtual knots and links.
\bigbreak

A {\em biquandle} \cite{CS,DVK,FJK,GEN,FB,FBu} is an algebra with 4 binary
operations written $a^b, a_b, a^{\overline{b}}, a_{\overline{b}}$ together
with some relations which we will indicate below. The {\em
fundamental} biquandle is associated with a link diagram and is
invariant under the generalized Reidemeister moves for virtual knots
and links.  The operations in this algebra are motivated by the
formation of labels for the edges of the diagram.  View Figure 57.  In
this Figure we have shown the format for the operations in a
biquandle. The overcrossing arc has two labels, one on each side of
the crossing. There is an algebra element labeling each {\em edge} of
the diagram.  An edge of the diagram corresponds to an edge of the
underlying plane graph of that diagram.  \bigbreak

Let the edges oriented toward a crossing in a diagram be called the {\em input} edges for the crossing, and
the edges oriented away from the crossing be called the {\em output} edges for the crossing.
Let $a$ and $b$ be the input edges for a positive crossing, with $a$ the label of the undercrossing
input and $b$ the label on the overcrossing input. In the biquandle, we label the undercrossing
output by $$c=a^{b},$$ while the overcrossing output is labeled
$$d= b_{a}.$$ \smallbreak

The labelling for the negative crossing is similar using the other two operations.

\noindent To form the fundamental biquandle, $BQ(K)$, we take one generator for each edge of the diagram and two
relations at each crossing (as described above).
\vvvv\vvv

\begin{center}
{\tt    \setlength{\unitlength}{0.96pt}
\begin{picture}(280,163)
\thicklines   \put(131,162){\line(0,-1){159}}
              \put(1,1){\framebox(278,161){}}
              \put(231,104){\makebox(21,19){$a^{\o{b}} = a\, \UL{\it b}$}}
              \put(167,59){\makebox(20,15){$b_{\o{a}} = b\, \LL{\it a}$}}
              \put(233,62){\makebox(20,16){$b$}}
              \put(212,38){\makebox(18,17){$a$}}
              \put(73,106){\makebox(23,23){$a^{b} = a\, \UR{\it b}$}}
              \put(77,56){\makebox(32,22){$b_{a} = b\, \LR{\it a}$}}
              \put(10,86){\makebox(20,16){$b$}}
              \put(52,43){\makebox(18,17){$a$}}
              \put(211,89){\vector(0,1){33}}
              \put(211,43){\vector(0,1){32}}
              \put(250,82){\vector(-1,0){78}}
              \put(51,90){\vector(0,1){33}}
              \put(51,43){\vector(0,1){34}}
              \put(11,83){\vector(1,0){80}}
\end{picture}}
\vvvv\vv

{\bf Figure 57 -- Biquandle Relations at a Crossing} \end{center}
\newpage

\noindent Another way to write this formalism for the biquandle is as follows:

$$a^{b} = a\, \UR{\it b}$$
$$a_{b} = a\, \LR{\it b}$$
$$a^{\o{b}} = a\, \UL{\it b}$$
$$a_{\o{b}} = a\, \LL{\it b}.$$

\noindent We call this the {\em operator formalism} for the biquandle.
\bigbreak

\noindent These considerations lead to the following definition.
\smallbreak

\noindent {\bf Definition.}  A {\em biquandle} $B$ is a set with four binary operations indicated 
above:  $a^{b} \,\mbox{,} \, a^{\o{b}} \, \mbox{,} \,  a_{b} \,\mbox{,} \, a_{\o{b}}.$ We shall refer to the operations with barred
variables as the {\em left} operations and the operations without barred variables as the {\em right} operations. The biquandle is
closed under these operations and the following axioms are satisfied:

\begin{enumerate}
\item Given an element $a$ in $B$, then there exists an $x$ in the biquandle such that $x=a_{x}$ and
$a = x^{a}.$ There also exists a $y$ in the biquandle such that $y=a^{\o{y}}$ and
$a = y_{\o{a}}.$

\item   For any elements $a$ and $b$ in $B$ we have

$$a = a^{b \o{b_{a}}}  \quad \mbox{and} \quad  b= b_{a \o{a^{b}}} \quad \mbox{and}$$

$$a = a^{\o{b}b_{\o{a}}}  \quad \mbox{and} \quad  b= b_{\o{a} a^{\o{b}}}.$$

\item 
\noindent Given elements $a$ and $b$ in $B$ then there exist elements $x, y,
z, t$ such that $x_{b}=a$, $y^{\overline{a}}= b$, $b^x=y$,
$a_{\overline{y}}=x$ and $t^a=b$, $a_t=z$, $z_{\overline{b}}=a$,
$b^{\overline{z}}=t$. The biquandle is called {\em strong} if $x, y,
z, t$ are uniquely defined and we then write $x=a_{b^{-1}},
y=b^{\overline{a}^{-1}}, t=b^{a^{-1}}, z=a_{\overline{b}^{-1}}$,
reflecting the inversive nature of the elements.

\item For any $a$ , $b$ , $c$ in $B$ the following equations hold and the same equations hold when all right operations are
replaced in these equations by left operations.

$$a^{b c} = a^{c_{b} b^{c}} \mbox{,} \quad c_{b a} = c_{a^{b} b_{a}} \mbox{,} \quad (b_{a})^{c_{a^{b}}} = (b^{c})_{a^{c_{b}}}.$$

\end{enumerate}

These axioms are transcriptions of the Reidemeister moves.The first axiom transcribes the first Reidemeister move.
The second axiom transcribes the directly oriented second Reidemeister move.
The third axiom transcribes the reverse oriented Reidemeister move. The fourth axiom transcribes the third Reidemeister move.
Much more work is needed in exploring these  algebras and their applications to knot theory.
\bigbreak

We may simplify the appearance of these conditions by defining
$$S(a,b)=(b_a,a^b),\quad
\overline{S}(a,b)=(b^{\overline{a}},a_{\overline{b}})$$ and in the case
of a strong biquandle,$$S^+_-(a,b)=(b^{a_{b^{-1}}},a_{b^{-1}}),\quad
S^-_+(a,b)=(b^{a^{-1}},a_{b^{a^{-1}}})$$ and 
$${\overline{S}}^{\lower.5ex\hbox{\ $\scriptstyle+$}}_{\ -}(a,b)=
(b_{\ \overline{a^{\overline{b}^{-1}}}}\ , \ a^{\overline{b}^{-1}})=
(b_{\ \overline{a^{b_{a^{-1}}}}} \ , \ a^{b_{a^{-1}}})$$ \ 
and \ $${\overline{S}}^{\lower.5ex\hbox{\ $\scriptstyle-$}}_{\ +}(a,b)=
(b_{\ \overline{a}^{-1}} \ , \ a^{\overline{b_{\overline{a}^{-1}}}})=
(b_{{a}^{b^{-1}}} \ , \ a^{\overline{b_{{a}^{b^{-1}}}}})$$

which we call the {\em sideways} operators. The conditions
then reduce to $$S\overline{S}=\overline{S}S=1,$$ $$ (S\times 1)
(1\times S) (S\times 1) = (1\times S) (S\times 1) (1\times S)$$
$$\overline{S}^-_+S^+_-=S^-_+\overline{S}^+_-=1$$ and finally all
the sideways operators leave the diagonal $$\Delta=\{(a,a)|a\in
X\}$$ invariant.

\subsection{The Alexander Biquandle}

It is not hard to see that
the following equations in a module over $Z[s,s^{-1},t,t^{-1}]$ give a biquandle
structure.

$$a^b=a\,\UR{\it b} = ta + (1-st)b \, \mbox{,} \quad a_b=a\,\LR{\it b} = sa$$
$$a^{\overline{b}}=a\,\UL{\it b} = t^{-1}a + (1-s^{-1}t^{-1})b \, \mbox{,} 
\quad a_{\overline{b}}=a\,\LL{\it b} = s^{-1}a.$$

\noindent We shall refer to this
structure, with the equations given above, as the {\em Alexander Biquandle}.
\vspace{3mm}

Just as one can define the Alexander Module of a classical knot, we have the Alexander Biquandle of
a virtual knot or link, obtained by taking one generator for each {\em edge} of the projected graph of the  knot diagram and
taking the module relations in the above linear form. Let $ABQ(K)$ denote this module structure for an
oriented link $K$. That is, $ABQ(K)$ is the module generated by the edges of the diagram, factored by the submodule
generated by the relations. This module then has a biquandle structure specified by the operations defined above for an
Alexander Biquandle.
\newpage
\vspace*{-.2in}

The determinant of the matrix of relations obtained from the crossings of a diagram
gives a polynomial invariant (up to multiplication by $\pm s^{i}t^{j}$ for integers $i$ and $j$)
of knots and links that we denote by $G_{K}(s,t)$ and call the {\em generalized Alexander polynomial}.
{\bf This polynomial vanishes on classical knots, but is remarkably successful at detecting virtual knots and
links.} In fact $G_{K}(s,t)$ is the same as the polynomial invariant of virtuals
of Sawollek \cite{SAW} and defined by an alternative method by Silver and Williams \cite{SW}.
It is a reformulation of the invariant for knots in surfaces due to the principal investigator, Jaeger and
Saleur \cite{JKS,KS}.
\vspace{3mm}

We end this discussion of the Alexander Biquandle with two
examples that show clearly its limitations. View Figure 58. In this
Figure we illustrate two diagrams labeled $K$ and $KI.$ It is not
hard to calculate that both $G_{K}(s,t)$ and $G_{KI}(s,t)$ are
equal to zero. However, The Alexander Biquandle of $K$ is
non-trivial -- it is isomorphic to the free module over $Z[s,
s^{-1},t, t^{-1}]$ generated by elements $a$ and $b$ subject to
the relation $(s^{-1} - t -1)(a-b) =0.$ Thus $K$ represents a
non-trivial virtual knot. This shows that it is possible for a
non-trivial virtual diagram to be a connected sum of two trivial
virtual diagrams. However, the diagram $KI$ has a trivial
Alexander Biquandle.  In fact the diagram $KI$, discovered by Kishino \cite{P}, is
now known to be knotted
and its general biquandle is non-trivial. The Kishino diagram has
been shown non-trivial by a calculation of the three-strand Jones
polynomial \cite{KiSa},  by the surface bracket polynomial of Dye and Kauffman
\cite{Dye,DKMin}, by the $\Xi$-polynomial (the surface
generalization of the Jones polynomial of Manturov
\cite{MACURVES}, and its biquandle has been shown to be non-trivial by
a quaternionic biquandle representation \cite{FB} of Fenn and Bartholomew which
we will now briefly describe.
\vv

Referring back to the previous section define the linear biquandle by
$$S=\pmatrix{
1+i&jt\cr -jt^{-1}&1+i\cr},$$
where $i,j$ have their usual meanings as quaternions and $t$ is a central variable. Let $R$ denote the
ring which they determine. Then as in the Alexander case considered above, for each diagram there is a
square presentation of an $R$-module. We can take the (Study) determinant of the presentation matrix.
In the case of the Kishino knot this is zero. However the greatest common divisor of the codimension 1 
determinants is $2+5t^2+2t^4$ showing that this knot is not classical.
\vvv\vv

\begin{center} {\tt    \setlength{\unitlength}{0.92pt}
\begin{picture}(264,109)
\thicklines   \put(186,1){\makebox(29,29){$KI$}}
              \put(50,1){\makebox(32,31){$K$}}
              \put(241,70){\circle{20}}
              \put(159,70){\circle{20}}
              \put(105,71){\circle{20}}
              \put(24,70){\circle{20}}
              \put(219,49){\line(-1,0){9}}
              \put(223,89){\line(-1,0){13}}
              \put(180,50){\line(1,0){23}}
              \put(179,90){\line(1,0){24}}
              \put(205,35){\line(1,0){55}}
              \put(206,105){\line(0,-1){70}}
              \put(261,105){\line(-1,0){55}}
              \put(139,90){\line(1,-1){40}}
              \put(223,89){\line(1,-1){37}}
              \put(219,50){\line(1,1){40}}
              \put(140,90){\line(0,1){13}}
              \put(140,104){\line(1,0){52}}
              \put(192,104){\line(0,-1){10}}
              \put(192,86){\line(0,-1){33}}
              \put(139,49){\line(0,-1){12}}
              \put(139,36){\line(1,0){53}}
              \put(192,37){\line(0,1){7}}
              \put(261,91){\line(0,1){13}}
              \put(259,36){\line(0,1){17}}
              \put(139,51){\line(1,1){39}}
              \put(3,52){\line(1,1){39}}
              \put(125,37){\line(0,1){17}}
              \put(75,35){\line(1,0){50}}
              \put(75,43){\line(0,-1){8}}
              \put(75,87){\line(0,-1){34}}
              \put(75,106){\line(0,-1){10}}
              \put(126,106){\line(-1,0){49}}
              \put(125,92){\line(0,1){13}}
              \put(56,38){\line(0,1){7}}
              \put(3,37){\line(1,0){53}}
              \put(3,50){\line(0,-1){12}}
              \put(56,87){\line(0,-1){33}}
              \put(56,105){\line(0,-1){10}}
              \put(4,105){\line(1,0){52}}
              \put(4,91){\line(0,1){13}}
              \put(83,51){\line(1,1){40}}
              \put(43,50){\line(1,0){40}}
              \put(87,90){\line(1,-1){37}}
              \put(44,92){\line(1,0){41}}
              \put(3,91){\line(1,-1){40}}
\end{picture}}
\end{center}

\begin{center} {\bf  Figure 58 -- The Knot $K$ and the Kishino Diagram $KI$ } \end{center}
\newpage

\subsubsection{Virtual quandles}

There is another generalization of quandle \cite{Ma1} by means of
which one can obtain the same polynomial as in \cite{SW,SAW} from the
other point of view \cite{MAPOLY,MAPOLY'}. Namely, the formalism is the
same as in the case of quandles at classical crossings but one adds
a special structure at virtual crossings.
The fact that these approaches give the same result in the linear
case was proved recently by Roger Fenn and Andrew Bartholomew.
\bigbreak

Virtual quandles (as well as biquandles) yield generalizations of
the fundamental group and some other invariants. Also, virtual
biquandles admit a generalization for multi-variable polynomials
in the case of multicomponent links, see \cite{MAPOLY}. One can
extend these definitions by bringing together the virtual quandle
(at virtual crossings) and the biquandle (at classical crossings)
to obtain what is called a {\em virtual biquandle}; this work is
now in process, \cite{KaMa}. \bigbreak

\subsection{A Quantum Model for $G_{K}(s,t),$ Oriented and Bi-oriented Quantum Algebras}

We can understand the structure of the invariant $G_{K}(s,t)$ by rewriting it as a quantum
invariant and then analysing its state summation.  The quantum model for this invariant is obtained in a fashion analogous to the construction of a
quantum model of the Alexander polynomial in \cite{KS,JKS}. The strategy in those papers was to take the basic
two-dimensional matrix of the Burau representation, view it as a linear transformation  $T: V \longrightarrow V$ on
a two-dimensional module $V$, and them take the induced linear transformation
$\hat{T}: \Lambda^{*}V \longrightarrow  \Lambda^{*}V$ on the exterior algebra of $V$. This gives a transformation on a
four-dimensional module that is a solution to the Yang-Baxter equation.  This solution of the Yang-Baxter equation
then becomes the building block for the corresponding quantum invariant.  In the present instance, we have a generalization
of the Burau representation, and this same procedure can be applied to it.
\bigbreak

\noindent The normalized state summation $Z(K)$ obtained by the above process satisfies a skein relation that is just like that of the Conway polynomial:
$Z(K_{+}) - Z(K_{-}) = zZ(K_{0}).$
The basic result behind the correspondence of $G_{K}(s,t)$ and $Z(K)$ is the
 \bigbreak

 \noindent {\bf Theorem \cite{GEN}.} {\em For a (virtual) link $K$, the invariants $Z(K)(\sigma = \sqrt{s}, \tau = 1/\sqrt{t})$ and $G_{K}(s,t)$ are equal
up to
 a multiple of $\pm s^{n}t^{m}$ for integers $n$ and $m$ (this being the well-definedness criterion for $G$).}
 \bigbreak

It is the purpose of this section to place our work with the generalized Alexander polynomial in a context
of bi-oriented quantum algebras and to
introduce the concept of an oriented quantum algebra. In \cite{KRO,KRCAT} Kauffman and
Radford introduce the concept and show that
{\em oriented quantum algebras encapsulate the notion of an oriented quantum link invariant.}
\smallbreak

An {\em oriented quantum algebra} $(A, \rho, D, U)$ is an abstract
model for an oriented quantum invariant of classical links \cite{KRO,KRCAT}. This model is based
on a solution to the Yang-Baxter equation.   The
definition of an oriented quantum algebra is as follows:   We are
given an algebra $A$ over a base ring $k$, an invertible solution
$\rho$ in $A \otimes A$ of the Yang-Baxter equation (in the
algebraic formulation of this equation -- differing from a braiding operator by a transposition), and commuting automorphisms
$U,D:A \longrightarrow A$  of the algebra, such that

$$(U \otimes U)\rho = \rho,$$

$$(D \otimes D)\rho = \rho,$$

$$[(1_{A} \otimes U)\rho)][(D \otimes 1_{A^{op}})\rho^{-1}]
= 1_{A \otimes A^{op}},$$

\noindent and

$$[(D \otimes 1_{A^{op}})\rho^{-1}][(1_{A} \otimes U)\rho)]
= 1_{A \otimes A^{op}}.$$
\v

\noindent The last two equations say that $[(1_{A} \otimes U)\rho)]$ and $[(D \otimes 1_{A^{op}})\rho^{-1}]$
are inverses in the algebra $A \otimes A^{op}$ where $A^{op}$ denotes the opposite algebra.
\vspace{3mm}

When $U=D=T$, then $A$ is said to be {\em balanced}.
In the case where $D$ is the identity mapping, we call the
oriented quantum algebra {\em standard}.  In  \cite{KRCAT} we show that
{\em the invariants defined by Reshetikhin and Turaev (associated with a
quasi-triangular Hopf algebra) arise from standard oriented quantum algebras.}
It is an interesting structural feature of algebras that we have elsewhere
\cite{GAUSS} called {\em quantum algebras} (generalizations of quasi-triangular Hopf algebras)
that they give rise to standard oriented quantum algebras. 
We are continuing research on the relationships of quantum link invariants and
oriented quantum algebras. In particular we are working on the reformulation of 
existing invariants such as the Links-Gould invariant \cite{KDL}
that are admittedly powerful but need a deeper understanding both topologically and algebraically.
\bigbreak

\noindent We now extend the concept of oriented quantum algebra  by adding a second solution to
the Yang-Baxter equation $\gamma$ that will take the role of the virtual crossing.
\smallbreak

\noindent {\bf Definition.} A {\em bi-oriented quantum algebra} is a quintuple  $(A, \rho, \gamma, D , U)$ such that
$(A, \rho, D, U)$ and $(A,\gamma, D, U)$ are oriented quantum algebras and $\gamma$ satisfies the following properties:

\begin{enumerate}
\item $\gamma_{12}\gamma_{21} = 1_{A \otimes A}.$ (This is the equivalent to the statement that the
braiding operator corresponding to $\gamma$ is its own inverse.)

\item Mixed identities involving $\rho$ and $\gamma$ are satisfied. These correspond to the
braiding versions of the virtual detour move of type three that involves two virtual crossings and one
real crossing. See \cite{GEN} for the details.

\end{enumerate}
\bigbreak

By extending the methods of \cite{KRCAT}, it is not hard to see that {\em a bi-oriented quantum algebra will always give rise to invariants
of virtual links up to the type one moves (framing and virtual framing).}
\bigbreak

\noindent
In the case of the generalized Alexander polynomial, the state model $Z(K)$ translates directly into a specific example of a bi-oriented balanced quantum
algebra $(A, \rho, \gamma, T).$  The main point about this bi-oriented quantum algebra is that the operator $\gamma$ for the virtual
crossing is {\em not} the identity operator; this non-triviality is crucial to the structure of the invariant.
We will investigate bi-oriented quantum  algebras and other examples of virtual invariants derived from them.
\bigbreak

We have taken a path to explain not only the evolution of a theory of invariants
of virtual knots and links, but also (in this subsection) a description of our oriented quantum algebra formulation of
the whole class of quantum link invariants. Returning to the case of the original Jones polynomial, we want to understand
its capabilities in terms of the oriented quantum algebra that generates the invariant.
\bigbreak

\subsection{Invariants of Three-Manifolds}
As is well-known, invariants of three-manifolds can be formulated in terms of Hopf algebras and quantum algebras and spin
recoupling networks. In formulating such invariants it is useful to represent the three-manifold via surgery on a framed link.
Two framed links that are equivalent in the Kirby calculus of links represent the same three-manifold and conversely.
To obtain invariants of three-manifolds one constructs invariants of framed links that are also invariant under the Kirby moves
(handle sliding, blowing up and blowing down).

A classical three-manifold is mathematically
the same as a Kirby equivalence class of a framed link. The fundamental group of the three-manifold associated with a link is
equal to the fundamental group of the complement of the link modulo the subgroup generated by the framing longtudes for the link.
We refer to the fundmental group of the three-manifold as the {\em three-manifold group}.
If there is a counterexample to the
classical Poincar\'{e} conjecture, then the counterexample would be represented by surgery on some link $L$ whose three-manifold
group is trivial, but
$L$ is not trivial in Kirby calculus (i.e. it cannot be reduced to nothing).
\bigbreak

Kirby calculus can be generalized to the class of virtual knots and links.
We define a {\em virtual three-manifold} to be a Kirby equivalence class of framed virtual links.
The three-manifold group generalizes via the
combinatorial fundamental group associated to the virtual link (the framing longitudes still exist for virtual links).
The {\em Virtual Poincar\'{e} Conjecture} to virtuals would say that a virtual three-manifold with trivial
fundamental group is trivial in Kirby calculus. However, {\bf The virtual Poincar\'{e} conjecture is false} \cite{DK}.  There exist
virtual links whose three-manifold group is trivial that are nevertheless not Kirby equivalent to nothing. The simplest example
is the virtual knot in Figure 59 .  We detect the non-triviality of the Kirby class of this knot by computing that it has
an
$SU(2)$ Witten invariant that is different from the standard three-sphere.
\bigbreak

\begin{center}
{\tt    \setlength{\unitlength}{0.92pt}
\begin{picture}(97,78)
\thicklines   \put(65,18){\line(-1,0){15}}
              \put(65,73){\line(0,-1){55}}
              \put(94,75){\line(-1,0){29}}
              \put(94,55){\line(0,1){20}}
              \put(74,55){\line(1,0){20}}
              \put(22,19){\circle{18}}
              \put(42,43){\line(-1,0){15}}
              \put(4,42){\line(1,0){14}}
              \put(3,20){\line(1,0){36}}
              \put(3,21){\line(0,1){21}}
              \put(21,56){\line(1,0){37}}
              \put(21,3){\line(0,1){53}}
              \put(42,3){\line(-1,0){21}}
              \put(43,43){\line(0,-1){40}}
\end{picture}}
\end{center}

\begin{center} {\bf Figure 59 -- A counterexample to the Poincar\'{e} Conjecture for Virtual Three-Manifolds} \end{center}
\bigbreak

This counterexample to the Poincar\'{e} conjecture in the virtual domain shows how a classical counterexample might behave in the
context of Kirby calculus. Virtual knot theory can be used to search for a counterexample to the classical Poincar\'{e} conjecture by
searching for virtual counterexamples that are equivalent in Kirby calculus to classical knots and links. This is a new and
exciting approach to the dark side of the classical Poincar\'{e} conjecture.
\bigbreak

\subsection{ Gauss Diagrams and Vassiliev Invariants}
The reader should recall the notion of a {\em Gauss diagram} for a knot. If $K$ is a knot
diagram, then $G(K),$ its Gauss diagram, is a circle comprising the Gauss code of the knot by arranging the traverse of the diagram from
crossing to crossing along the circle and putting an arrow (in the form of a chord of the circle) between the two appearances of the
crossing. The arrow points from the overcrossing segment to the undercrossing segment in the order of the traverse of the diagram. (Note: Turaev uses another convention, \cite{TURAEV}.)
Each chord is endowed with a sign that is equal to the sign of the corresponding crossing in the knot diagram. At the level of the Gauss diagrams,
a virtual crossing is simply the absence of a chord. That is, if we wish to transcribe a virtual knot diagram to a Gauss diagram, we ignore the
virtual crossings. Reidemeister moves on Gauss diagrams are defined by translation from the corresponding diagrams from planar representation.
Virtual knot theory is precisely the theory of {\em arbitrary Gauss diagrams}, up to the Gauss diagram Reidemeister moves. Note that an arbitrary
Gauss diagram is any pattern of directed, signed chords on an oriented circle. When transcribed back into a planar knot diagram, such a Gauss
diagram may require virtual crossings for its depiction.
\bigbreak

In \cite{GPV} Goussarov, Polyak and Viro initiate a very important program for producing Vassiliev invariants of finite type of virtual and classical
knots. The gist of their program is as follows.They define the notion of a semi-virtual crossing, conceived as a dotted, oriented, signed
chord in a Gauss diagram for a knot. An {\em arrow diagram} is a Gauss diagram all of whose chords are dotted. Let $\cal A$ denote
the collection of all linear combinations of arrow diagrams with integer coefficients. Let $\cal G$ denote the collection of all arbitrary Gauss
diagrams (hence all representatives of virtual knots). Define a mapping $$i:\cal G \longrightarrow \cal A$$ by expanding each chord of
a Gauss diagram $G$ into the sum of replacing the chord by a dotted chord and the removal of that chord. Thus
$$i(G) = \Sigma_{r \in R(G)} G^{r}$$ where $R(G)$ denotes all ways of replacing each chord in $G$ either by a dotted chord, or by nothing; and
$G^{r}$ denotes that particular replacement applied to $G.$
\bigbreak

Now let $\cal P$ denote the quotient of $\cal A$ by the subalgebra generated by the relations in $\cal A$ corresponding to the Reidemeister moves.
Each Reidemeister move is of the form $X = Y$ for certain diagrams, and this translates to the relation $i(X) - i(Y) = 0$ in $\cal P,$ where
$i(X)$ and $i(Y)$ are individually certain linear combinations in $\cal P.$ Let $$I:\cal G \longrightarrow \cal P$$ be the map induced by
$i.$ Then it is a formal fact that $I(G)$ is invariant under each of the Reidemeister moves, and hence that $I(G)$ is an invariant of the
corresponding Gauss diagram (virtual knot) $G.$ The algebra of relations that generate the image of the Reidemeister moves in $\cal P$ is
called the {\em Polyak algebra.}
\bigbreak

So far, we have only desribed a tautological and not a computable invariant. The key to obtaining computable invariants is to truncate.
Let ${\cal P}_{n}$ denote $\cal P$ modulo all arrow diagrams with more than $n$ dotted arrows. Now ${\cal P}_{n}$ is a finitely generated
module over the integers, and the composed map $$I_{n}:{\cal G} \longrightarrow {\cal P}_{n}$$ is also an invariant of virtual knots.
Since we can choose a specific basis for ${\cal P}_{n},$ the invariant $I_{n}$ is in principle computable, and it yields a large collection of
Vassiliev invariants of virtual knots that are of finite type. The paper by Goussarov, Polyak and Viro investigates specific methods for
finding and representing these invariants.  They show that every Vassiliev invariant of finite type for classical knots can be written as
a combinatorial state sum for long knots. They use the virtual knots as an intermediate in the construction.
\bigbreak

By directly constructing Vassiliev invariants of virtual knots from known invariants of virtuals, we can construct invariants that are not of
finite type in the above sense (See \cite{VKT}.) These invariants also deserve further investigation.
\bigbreak

\section{Other Invariants}
The {\em Khovanov Categorification of the Jones polynomial} \cite{DB}
is important for our concerns. This invariant is constructed by
promoting the states in the bracket summation to tensor powers of
a vector space $V$, where a single power of $V$ corresponds to a
single loop in the state. In this way a graded complex is
constructed, whose graded Euler characteristic is equal to the
original Jones polynomial, and the ranks of whose graded homology
groups are themselves invariants of knots. It is now known that
the information in the Khovanov construction exceeds that in the
original Jones polynomial. It is an open problem whether a
Khovanov type construction can generalize to virtual knots in
the general case. The construction for the Khovanov polynomial for virtuals
over ${\bf Z}_{2}$ was proposed in \cite{MAKH}.
Recent work by other authors related to knots in thickened surfaces promises to shed light on this
issue.\bigbreak

One of the more promising directions for
relating Vassiliev invariants to our present concerns is the theory of gropes \cite{CT1,CT2}, where one considers surfaces spanning
a given knot, and then recursively the surfaces spanning curves embedded in the given surface. This hierarchical structure
of curves and surfaces is likely to be a key to understanding the geometric underpinning of the original Jones polynomial.
The same techniques in a new guise could elucidate invariants of virtual knots and links.
\bigbreak

\section{The Bracket Polynomial and The Jones Polynomial}
It is an open problem whether there exist classical knots (single component loops) that are knotted and yet have unit Jones polynomial.
In other words, it is an open problem whether the Jones polynnomial can detect all knots. 
There do exist families of links whose linkedness is undectable by the Jones polynomial \cite{Morwen,EKT}. It is the purpose of this 
section of the paper to give a summary of some of the information that is known in this arena. We begin with a sketch of ways to calculate
the bracket polynomial model of the Jones polynomial, and then discuss how to construct classical links that are undetectable by the 
Jones polynomial.
\bigbreak

The formula
for the bracket model of the Jones polynomial \cite{state} can be indicated as follows: The letter chi, \mbox{\large $\chi$}, denotes a
crossing in a link diagram. The barred letter denotes the mirror image of this first crossing. A crossing in a diagram for the knot or link
is expanded into two possible states by either smoothing (reconnecting) the crossing horizontally,
\mbox{\large
$\asymp$}, or vertically $><$. Any closed loop (without
crossings) in the plane has value $\delta = -A^{2} - A^{-2}.$    
$$\mbox{\large $\chi$} = A \mbox{\large $\asymp$} + A^{-1} ><$$
$$\overline{\mbox{\large $\chi$}} = A^{-1} \mbox{\large $\asymp$} + A ><.$$
\noindent One useful consequence of these formulas is the following {\em switching formula}
$$A\mbox{\large $\chi$} - A^{-1} \overline{\mbox{\large $\chi$}} = (A^{2} - A^{-2})\mbox{\large $\asymp$}.$$ Note that 
in these conventions the $A$-smoothing of $\mbox{\large $\chi$}$ is $\mbox{\large $\asymp$},$ while the $A$-smoothing of
$\overline{\mbox{\large $\chi$}}$ is $><.$ Properly interpreted, the switching formula above says that you can switch a crossing and 
smooth it either way and obtain a three diagram relation. This is useful since some computations will simplify quite quickly with the 
proper choices of switching and smoothing. Remember that it is necessary to keep track of the diagrams up to regular isotopy (the 
equivalence relation generated by the second and third Reidemeister moves). Here is an example. View Figure 60.

\noindent You see in Figure 60, a trefoil diagram $K$, an unknot diagram $U$ and another unknot diagram $U'.$  Applying the switching formula,
we have $$A^{-1} K - A U = (A^{-2} - A^{2}) U'$$ and 
$U= -A^{3}$ and $U'=(-A^{-3})^2 = A^{-6}.$  Thus $$A^{-1} K - A(-A^{3}) = (A^{-2} - A^{2}) A^{-6}.$$ Hence
$$A^{-1} K = -A^4 + A^{-8} - A^{-4}.$$  Thus $$K = -A^{5} - A^{-3} + A^{-7}.$$ This is the bracket polynomial of the trefoil diagram $K.$
We have used to same symbol for the diagram and for its polynomial.
\bigbreak

\centerline{\includegraphics[scale=1.0]{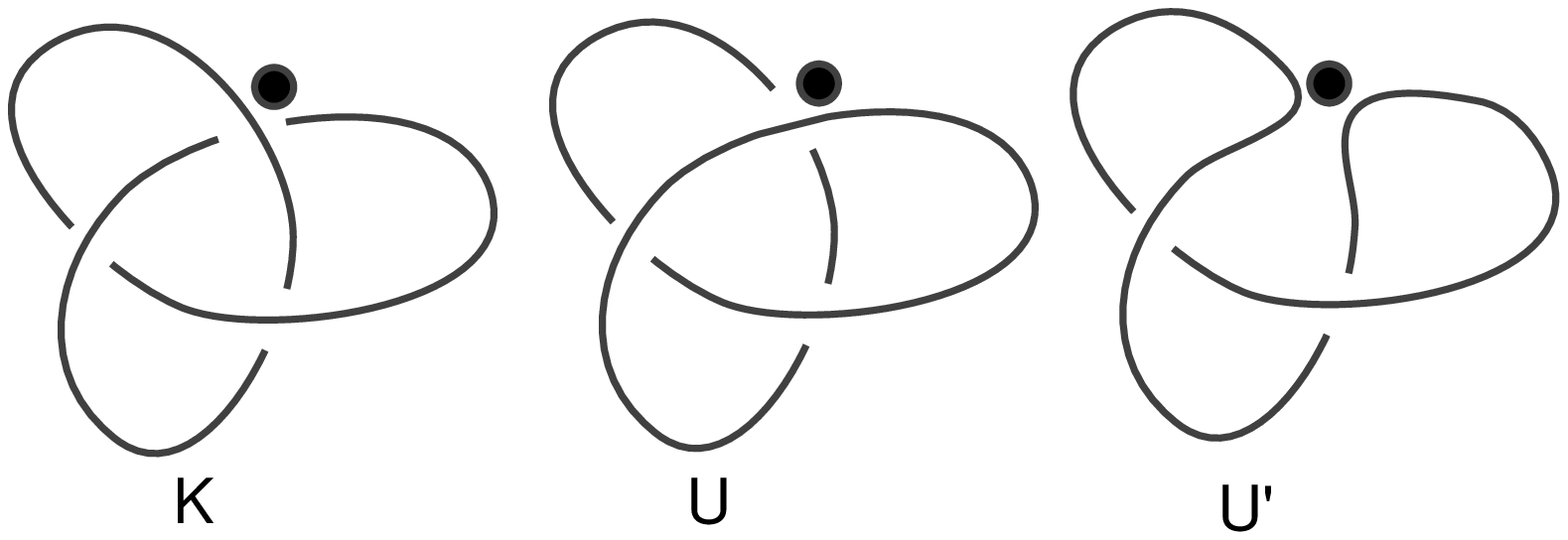}}
\vvvv

\begin{center} {\bf Figure 60 -- Trefoil and Two Relatives} \end{center}
\bigbreak

\centerline{\includegraphics[scale=1.0]{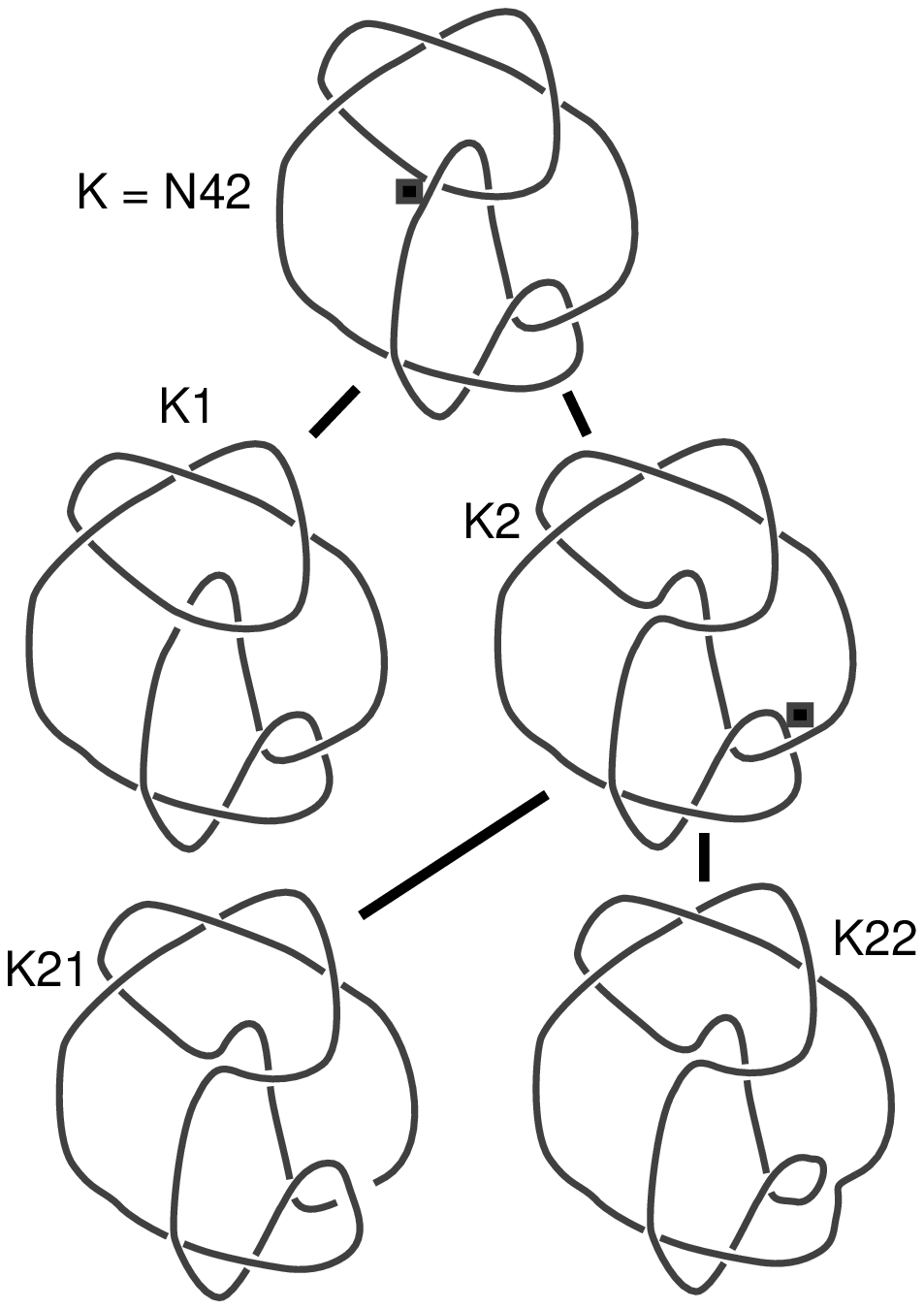}}
\vvvv

\begin{center} {\bf Figure 61 -- Skein Tree for $9_{42}$} \end{center}
\bigbreak

\noindent Since the trefoil diagram $K$ has writhe $w(K) = 3,$ we have the normalized polynomial 
$$f_{K}(A) = (-A^{3})^{-3}<K> = -A^{-9}(-A^{5} - A^{-3} + A^{-7}) = A^{-4} + A^{-12} - A^{-16}.$$ The asymmetry of this polynomial under the 
interchange of $A$ and $A^{-1}$ proves that the trefoil knot is not ambient isotopic to its mirror image.
\bigbreak

In Figure 61 you see the knot $K = N_{42} = 9_{42}$ (the latter being its standard name in the knot tables) and a skein tree for it via
switching and smoothing.  In Figure 62 we show simplified (via regular isotopy) representatives for the end diagrams in the skein tree.
\bigbreak

\centerline{\includegraphics[scale=1.0]{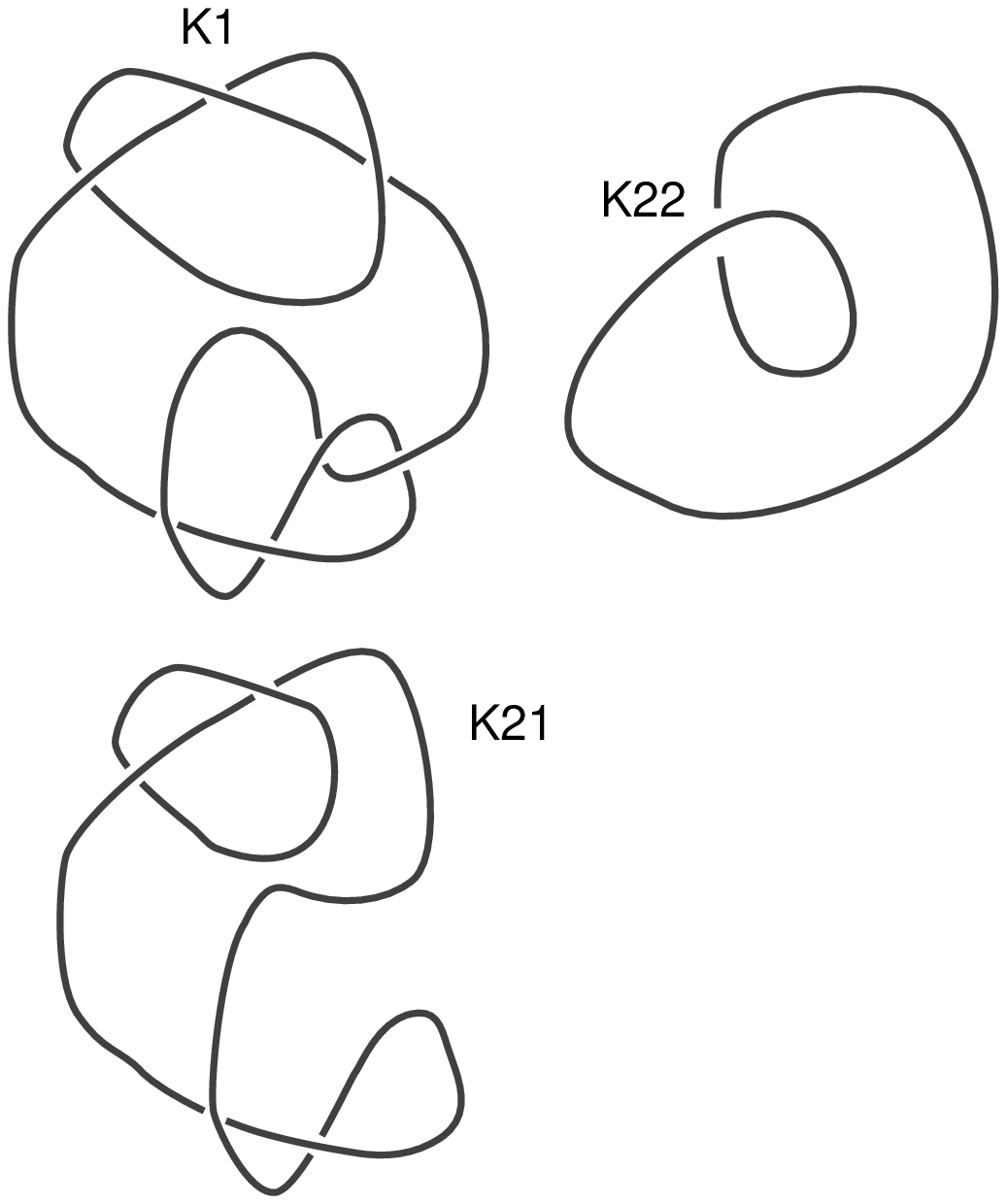}}
\vvvv

\begin{center} {\bf Figure 62 -- Regular Isotopy Versions of Bottom of Skein Tree for $9_{42}$} \end{center}
\bigbreak

\noindent It follows from the switching formula that for $K = 9_{42},$
$$ A^{-1} K - A K_{1} = (A^{-2} - A^{2})K_{2}$$
$$A K_{2} - A^{-1}K_{21} = (A^{2} - A^{-2})K_{22}$$ and that $K_{1}$ is a connected sum of a right-handed trefoil diagram and a figure eight
knot diagram, while $K_{21}$ is a Hopf link (simple link of linking number one) with extra writhe of $-2$ while $K_{22}$ is an unknot with 
writhe of $1.$ These formulas combine to give $$<K_{1}> = -A^{-9} + A^{-5} - A^{-1} + A^{3} - A^{7} + A^{11} - A^{15}.$$ Since $K$ has 
writhe one, we get $$f_{K} = A^{-12} - A^{-8} + A^{-4} - 1 + A^{4} - A^{8} + A^{12}.$$ This shows that the normalized bracket polynomial 
does not distinguish $9_{42}$ from its mirror image. This knot is, in fact chiral (inequivalent to its mirror image), a fact that can be 
verified by other means. The knot $9_{42}$ is the first chiral knot whose chirality is undetected by the Jones polynomial.

\noindent {\bf Remark.} In writing computer programs for calculating the bracket polynomial it is useful to use the coding method 
illustrated in Figure 63. In this method each edge from one classical crossing to another has a label. Virtual crossings do not
appear in the code. Each classical crossing has a four letter code in the form $[abcd]$ connoting a clockwise encirclement of the 
crossing, starting at an overcrossing edge. 
\bigbreak

\centerline{\includegraphics[scale=1.0]{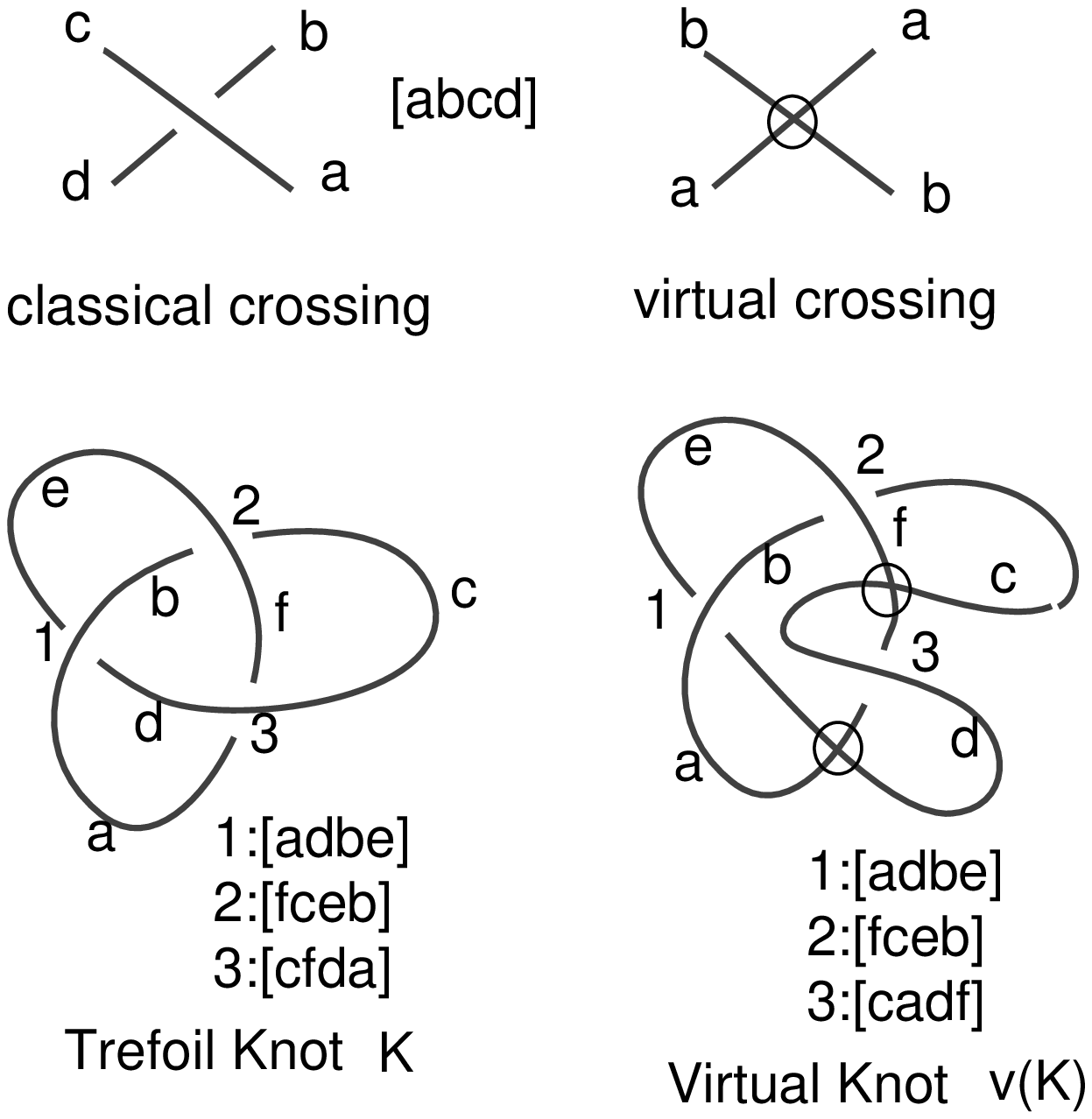}}
\vvvv

\begin{center} {\bf Figure 63 -- Coding a Link Diagram} \end{center}
\bigbreak

\subsection{Thistlethwaite's Example}
View Figure 64. Here we have a version of a link $L$ discovered by Morwen Thistlethwaite \cite{Morwen} in
December 2000. We discuss some theory behind this link in the next subsection. It is a link that is linked but whose linking is not 
detectable by the Jones polynomial. One can verify such properties by using a computer program, or by the algebraic techniques described
below.  

\centerline{\includegraphics[scale=1.0]{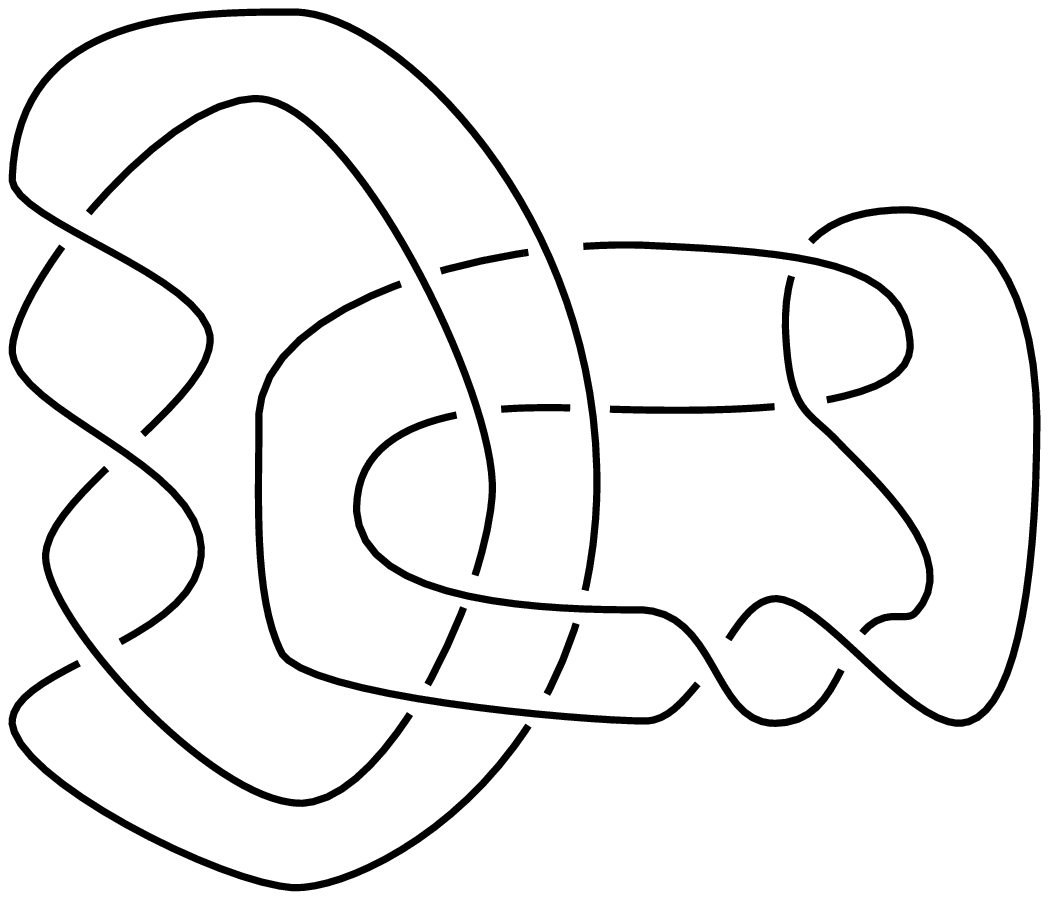}}
\vvvv

\begin{center} {\bf Figure 64 -- Thistlethwaite's Link} \end{center}
\bigbreak

\subsection{Present Status of Links Not Detectable by the Jones Polynomial}
In this section we give a quick review of the status of our work \cite{EKT}
producing infinite families of distinct links all evaluating as unlinks by the 
Jones polynomial.  
\bigbreak

A tangle (2-tangle) consists in an embedding of two arcs in a three-ball (and possibly some circles
embedded in the interior of the three-ball) such that the endpoints of the arcs are
on the boundary of the three-ball. One usually depicts the arcs as crossing the boundary transversely so that the
tangle is seen as the embedding in the three-ball augmented by four segments emanating from the ball, each from 
the intersection of the arcs with the boundary. These four segments are the {\em exterior edges} of the 
tangle, and are used for operations that form new tangles and new knots and links from given tangles. 
Two tangles in a given three-ball are said to be {\em topologically equivalent} if there is an ambient isotopy
from one to the other in the given three-ball, fixing the intersections of the tangles with the boundary.
\bigbreak

It is customary to illustrate tangles with a diagram that consists in a box  (within which are the arcs of the
tangle) and with the exterior edges emanating from the box  in the NorthWest (NW), NorthEast (NE), SouthWest (SW)
and SouthEast (SE) directions. Given tangles $T$ and $S$, one defines the {\em sum}, denoted $T+S$ by placing
the diagram for $S$ to the right of the diagram for $T$ and attaching the NE edge of $T$ to the NW edge of $S$, and
the SE edge of $T$ to the SW edge of $S$. The resulting tangle $T+S$ has exterior edges correponding to the 
NW and SW edges of $T$ and the NE and SE edges of $S$. 
There are two ways to create links associated to a tangle $T.$ The {\em numerator} $T^{N}$ is obtained, by attaching
the (top) NW and NE edges of $T$ together and attaching the (bottom) SW and SE edges together. 
The denominator $T^{D}$ is obtained,
by attaching the (left side) NW and SW edges together and attaching the (right side) NE and SE edges together.  
We denote by $[0]$ the tangle with only unknotted arcs (no embedded circles) with one arc connecting, within
the three-ball, the (top points) NW intersection point with the NE intersection point, and the other arc connnecting
the (bottom points) SW intersection point with the SE intersection point. 
A ninety degree turn of the tangle $[0]$ produces
the tangle $[\infty]$ with connections between NW and SW and between NE and SE. One then can prove the basic 
formula for any tangle $T$
$$<T> = \alpha_{T} <[0]> + \beta_{T}<[\infty]>$$
\noindent where $\alpha_{T}$ and $\beta_{T}$ are well-defined polynomial invariants (of regular isotopy) of the
tangle $T.$  From this formula one can deduce that 
$$<T^{N}> = \alpha_{T} d + \beta_{T}$$
\noindent and
$$<T^{D}> = \alpha_{T} + \beta_{T} d.$$
\bigbreak

We define the {\em bracket vector} of $T$ to be the ordered pair $(\alpha_{T}, \beta_{T})$ and denote it by
$br(T)$, viewing it as a column vector so that $br(T)^{t} = (\alpha_{T}, \beta_{T})$ where $v^{t}$ denotes the 
transpose of the vector $v.$ With this notation the two formulas above for the evaluation for numerator and 
denominator of a tangle become the single matrix equation

$$\left[
 \begin{array}{c}
      <T^{N}> \\
      <T^{D}>
 \end{array}
 \right] =  \left[
 \begin{array}{cc}
      d & 1 \\
      1 & d
 \end{array}
 \right] br(T).$$
\bigbreak
 
 We then use this formalism to express the bracket polynomial for our examples. The class of examples that 
 we considered are each denoted by $H(T,U)$ where $T$ and $U$ are each tangles and $H(T,U)$ is a satellite
 of the Hopf link that conforms to the pattern shown in Figure 65, formed by clasping together the numerators of the
 tangles $T$ and $U.$ Our method is based on a transformation $H(T,U) \longrightarrow H(T,U)^{\omega}$,
 whereby the tangles $T$ and $U$ are cut out and reglued by certain specific homeomorphisms of the tangle 
 boundaries. This transformation can be specified by a modification described by a specific rational tangle
 and its mirror image. Like mutation, the transformation $\omega$ preserves the bracket polynomial. However,
 it is more effective than mutation in generating examples, as a trivial link can be transformed to a prime
 link, and repeated application yields an infinite sequence of inequivalent links.
\bigbreak

\centerline{\includegraphics[scale=1.0]{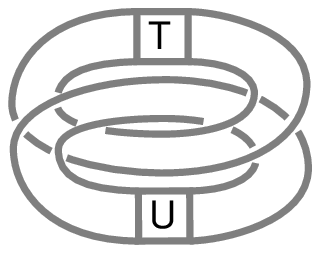}}
\vvvv

\begin{center} {\bf Figure 65 -- Hopf Link Satellite H(T,U)} \end{center}
\bigbreak

\centerline{\includegraphics[scale=1.0]{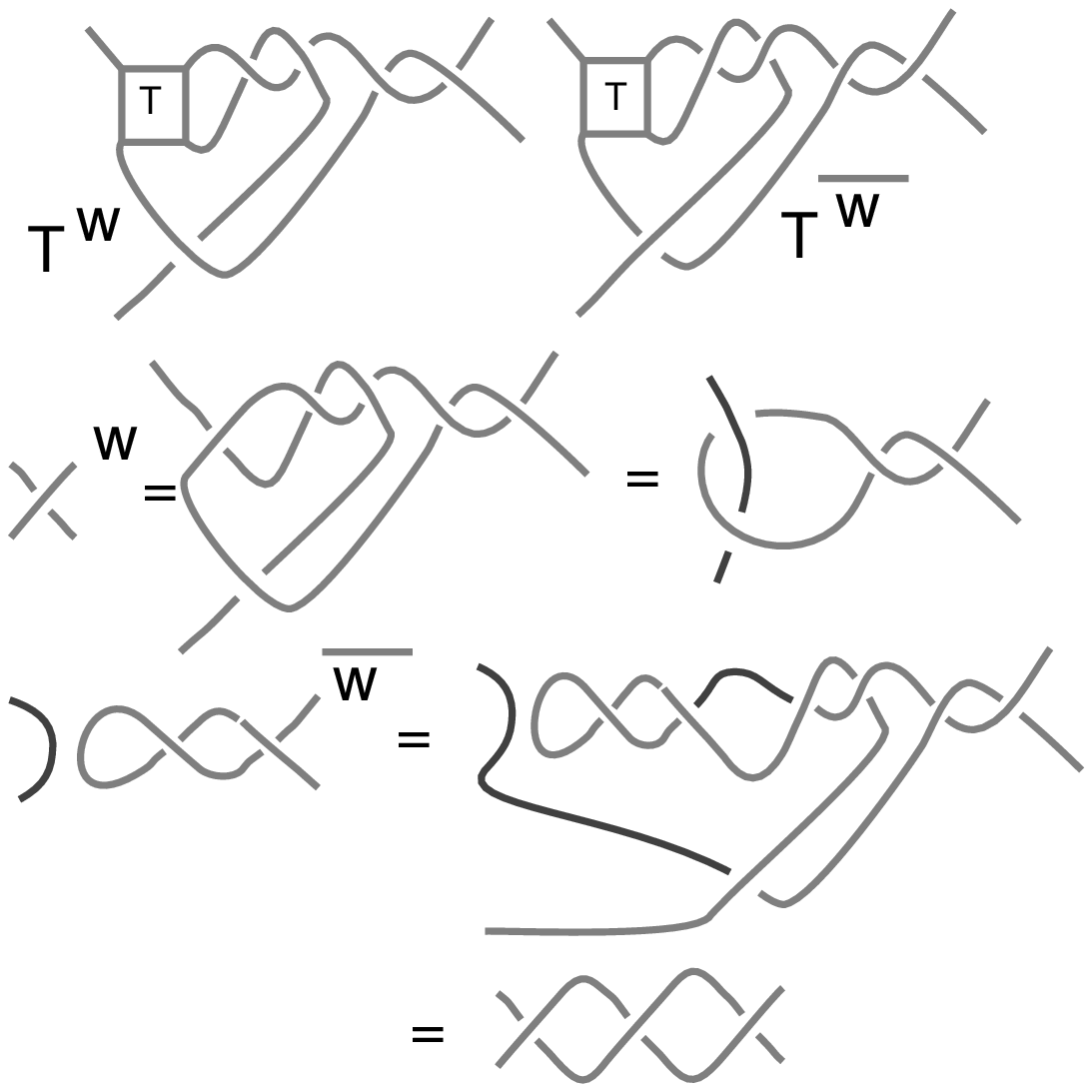}}
\vvvv

\begin{center} {\bf Figure 66 -- The Omega Operations} \end{center}
\bigbreak

\centerline{\includegraphics[scale=1.0]{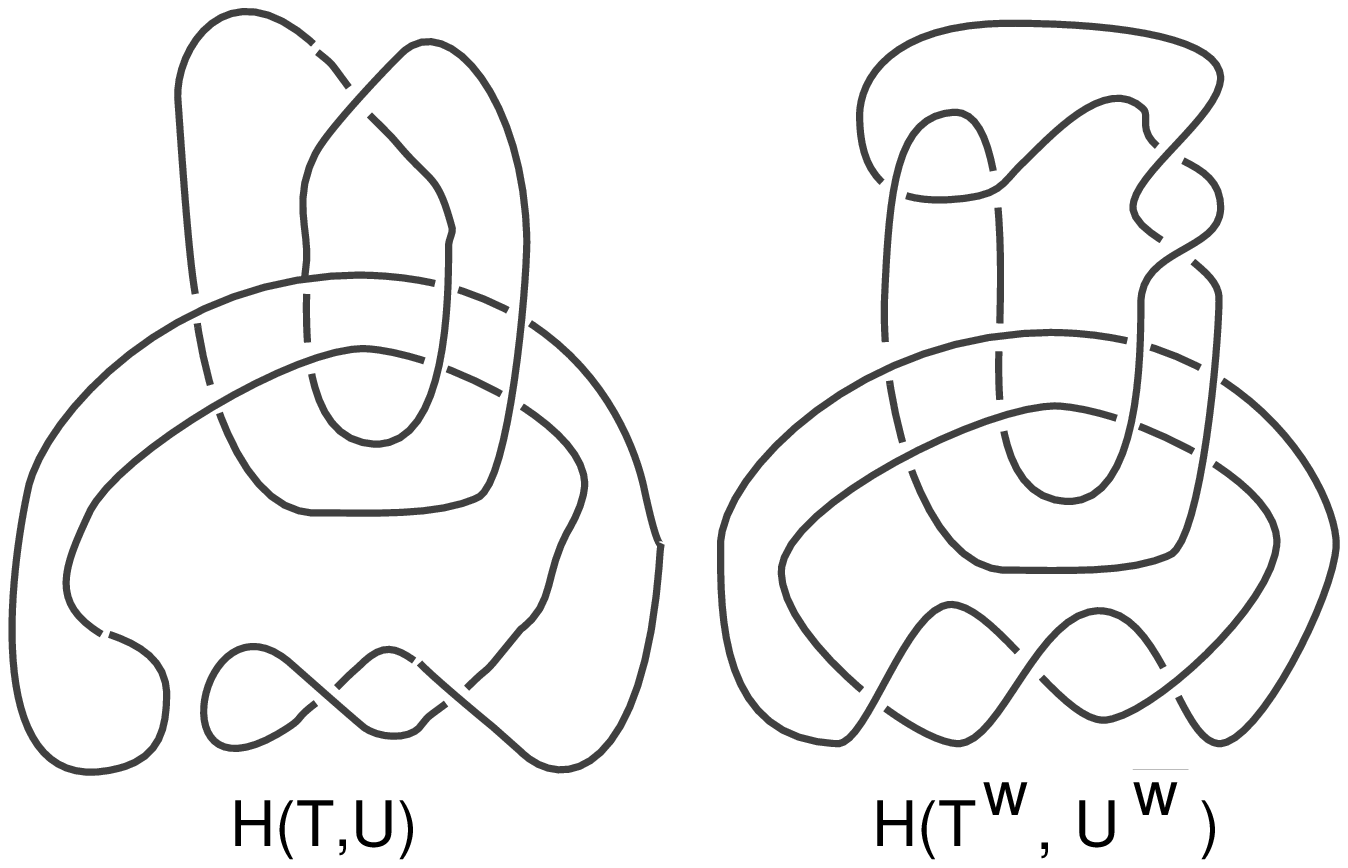}}
\vvvv

\begin{center} {\bf Figure 67 -- Applying Omega Operations to an Unlink} \end{center}
\vvvv

 Specifically, the transformation $H(T,U)^{\omega}$ is given by the formula 
 $$H(T,U)^{\omega} = H(T^{\omega}, U^{\bar{\omega}})$$ 
 \noindent where the tangle operations $T^{\omega}$ and  $U^{\bar{\omega}})$ are as shown in Figure 66.
 By direct calculation, there is a matrix $M$ such that 
 $$<H(T,U)> = br(T)^{t} M br(U)$$
 \noindent and there is a matrix $\Omega$ such that 
 $$br(T^{\omega}) = \Omega br(T)$$
 \noindent and 
 $$br(T^{\bar{\omega}}) = \Omega^{-1} br(T).$$
 \noindent One verifies the identity
 $$\Omega^{t} M \Omega^{-1} = M$$
 \noindent from which it follows that $<H(T,U)>^{\omega} = <H(T,U)>.$
 This completes the sketch of our method for obtaining links that whose linking cannot be seen by the Jones polynomial.
Note that the link constructed as $H(T^{\omega}, U^{\bar{\omega}})$ in Figure 67 has the same Jones polynomial as an unlink of
two components. This shows how the first example found by Thistlethwaite fits into our construction.
\bigbreak

\subsection{Switching a Crossing}
If in Figure 67, we start with $T$ replaced by $Flip(T)$, switching the crossing, the resulting link $L = H(Flip(T)^{\omega}, U^{\bar{\omega}})$ 
will still have Jones polynomial the same as the unlink, but the link $L$ will be distinct from the link $H(T^{\omega}, U^{\bar{\omega}})$
of Figure 67.  We illustrate this process in Figure 68.

\centerline{\includegraphics[scale=0.95]{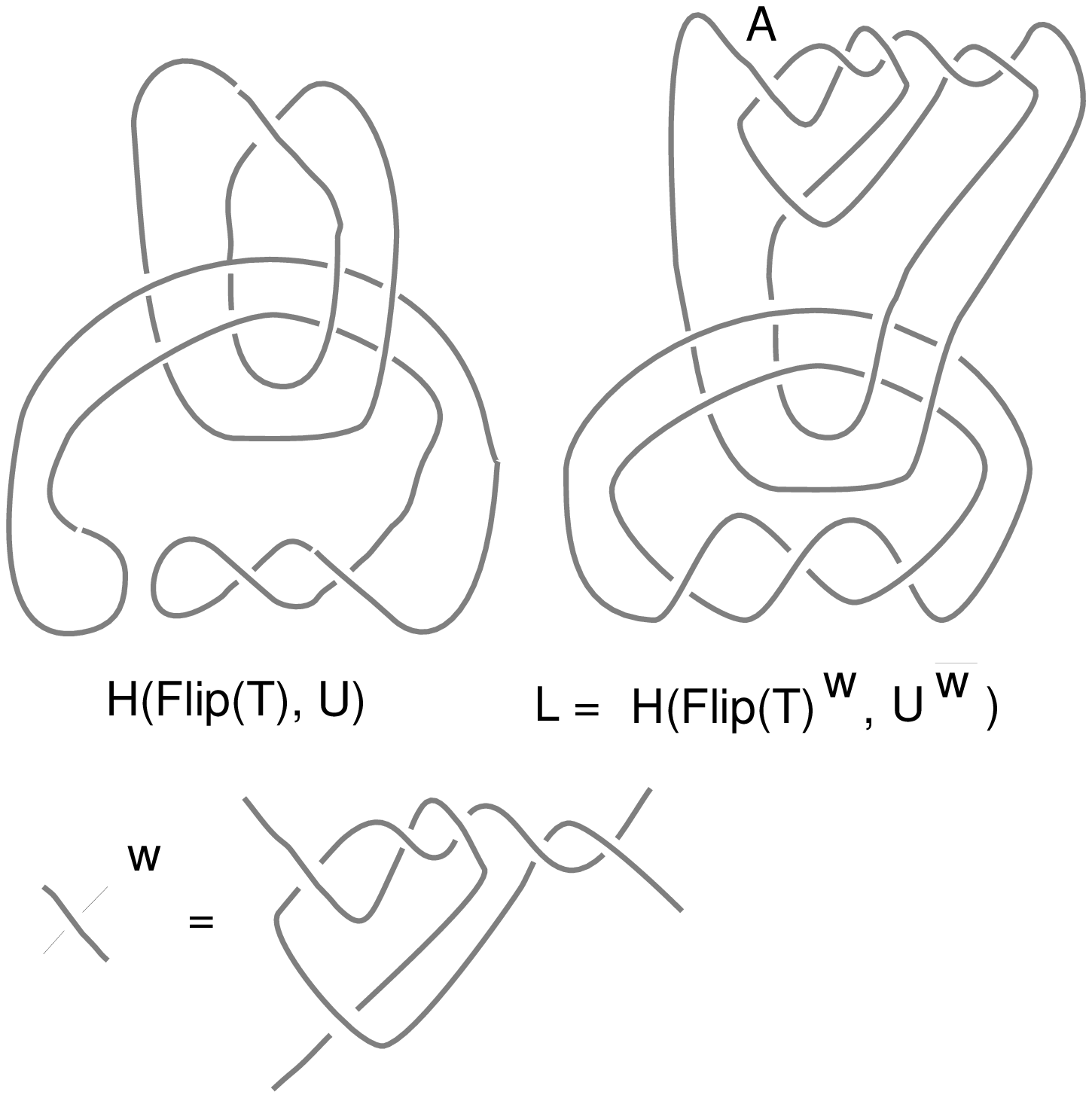}}
\vvvv

\begin{center} {\bf Figure 68 -- Applying Omega Operations to an Unlink with Flipped Crossing} \end{center}
\vvvv
 
The link $L$ has the remarkable property that both it and the link obtained from it by flipping the crossing labeled $A$ in Figure 68 have 
Jones polynomial equal to the Jones polynomial of the unlink of two components. (We thank Alexander Stoimenow for pointing out the possibility
of this sort of construction.)
\bigbreak

Now let's think about a link $L$ with the property that it has the same Jones polynomial as a link $L'$ obtained from $L$ by switching a single
crossing. We can isolate the rest of the link that is not this crossing into a tangle $S$ so that (without any loss of 
generality) $L = N(S + [1])$ and $L' = N(S + [-1]).$ Lets assume that orientation assignment to $L$ and $L'$ is as shown in Figure 69.
\bigbreak

\centerline{\includegraphics[scale=1.0]{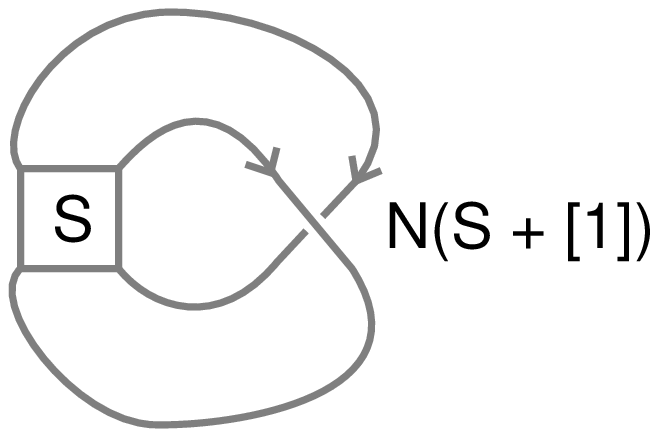}}
\vvvv

\begin{center} {\bf Figure 69 -- N(S + [1])} \end{center}
\vvvv

Since we are told that $L$ and $L'$ have the same Jones polynomial, it follows that 
$<L> = -A^{-3} \kappa$ and $<L'> = -A^{3} \kappa$ for some non-zero Laurent polynomial $\kappa.$
Now suppose that $<S> = \alpha <[0]> + \beta <[\infty]>.$ Then 
$$<L> = \alpha (-A^3) + \beta (-A^{-3})$$ and 
$$<L'> = \alpha (-A^{-3}) + \beta (-A^3).$$ From this it follows that 
$$ \kappa = \alpha A^6 + \beta$$ and 
$$\kappa = \alpha A^{-6} + \beta.$$ Thus
$$\alpha = 0$$ and $$\beta = \kappa.$$  We have shown that 
$$<S> = \kappa <[\infty]>.$$ This means that we can, by using the example described above, produce a tangle $S$ that is not splittable
and yet has the above property of having one of its bracket coefficients equal to zero. The example is shown in Figure 70.

\centerline{\includegraphics[scale=1.0]{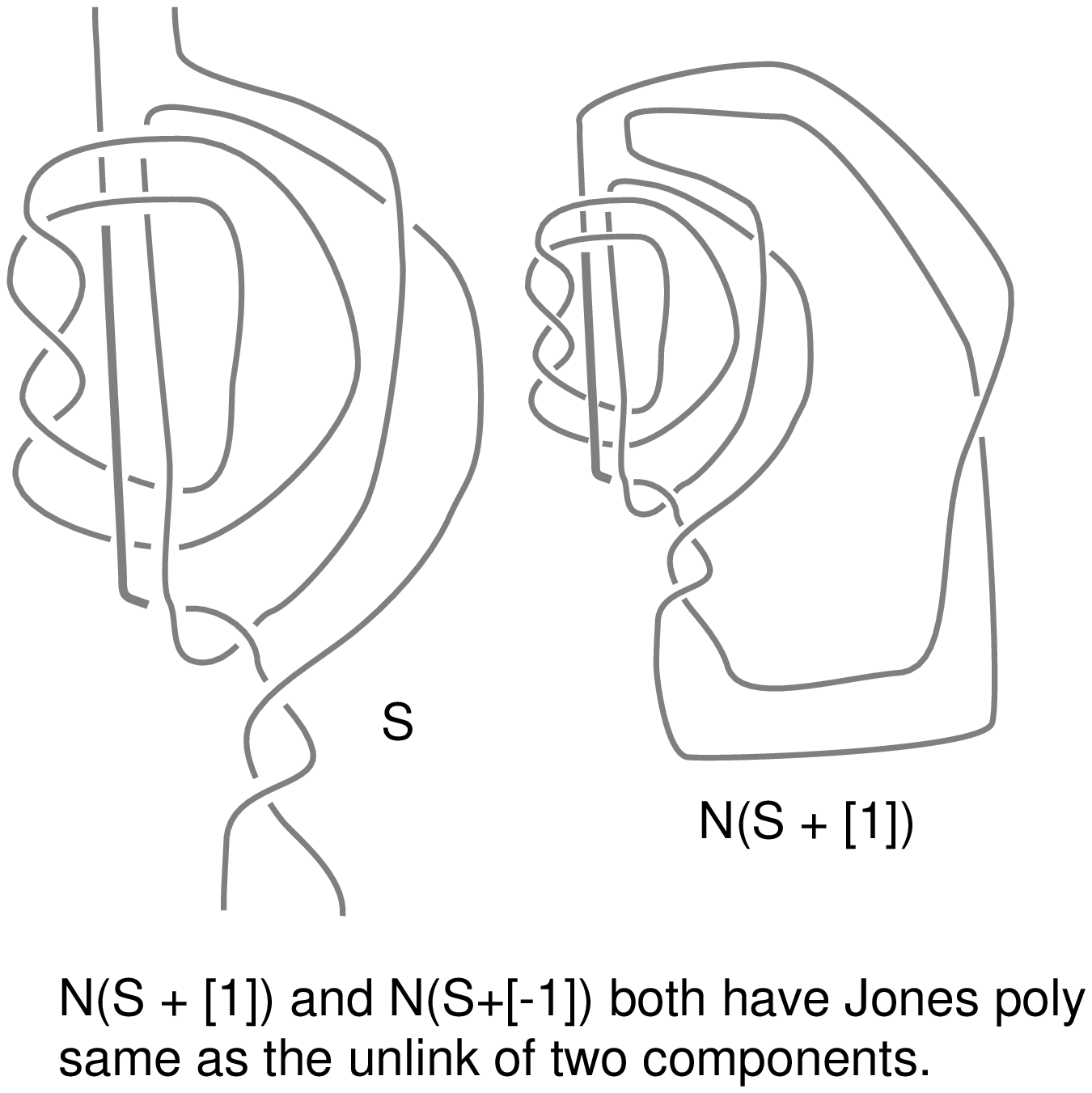}}
\vvvv

\begin{center} {\bf Figure 70 -- The Tangle S} \end{center}
\vvvv\vvv

Finally, in Figure 71 we show $L$ and the link $v(L)$ obtained by virtualizing the crossing corresponding to $[1]$ in the decomposition
$L = N(S + [1]).$  The virtualized link $v(L)$ has the property that it also has Jones polynonial the same as an unlink of two components.
We wish to prove that $v(L)$ is not isotopic to a classical link. The example has been designed so that surface bracket techniques will be
difficult to apply.
\newpage

\centerline{\includegraphics[scale=1.0]{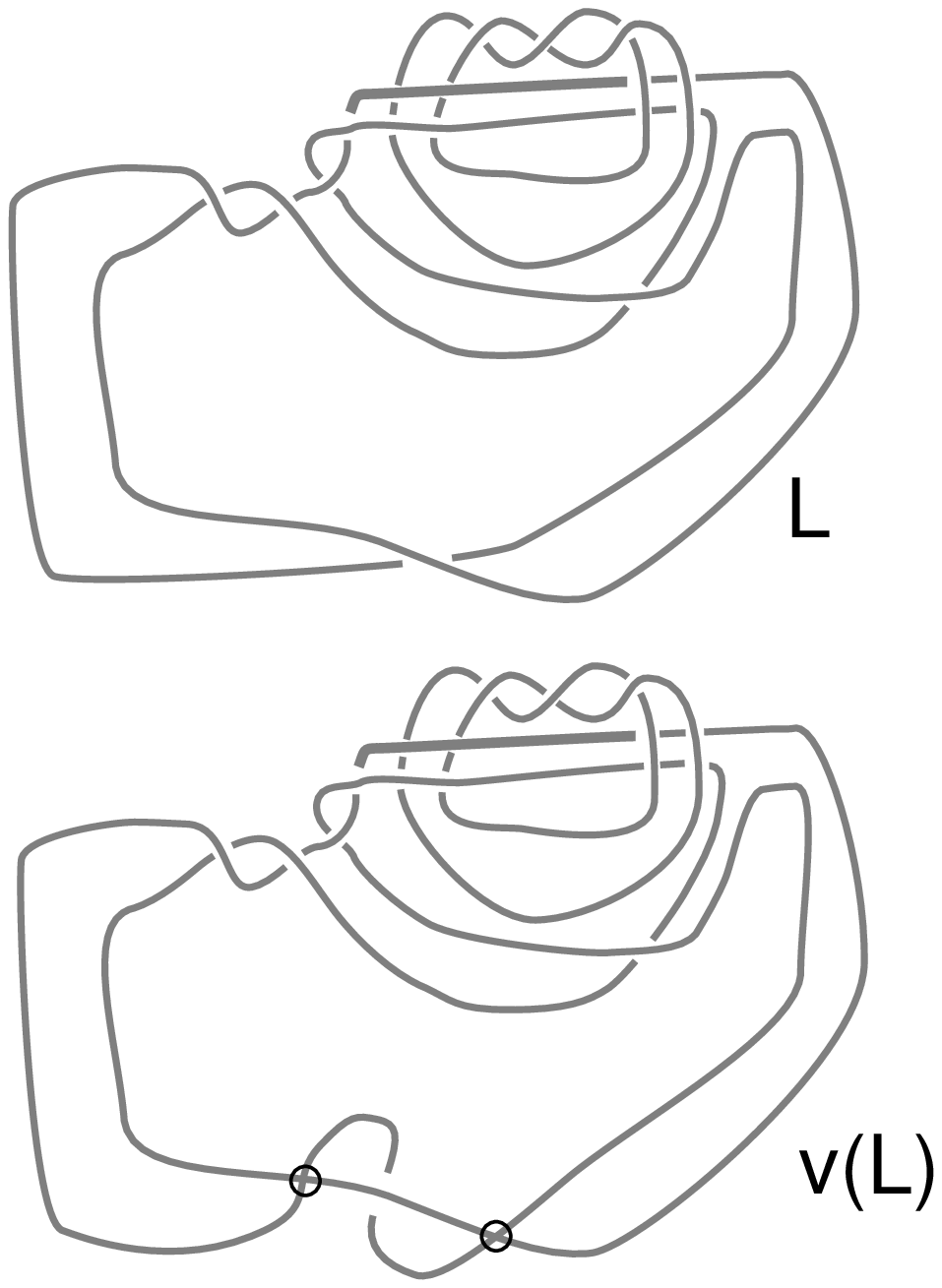}}
\vvvv

\begin{center} {\bf Figure 71 -- The Virtual Link v(L)} \end{center}

\end{document}